\documentclass[a4paper]{amsart}
\usepackage{graphicx,amsmath,amssymb,cases} 
\usepackage{color}
\usepackage{tikz}
\usepackage{tikz-cd}
\usetikzlibrary{cd,positioning}
\usepackage{array,booktabs,longtable}
\usepackage{verbatim}
\usepackage[dvipsnames]{xcolor}
\usepackage{ytableau}
\usepackage[all]{xy}
\SelectTips{cm}{12}
\objectmargin={1pt}
\newsavebox{\AFGDTableBox}
\newsavebox{\AFGQTableBox}
\usepackage[
  a4paper,
  left=30mm,
  right=30mm,
  top=25mm,
  bottom=25mm
]{geometry}
 
\usepackage{geometry}

\usepackage{graphicx}
\usepackage{booktabs}
\usepackage{pdflscape}

\usepackage[mathscr]{euscript}

\title[Equivariant $K$-homology of Affine Grassmannians]{Equivariant $K$-homology of Affine Grassmannians of orthogonal groups}
\author[K. Brahma]{Koushik Brahma}
\address{Waseda University, Faculty of Science and Engineering,
3-4-1 Okubo, Shinjuku-ku
Tokyo 169-8555, Japan}
\email{koushikbrahma95@gmail.com}

\author[T. Ikeda]{Takeshi Ikeda}
\address{Waseda University, Faculty of Science and Engineering,
3-4-1 Okubo, Shinjuku-ku
Tokyo 169-8555, Japan}
\email{gakuikeda@waseda.jp}

\author[M. Shimozono]{Mark Shimozono}
\address{Department of Mathematics, Virginia Polytechnic Institute and State University, Blacksburg, VA 24061-0123 USA}
\email{mshimo@math.vt.edu}

\date{\today}

\allowdisplaybreaks[4] 

\usepackage[
  colorinlistoftodos,
  textwidth=25mm
]{todonotes}

\usepackage[colorinlistoftodos]{todonotes}

\newcommand{\Gr}{\mathrm{Gr}}
\newcommand{\af}{\mathrm{af}}
\DeclareMathOperator{\Ker}{Ker}
\newcommand{\eps}{\varepsilon}
\newcommand{\G}{\mathfrak{G}}
\newcommand{\Q}{\mathbb{Q}}
\newcommand{\Sc}{{S_\Q}}
\newcommand{\sco}{\mathrm{sc}}
\newcommand{\gGtil}{g\widetilde{\Gamma}}
\newcommand{\gGhat}{g\widehat{\Gamma}}

\newcommand{\Z}{\mathbb{Z}}

\newcommand{\DD}{\widehat{D}_{w}}
\newcommand{\SP}{\mathscr{SP}}
\newcommand{\half}{\frac{1}{2}}
\newcommand{\hL}{\hat{\Lambda}}
\DeclareMathOperator{\Spin}{Spin}

 \newcommand{\htK}[1]{\widehat{\Gamma}_{#1}}
  \newcommand{\htL}[1]{\widetilde{\Gamma}_{#1}}

\newcommand{\sR}{\mathscr{R}}
\newcommand{\RT}{R(T)}
\newcommand{\Winf}{W_{\mathrm{D}_\infty}}
\newcommand{\Winfgr}{\Winf^0}
\newcommand{\Waf}{W_\af}
\newcommand{\Wafex}{\widetilde{W}_\af}

\newcommand{\Wafexgr}{\Wafex^0}
\newcommand{\Wafgr}{W_\af^0}

\newcommand{\KGd}{\mathbb{K}_G^\Delta}
\newcommand{\KG}{\mathbb{K}_G}
\newcommand{\tKGd}{\widetilde{\mathbb{K}}_G^\Delta}
\newcommand{\tKG}{\widetilde{\mathbb{K}}_G}
\newcommand{\LG}{\mathbb{L}_G}
\newcommand{\tLG}{\widetilde{\mathbb{L}}_G}

\newcommand{\A}[1]{\mathrm{A}_{#1}}
\newcommand{\B}[1]{\mathrm{B}_{#1}}
\DeclareRobustCommand{\D}[1]{\mathrm{D}_{#1}}%

\newcommand{\gq}[2]{gq^{(#2)}_{#1}}
\newcommand{\gd}[2]{{gd}^{(#2)}_{#1}}
\newcommand{\tgq}[2]{\widetilde{gq}^{(#2)}_{#1}}
\newcommand{\tgd}[2]{\widetilde{gd}^{(#2)}_{#1}}

\newcommand{\tf}[1]{\eta_{#1}}

\newcommand{\tgx}[2]{\widetilde{gx}^{(#2)}_{#1}}
\newcommand{\gx}[2]{{gx}^{(#2)}_{#1}}

\newcommand{\Iwa}{\mathscr{I}}

\newcommand{\WBgr}[1]{W_{\mathrm{B}_{#1}^{(1)}}^0}
\newcommand{\WDgr}[1]{W_{\mathrm{D}_{#1}^{(1)}}^0}
\newcommand{\WBexgr}[1]{\widetilde{W}_{\mathrm{B}_{#1}^{(1)}}^0}
\newcommand{\WDexgr}[1]{\widetilde{W}_{\mathrm{D}_{#1}^{(1)}}^0}

\newcommand{\WB}[1]{W_{\mathrm{B}_{#1}}}
\newcommand{\WD}[1]{W_{\mathrm{D}_{#1}}}

\newcommand{\Aut}{\mathrm{Aut}}
\newcommand{\auto}{\tau_{01}}
\newcommand{\autn}{\tau_{n-1,n}}
\newcommand{\autDrev}[1]{\sigma^{\mathrm{D}_{#1}^{(1)}}}
\newcommand{\ba}[1]{\overline{#1}}
\newcommand{\C}{\mathbb{C}}

\newcommand{\cl}{\mathrm{cl}}
\newcommand{\cO}[1]{\mathscr{O}_{#1}}

\newcommand{\Fl}{\mathrm{Fl}}

\newcommand{\Fund}{\Sigma}

\newcommand{\Gsc}{G_\sco}
\newcommand{\Gad}{G}
\newcommand{\Tsc}{T_{\sco}}
\newcommand{\Tad}{T}

\newcommand{\R}{\mathscr{R}}

\newcommand{\Rsc}{R(\Tsc)}
\newcommand{\Rad}{R(\Tad)}

\newcommand{\id}{\mathrm{id}}

\newcommand{\la}{\lambda}
\newcommand{\pair}[2]{\langle #1\,,\,#2\rangle}

\newcommand{\pt}{\mathrm{pt}}

\newcommand{\Wx}{\widetilde{W}_\af}
\newcommand{\Wxgr}{\Wx^0}
\newcommand{\Xinf}{X_\infty}
\newcommand{\Xvinf}{\Xinf^\vee}

\newcommand{\Y}{\mathscr{P}}

\newcommand{\cOO}[1]{[\mathscr{O}^{#1}]}

\newtheorem*{thmA}{Main Theorem}

\newtheorem{thm}{Theorem}[section]
\newtheorem{lem}[thm]{Lemma}
\newtheorem{prop}[thm]{Proposition}
\newtheorem{cor}[thm]{Corollary}
\newtheorem{conj}[thm]{Conjecture}

\theoremstyle{remark}
\newtheorem{example}[thm]{Example}
\newtheorem{defn}[thm]{Definition}
\newtheorem{remark}[thm]{Remark}

\numberwithin{equation}{section}

\begin{document}


\begin{abstract}
  We construct explicit symmetric-function realizations of the torus-equivariant $K$-homology of affine Grassmannians of types B and D. More precisely, we realize the extended $K$-theoretic Peterson algebras as subalgebras of a completed ring of equivariant symmetric functions and identify their Schubert classes in terms of symmetric functions. We introduce two new families of equivariant symmetric functions for each of types B and D and prove that these functions represent the classes of the structure sheaves of Schubert varieties and the ideal sheaves of their boundaries, respectively.
  We establish several fundamental properties of these functions, including translation factorization formulas and explicit expressions for the classes indexed by small affine Grassmannian elements. In type D, we further describe the action of { the affine Dynkin diagram automorphism interchanging the two spin nodes} and show how it intertwines the Schubert bases with a nontrivial transformation of the equivariant parameters. 
\end{abstract}

\maketitle


\section{Introduction}

The purpose of the present paper is to construct explicit symmetric-function realizations of the torus-equivariant \(K\)-homology of affine Grassmannians in the orthogonal types. The main task is to construct the appropriate equivariant symmetric functions in types $\B{}$ and $\D{}$ and
to identify them with the structure-sheaf and boundary ideal-sheaf Schubert bases. Our construction is based on the extended affine Weyl group, the extended $K$-theoretic Peterson algebra and affine Demazure operators. In particular, it applies naturally to the simply connected and adjoint types of the corresponding algebraic groups.

One of the principal motivations for studying \(K_*^T(\Gr_G)\) comes from Peterson's ``quantum equals affine'' philosophy. In equivariant homology, Peterson constructed, after suitable localization, an isomorphism between the equivariant homology of the affine Grassmannian and the equivariant quantum cohomology of the flag variety \(G/B\); see \cite{Pet,LS}. The \(K\)-theoretic analogue was conjectured in \cite{LLMS} and subsequently established in general type by Kato \cite{Kato}; see also \cite{CL}. Thus the Schubert calculus of \(K_*^T(\Gr_G)\) is closely related to the equivariant quantum \(K\)-theory of \(G/B\). A concrete algebraic realization of affine \(K\)-homology is therefore expected to provide useful tools for the study of quantum \(K\)-theory.

In type \(A\), this geometry is closely connected with symmetric functions. The Schubert basis of the homology of the affine Grassmannian is represented by the \(k\)-Schur functions of Lapointe, Lascoux, and Morse \cite{LLM}, as proved by Lam \cite{Lam}. Their \(K\)-theoretic analogues were introduced by Lam, Schilling and Shimozono in \cite{LSS:K}; the \(K\)-theoretic \(k\)-Schur functions represent the ideal-sheaf Schubert basis. 
The equivariant theory was developed by Lam and Shimozono \cite{LS:A} in homology and, more recently, by Ikeda, Shimozono, and Yamaguchi in \(K\)-homology \cite{ISY}. In the latter work, the two Schubert bases are realized by the \(K\)-theoretic double \(k\)-Schur functions and their closed versions, where the closed versions represent the structure-sheaf basis.

In type \(\mathrm{C}\), the non-equivariant homology of the affine Grassmannian was described by Lam, Schilling, and Shimozono using analogues of Schur  \(P\)-functions \cite{LSS:C}. Recently, a symmetric-function realization of equivariant homology of affine Grassmannians of type $\mathrm{C}$ was developed in \cite{IIS} and the equivariant $K$-homology counterpart in \cite{ISY2}.
For types $\B{}$ and $\D{}$, Pon \cite{Pon} constructed symmetric-function models in non-equivariant homology, whereas corresponding equivariant symmetric-function realization had remained unavailable.

The type $\D{}$ case presents additional subtleties.
As Lam observed in \cite[Remark 1] {Lam2} (see also \cite[Conjecture 2.8]{Pon}), the homology of the
affine Grassmannian of $\operatorname{Spin}(4m)$ does not
admit a graded Hopf algebra embedding into the ordinary
ring of symmetric functions.
Indeed, over $\mathbb{Q}$, its primitive subspace in
homological degree $4m-2$ is two-dimensional, whereas
the primitive subspace of the symmetric-function Hopf
algebra has dimension at most one in each degree.
Thus a symmetric-function realization compatible with
the graded Hopf algebra structure cannot be expected
in non-equivariant homology in general.
Our main result shows that an injective symmetric-function
realization can nevertheless be constructed for
torus-equivariant $K$-homology in type $\D{}$.
The equivariant parameters are essential: specializing
them to zero makes distinct Schubert classes coincide,
as we discuss in Section~\ref{sec:non-eq}.

Together with the theories in types $\A{}$ and $\mathrm{C}$, the present work completes the symmetric-function realizations of the torus-equivariant $K$-homology of affine Grassmannians in all classical types. Beyond providing explicit representatives for the Schubert bases, these realizations open the way to studying Pieri and Giambelli formulas, positivity of equivariant \(K\)-homology structure constants, and the quantum \(K\)-theory of orthogonal flag varieties through the ``quantum equals affine'' correspondence \cite{Kato}.

\subsection{Torus-equivariant $K$-homology of the affine Grassmannian}
Let $\Gad$ be a complex connected semisimple linear algebraic group
of adjoint type, and $\Gsc$ the simply connected 
group associated to $\Gad$. Then $\Gad=\Gsc/Z(\Gsc)$, where 
$Z(\Gsc)$ is the center of $\Gsc$.  We choose maximal tori
$\Tad\subset \Gad$ and $\Tsc\subset \Gsc$ that correspond under the quotient map \(\Gsc\to \Gad\).
The affine Grassmannian of $G$ is $\Gr_G = G(\mathbb{F})/G(\mathbb{O})$ where $\mathbb{F}={\C}(\!(t)\!)$ and $\mathbb{O}=\C[\![t]\!]$.
Let $K^T_*(\Gr_G)$ be the $T$-equivariant $K$-homology of $\Gr_G$, the Grothendieck group of $T$-equivariant coherent sheaves on $\Gr_G$ with proper support.
It has the Pontryagin product, making it a commutative algebra over the representation ring $R(T)$ of $T$. We refer the reader to \cite{LSS:K,ISY} for more details.


\subsection{The weight lattice and the representation ring}\label{sec:weights}
Let \(P\) and \(Q\) denote the weight and root lattices of $G$, respectively, and let \(P^\lor\) and \(Q^\lor\) denote its coweight and coroot lattices, respectively.
Letting $X^*(S)=\mathrm{Hom}_{\text{alg. gp.}}(S,\C^\times)$ be the character lattice and $X_*(S)=\mathrm{Hom}_{\text{alg. gp.}}(\C^\times,S)$ the cocharacter lattice 
for any algebraic torus $S$, we have
\[
P \cong X^*(\Tsc),\quad Q^\lor \cong X_*(\Tsc),
\quad
Q\cong X^*(\Tad),\quad P^\lor\cong X_*(\Tad).
\]
{
Denote by $\pair{\cdot}{\cdot}: P^\lor \otimes_{\Z} Q \to \Z$ the 
natural pairing given by function composition: for $\mu\in P^\lor$ and $\alpha\in Q$,
$\pair{\mu}{\alpha} = \alpha \circ \mu \in \mathrm{Hom}_{\text{alg. gp.}}(\C^\times,\C^\times)\cong \Z$.}
The representation ring $\Rad$ is a subring of 
$R(\Tsc)$ with natural identification
$R(\Tsc)\cong\Z[P]$, $\Rad\cong\Z[Q]$. Note that $\Rsc$ is a finite free $\RT$-module, with basis
indexed by the group $\Sigma \cong P/Q$ (see \S \ref{SS:extended}).

Since the Schubert structure constants lie in $\RT=\Z[Q]$,
we retain $\RT$ as the base ring for our Schubert calculus.
For the construction of their symmetric-function representatives,
we use the larger ring $\R$ introduced in Section~1.6
as an ambient coefficient ring.
\subsection{Extended affine Weyl groups}
\label{SS:extended}
The extended affine Weyl group $\Wafex$ is defined by 
\begin{align}\label{E:Wx as Pvee W}
\Wafex := P^\lor\rtimes W.
\end{align}
There is an injective group homomorphism $P^\lor \to \Wafex$ denoted $\mu\mapsto t_\mu$. 
The finite Weyl group $W$ is generated by simple reflections $s_i$ for $i$ in the classical Dynkin node set $I$. The group $W$ acts naturally on $P^\lor$ by \eqref{E:simple reflection on coweight}. 
In $\Wx$ we have the relation
$w t_\mu w^{-1} = t_{w(\mu)}$ for $w\in W$ and $\mu\in P^\lor$.

The affine Weyl group $\Waf$ is the subgroup of $\Wafex$ defined by $\Waf = Q^\lor \rtimes W$.  Let $I_\af := I\cup\{0\}$ denote the set of nodes for the untwisted affine Dynkin diagram.
Then $\Waf$ is a Coxeter group with generators $\{s_i\mid i\in I_\af\}$. 
If \(\theta\in Q\) denotes the highest root and $\theta^\vee\in Q^\vee$ is its associated coroot, then the affine simple reflection \(s_0\) can be written as
\begin{align}
s_0 = s_\theta t_{-\theta^\lor} = t_{\theta^\lor} s_\theta.
\end{align}
{
Let $\Gamma$ be an affine Dynkin diagram and $\Aut(\Gamma)$ the group of automorphisms of $\Gamma$, realized as a subgroup of permutations of the affine Dynkin node set $I_\af$. 

Define the subset $I_\af^s\subset I_\af$ of special (cominuscule) nodes, to be the orbit $\Aut(\Gamma)\cdot 0$ of the zero node $0\in I_\af$ under the group $\Aut(\Gamma)$.
There is a bijection $I_\af^s\to P^\lor/Q^\lor$ defined by
$i\mapsto \varpi_i^\lor+Q^\lor$ where $\varpi_i^\vee\in P^\vee$ is the fundamental coweight, with the convention $\varpi_0^\lor=0$. For $i\in I_\af^s$, the operation of addition by $\varpi_i^\lor$ defines a permutation of $P^\lor/Q^\lor$ and therefore induces a permutation of $I_\af^s$, which extends uniquely to an element $\pi_i\in \Aut(\Gamma)$ that sends $0$ to $i$. The elements $\Fund = \{\pi_i\mid i\in I_\af^s\}$ form a subgroup of $\Wx$ that is isomorphic to $P^\lor/Q^\lor$. 
}
The extended affine Weyl group $\Wx$ has a presentation 
\begin{align}\label{E:extended by fund}
\Wx\cong \Fund {\ltimes} W_\af
\end{align}
where $\Fund$ acts on $W_\af$ by affine Dynkin automorphisms. In $\Wx$ we have
$\sigma w \sigma^{-1} = \sigma(w)$ for $\sigma\in\Sigma$ and $w\in \Waf$.

Comparing the two presentations for $\Wx$ given by \eqref{E:Wx as Pvee W} and
\eqref{E:extended by fund} we have
\begin{align}\label{eq_act_pi}
   \pi_i &= t_{\varpi_i^\lor} u_i\qquad\text{for $i\in I_\af^s$}
\end{align}
where $u_i\in W$ is the unique shortest element such that $u_i^{-1}(\varpi_i^\lor)$ is antidominant. Note that
$u_i$ has the unique left descent $s_i$ for $i\neq 0$. We define $u_0=e$.

{We extend the Coxeter length function to \(\widetilde W_{\mathrm{af}}\) by \(\ell(\sigma w)=\ell(w)\) for \(\sigma\in\Sigma\), \(w\in W_{\mathrm{af}}\).}
The set of affine Grassmannian elements is defined by
\[
\Wafexgr:=\{w\in \Wafex\mid \ell(ws_i)>\ell(w)\quad
\text{for all $i\in I$}\}.
\]
Defining $\Wafgr=\Wafexgr\cap \Waf$, we have 
$\Wafexgr=\Sigma\cdot  \Wafgr$.

\subsection{Schubert basis}
\label{SS:Schubert basis}
Let $\Iwa$ denote the Iwahori subgroup, which is the 
preimage of $B\subset G$ under the evaluation map
$G(\mathbb{O})\rightarrow G$ at $t=0.$
For $w\in \Wafexgr$, the \emph{Schubert variety} is defined by
$X_w=\overline{\Iwa w G(\mathbb{O})/G(\mathbb{O})}$.
The variety $X_w$ is projective of complex dimension $\ell(w)$.
We have $$\Gr_G=\bigcup_{w\in \Wafexgr} X_w.$$  
{
The torus-equivariant $K$-homology $K_*^T(\Gr_G)$ of $\Gr_G$ is a commutative $\RT$-algebra with natural $\RT$-basis consisting of 
 the classes $\cO{w}:=[\cO{X_w}]$ of structure
sheaves $\cO{X_w}$ of $X_w$. 
The classes
$
\mathscr{I}_w:=[\mathscr{I}_{\partial X_w}]
$
of the ideal sheaves of the boundaries $\partial X_w\subset X_w$ form another $\RT$-basis of $K_*^T(\Gr_G)$. The algebra $K_*^T(\Gr_G)$ carries the Pontryagin product, which is most easily described on 
the fixed-point basis 
of \(\operatorname{Frac}(R(T))\otimes_{R(T)}K_*^T(\mathrm{Gr}_G)\)
 given by the equivariant $K$-classes $[e_x]$ of 
$T$-fixed points $x$. Using coset representatives given by translation elements, the product is given by
\begin{align}\label{E:pontryagin}
[e_{t_\lambda}] [e_{t_\mu}] = [e_{t_{\lambda+\mu}}] \qquad\text{for all $\lambda,\mu\in P^\vee$.}
\end{align}
}

\subsection{Connected components}
\label{SS:components}
{The finite group $\Sigma$ of \S \ref{SS:extended} is identified with the fundamental group of $G$, and also parametrizes the connected components of $\Gr_G$ (see \cite[Theorem 1.3.11 (1)]{Zhu}):
\begin{equation}
\Sigma\cong
\pi_0(\Gr_{\Gad})
\cong
\pi_1(\Gad)
\cong
P^\lor/Q^\lor.
\end{equation}
The $\RT$-algebra structure on $K_*^T(\Gr_G)$ 
is $\Sigma$-graded, 
and the component corresponding to $\bar{0}\in P^\lor/Q^\lor$
forms an $\RT$-subalgebra, such that
\begin{equation}
\label{eq:K_Gr_Ts}
K_{*}^{\Tsc}(\Gr_{\Gsc})\cong
R(\Tsc)
\otimes_{\RT}K_*^{T}(\Gr_G^{\bar{0}})    
\end{equation}
as $R(\Tsc)$-algebras.
This follows from the isomorphism
$\Gr_{\Gsc}\cong\Gr_G^{\bar{0}}$
induced by the covering $\Gsc\to G$
\cite[Theorem~1.3.11(3)]{Zhu},
together with extension of scalars along $\RT\to R(\Tsc)$.}
\subsection{Affine double $gq$-and $gd$-functions}
\label{SS:introBD}
Now we restrict ourselves to the cases $\B{n}$ and $\D{n}.$
We use a standard embedding $Q\hookrightarrow \bigoplus_{i=1}^n \Z a_i$ (see \S \ref{SS:root data})
and realize the representation ring as $$\RT=  \Z[Q] \subset \sR:=\Q[e^{\pm a_i} \mid 1\le i\le n] =
\Q[b_1,\bar{b}_1,\dotsc,b_n,\bar{b}_n]$$ where $b_i := 1 - e^{-a_i}$ for $1\le i\le n$ and $\bar{x} = -x/(1-x)$.

Let $\Lambda_{\sR}$ denote the ring of symmetric functions in the variables $y=(y_1,y_2,\ldots)$ with coefficients in $\sR$. We write $\hat{\Lambda}_{\sR} = \prod_{d\ge0} (\Lambda_{\sR})_d$ for its completion
where $(\Lambda_{\sR})_d$ is the component of degree $d$ in the $y$ variables. The action of $W$ on the coefficient ring $\sR$ extends naturally to an action on $\hat{\Lambda}_{\sR}$: $W$ acts on the coefficients and leaves the variables $y_i$ fixed. We extend the action of $W$ on $\hL_{\sR}$ to an action of $\Wafex$: for $\lambda\in P^\lor$, the translation element $t_\lambda$ acts by multiplication by $\tf{\lambda}\in \hL_{\sR}$ defined as follows.

For each $i\in I$, define the element
$\Omega(b_i|y)\in\hat{\Lambda}_{\sR}$ by
\begin{equation}
\Omega(b_i|y) =\prod_{j=1}^\infty \frac{1-\bar{b}_iy_j}{
1-b_i y_j
}.
\end{equation}
Let $\Omega(b_i| y)^{1/2}$ denote the unique element of
$\hL_{\sR}$ with constant term $1$ whose square is
$\Omega(b_i| y)$. 
In particular,
$\Omega(b_i| y)^{\pm 1/2}\in\hL_{\sR}$.
For an element $\la=\sum_{i=1}^n m_i \eps_i$ of the coweight lattice $P^\lor\subset \bigoplus_{i=1}^n (1/2) \Z \eps_i$ 
(see \S \ref{SS:root data}) define 
\begin{align}
  \tf{\la} = \prod_{i=1}^n \Omega(b_i| y)^{m_i}.
\end{align}
It is clear that $\tf{\la+\mu} = \tf{\la} \tf{\mu}$ for $\la,\mu\in P^\lor$
and $u(\tf{\la})=\tf{u(\la)}$ for $u\in W$.


Set $\alpha_0=-\theta$ where $\theta$ is the highest root.
Define $\Q$-linear operators $\hat{T}_i$ and
$\hat{D}_i$ on $\hL_{\R}$ for $i\in I_{\af}$ by
\begin{equation}\label{E:TDc}
\hat{T}_i=c(\alpha_i)^{-1}(s_i-1),\quad
\hat{D}_i=\hat{T}_i+1,\,\quad  c(\alpha)=1-e^{\alpha}.
\end{equation}
Both families of operators 
$\{\hat{T}_i\}$, $\{\hat{D}_i\}$ satisfy the braid relations, so $\hat{T}_w, \hat{D}_w$ for $w\in W_\af$ are well-defined. By 
\eqref{E:extended by fund} an arbitrary element of $\Wx$ has the form
$\sigma w$ for $\sigma\in\Sigma$ and $w\in \Waf$.
The operators $\hat{T}_{\sigma w}:=\sigma\hat{T}_w$ and $\hat{D}_{\sigma w}:=\sigma\hat{D}_w$
are well defined on $\hL_{\sR}$ and satisfy $\hat{T}_{uv}=\hat{T}_u\hat{T}_v$ and
$\hat{D}_{uv}=\hat{D}_u\hat{D}_v$ for all $u,v\in\Wx$ such that $\ell(uv)=\ell(u)+\ell(v)$.

For $w\in \Wxgr$, we define
\begin{equation}
\gx{w}{n}(y| b)=\hat{T}_w (1_{\hL_{\sR}}),
\quad
\tgx{w}{n}(y| b)=\hat{D}_w (1_{\hL_{\sR}}),
\end{equation}
where $\gx{w}{n}(y| b)$ (resp. $\tgx{w}{n}(y| b)$) is
$\gq{w}{n}(y| b)$ (resp. $\tgq{w}{n}(y| b)$) for type $\B{n}$
and $\gd{w}{n}(y| b)$ (resp. $\tgd{w}{n}(y| b)$) for type $\D{n}.$ We call \(gq^{(n)}_w(y| b)\) and \(\widetilde{gq}^{(n)}_w(y| b)\) the \emph{affine double \(gq\)-function} and the \emph{closed affine double \(gq\)-function}, respectively. Similarly, \(gd^{(n)}_w(y| b)\) and \(\widetilde{gd}^{(n)}_w(y| b)\) are called the \emph{affine double \(gd\)-function} and the \emph{closed affine double \(gd\)-function}, respectively. Let $\gGtil_{G}$ denote the $\RT$-span of $\gx{w}{n}(y| b)$ for $w\in \Wafexgr$ of the corresponding types.

\begin{thmA}[Theorem \ref{T:Peterson to symm funcs}]\label{thm_main} 
For types $\B{n}$ and $\D{n}$ there is an $R(T)$-algebra isomorphism
$$K_*^T(\Gr_G)\cong \gGtil_{G},$$ identifying the 
Schubert structure sheaf $\cO{w}$ 
with
$\tgx{w}{n}(y| b)$, 
and the ideal sheaf $\mathscr{I}_w$ of the boundary with $\gx{w}{n}(y| b)$ for $w$ in $\Wafexgr$ of the corresponding types.
\end{thmA}
\begin{remark}
The isomorphism in Theorem~\ref{T:Peterson to symm funcs} is compatible with a natural action of the equivariant $K$-theoretic affine nil-Hecke algebra. We do not introduce this algebra here, since it is not needed for the statement of the theorem. Its definition and the precise compatibility of the above isomorphism with the nil-Hecke action will be given later; see Proposition~\ref{prop_kg_mod}.
\end{remark}

We denote by $\gGhat_{G}$ the $R(T)$-span of $\gx{w}{n}(y| b)$ for $w\in \Wafgr.$  Then $\gGhat_{G}$ is naturally identified with $K_*^{\Tad}(\Gr_{\Gad}^{\bar{0}})$ and we have an isomorphism of $R(\Tsc)$-algebras
\begin{equation*}
K_*^{\Tsc}(\Gr_{\Gsc})\cong R(\Tsc)\otimes_{R(\Tad)}{\gGhat_{G}}
\end{equation*}
with the Schubert structure sheaf $\cO{w}$
identified with $\widetilde{gx}_w^{(n)}(y| b)$ and the ideal sheaf $\mathscr{I}_w$ of the boundary with $\gx{w}{n}(y| b)$ for all $w\in \Wafgr$ of the corresponding types. 

In type $\B{n}$, the $\Z$-span of the functions $\tgq{w}{n}(y| 0)$ obtained as specialization at $b=0$ form a Schubert basis of the non-equivariant $K$-homology ring $K_*(\mathrm{Gr}_{\Spin_{2n+1}(\mathbb C)})$ and we have a $\Z$-algebra isomorphism  
\begin{equation*}
K_*(\Gr_{\Spin_{2n+1}(\mathbb C)})\cong \bigoplus_{w\in \Wafgr}\Z\,\tgq{w}{n}(y| 0).
\end{equation*}
 In type $\D{n}$, specialization at $b=0$ identifies the symmetric
functions indexed by the two colors of the same partition
(Theorem~\ref{T:gd auto}), so these specialized functions do not form a
Schubert basis of the non-equivariant $K$-homology ring.
Their $\mathbb Z$-span nevertheless forms a quotient of
$K_*(\mathrm{Gr}_{\operatorname{Spin}_{2n}(\mathbb C)})$.
We conjecture that the kernel is generated by the difference
of the two Schubert classes indexed by the twin Pieri factors
(Conjecture~\ref{conj:D}).
\subsection{Relation to $\D{\infty}$ case for ``small'' classes}
\begin{table}[t]
\centering
\small
\setlength{\tabcolsep}{5pt}
\renewcommand{\arraystretch}{1.2}
\begin{tabular}{@{}p{0.34\linewidth}p{0.49\linewidth}l@{}}
\toprule
Notation & Description & Reference\\
\midrule
\multicolumn{3}{@{}l}{\textit{Coefficient rings}}\\[2pt]
$\RT=\Z[Q]$
  & Representation ring of the torus of the adjoint group
  & \S\ref{sec:weights}\\
$\Rsc=\Z[P]$
  & Representation ring of the torus of the simply connected group
  & \S\ref{sec:weights}\\
$\sR=\Q[e^{\pm a_1},\ldots,e^{\pm a_n}]$
  & Ambient coefficients for $K$-homology
  & \S\ref{SS:introBD}\\
$R_\infty=\Z[b_i,\bar b_i\mid i\ge1]$
  & Coefficients for type $D_\infty$ $K$-theory
  & \S\ref{subsec_D_infty}\\
$S=\operatorname{Sym}_{\Z}(Q)$
  & Equivariant homology coefficient ring
  & App.~\ref{S:eq_hom}\\
$S_{\sco}=\operatorname{Sym}_{\Z}(P)$
  & Simply connected coefficient ring
  & App.~\ref{S:eq_hom}\\
$\Sc=\Q[a_1,\ldots,a_n]$
  & Ambient coefficients for homology
  & App.~\ref{S:eq_hom}\\
$S_\infty=\Z[a_i\mid i\ge1]$
  & Coefficients for type $D_\infty$ homology
  & App.~\ref{S:eq_hom}\\
\midrule
\multicolumn{3}{@{}l}{\textit{Rings of symmetric functions}}\\[2pt]
$\Lambda_{\sR},\ \hL_{\sR}$
  & Symmetric functions and their $y$-degree completion
  & \S\ref{SS:introBD}\\
$\hL_{R_\infty}$
  & Completed symmetric functions of type $D_\infty$
  & \S\ref{subsec_D_infty}\\
$\gGtil_G,\ \gGhat_G$
  & Extended and neutral-component Schubert algebras
    in $K$-homology
  & \S\ref{S:sym_fun_rea}\\
$\hL_{\Sc}$
  & Completed symmetric functions for equivariant homology
  & App.~\ref{S:eq_hom}\\
$\htL{G},\ \htK{G}$
  & Extended and neutral-component Schubert algebras
    in homology
  & App.~\ref{S:eq_hom}\\
\midrule
\multicolumn{3}{@{}l}{\textit{Schubert functions}}\\[2pt]
$gq_w^{(n)},\ gd_w^{(n)}$
  & Boundary ideal-sheaf classes in types B and D
  & \S\ref{S:sym_fun_rea}\\
$\widetilde{gq}_w^{(n)},\
 \widetilde{gd}_w^{(n)}$
  & Structure-sheaf classes in types B and D
  & \S\ref{S:sym_fun_rea}\\
$gd_\lambda(y|b)$
  & Type $D_\infty$ dual factorial functions
  & \S\ref{subsec_D_infty}\\
$\widehat B_w^{(n)},\
 \widehat D_w^{(n)}$
  & Equivariant homology Schubert functions
    in types B and D
  & App.~\ref{S:eq_hom}\\
$\widehat D_\lambda(y|a)$
  & Type $D_\infty$ dual factorial homology functions
  & App.~\ref{S:eq_hom}\\
\bottomrule
\end{tabular}
\caption{Notation for coefficient rings, symmetric-function
algebras, and Schubert functions.}
\label{tab:notation}
\end{table}

Nakagawa and Naruse \cite{NN:Uni}
introduced a family of symmetric functions $gd_\lambda(y| b)$ indexed by the strict partitions $\lambda$ depending on infinitely many equivariant parameters $(b_1,b_2,\ldots).$
These functions can be
characterized by
the Cauchy identity
\[
\sum_{\lambda}
GD_\lambda(x| b)gd_\lambda(y| b)=
\prod_{i,j}\frac{1-\bar{x}_i y_j}{1-x_iy_j },
\]
where $\lambda$ runs over strict partitions, and $GD_\lambda(x| b)$ (see \S \ref{subsec_D_infty}) can be obtained from the $K$-theoretic Schur $P$-function $GP_\lambda(x| b)$ introduced by Ikeda and Naruse \cite{IN}. 

We prove that, for ``small" $w\in \Wafgr$, the affine double $gq$-function $\gq{w}{n}(y| b)$ and the affine double $gd$-function $\gd{w}{n}(y| b)$ can be expressed using $gd_\lambda(y| b)$ after an appropriate specialization of equivariant parameters; see Theorem \ref{small class type BD}, Theorem \ref{thm_sm_cl_b} and Theorem \ref{thm_sm_cl_type_D} for the precise statements.

\subsection{Factorization with respect to $\kappa_i$  classes}
We also give an explicit factorization formula for the Schubert structure sheaves in $K_*^T(\Gr_G)$ following the general framework of \cite{ISY}.
The underlying factorization phenomenon goes back to Magyar \cite{Mag}: every equivariant affine Grassmannian Schubert class, in either homology or \(K\)-homology, can be expressed as the product of one of the \(|W|/|\Sigma|\) distinguished ``irreducible'' Schubert classes and a collection of Schubert classes corresponding to \(\{\kappa_i: i\in I\}\), possibly occurring with multiplicities. Under Peterson's ``quantum equals affine'' correspondence \cite{Pet,LS}, the Schubert classes for \(\kappa_i\) can be written explicitly as Laurent monomials in the quantum parameters.
For the non-equivariant homology of type $\mathrm{A}$ affine Grassmannian, the Schubert basis is given by the (ungraded or $t=1$) $k$-Schur functions of \cite{LLM}, \cite{LM} as proved by Lam \cite{Lam}. In that setting the factoring was observed by Lapointe and Morse \cite{LM}, who called the classes $\kappa_i$ ``$k$-rectangles". Yet in the equivariant setting the structure giving rise to this factorization does not become clear until one works in the ``extended" setting, meaning using a group $G$ of adjoint type and the extended affine Weyl group.

We study the classes $\kappa_i$ in types $\B{}$ and $\D{}$ and prove that each closed affine double $gq$ function $\tgq{w}{n}(y| b)$ and each closed affine double $gd$ function $\tgd{w}{n}(y| b)$ can be factorized into a $\kappa$-irreducible class times a collection of $\kappa_i$-classes (with multiplicity) for $i\in I$; see Theorem~\ref{thm_fac_type_B} and Theorem~\ref{thm_fac_type_D} respectively. 
We further characterize the \(\kappa\)-irreducible affine Grassmannian elements combinatorially in terms of their associated partitions. More precisely, we derive explicit criteria classifying the \(\kappa\)-irreducible elements in types \(\B{}\) and \(\D{}\); see Theorem~\ref{thm:kappa-irr-type-B} and Theorem~\ref{thm:kappa-irr-type-D} respectively.
 \subsection{Future problems}
 The symmetric-function realizations developed here also provide
a framework for studying branching under the natural inclusions
$B_n\subset B_{n+1}$, $D_n\subset D_{n+1}$, and
$D_n\subset B_n\subset D_{n+1}$.
It would be interesting to describe the induced maps on
equivariant $K$-homology explicitly in the Schubert bases.


\subsection{Organization}
The paper is organized as follows. In \S \ref{sec_alg_cons}, we present a general algebraic construction of the equivariant \(K\)-homology of affine Grassmannians and recall the basic properties of the level-zero affine nil-Hecke algebra. In \S \ref{sec_root system}, we review the root systems, extended affine Weyl groups, and Pieri factors of types \(\B{n}^{(1)}\), \(\D{n}^{(1)}\), and \(\D{\infty}\). We also study the corresponding affine Grassmannian elements and describe their parametrization in terms of partitions.

In \S \ref{S:translation-factorization}, we investigate the \(\kappa_i\)-classes in types \(\B{}\) and \(\D{}\), and express anti-dominant translation elements in terms of elements in $\Sigma$ and the \(\kappa_i\)-classes. In \S \ref{S:sym_fun_rea}, we prove our main result, Theorem~\ref{T:Peterson to symm funcs}, and construct the symmetric-function realization associated with anti-dominant translation elements. In addition, for type $\D{}$, \S \ref{subsec_dyn_sym} describes the action of the Dynkin diagram automorphism on these symmetric functions.

In \S \ref{subsec_D_infty}, we introduce the type \(D_\infty\) symmetric functions \(gd_\lambda(y| b)\), indexed by strict partitions \(\lambda\). The main result of this section gives an expansion of \(gd_\lambda(y| b)\) in terms of \(gd_\mu(y| c)\), where coefficients are expressed in terms of double Grothendieck polynomials. In \S \ref{S: small_cl}, we prove that, for certain \(w\in\Wafgr\), the symmetric functions \(\gq{w}{n}(y| b)\) and \(\gd{w}{n}(y| b)\) admit explicit expressions in terms of the functions \(gd_\lambda(y| b)\). In \S \ref{S:ka_irr_elt}, we derive explicit combinatorial criteria on partitions characterizing the \(\kappa\)-irreducible elements of types \(\B{}\) and \(\D{}\).

Finally, Appendix~\ref{Ap_qu_K_th} discusses the computation of structure constants in the quantum \(K\)-theory of flag varieties via the \(K\)-homology of affine Grassmannians. Appendix~\ref{S:eq_hom} develops a symmetric-function realization of the equivariant homology of affine Grassmannians, extending the non-equivariant realization of Pon \cite{Pon}. In Appendix~\ref{S: expansion}, we formulate conjectural expansions of the functions \(\tgq{w}{n}(y)\) and \(\tgd{w}{n}(y)\) for all $w\in \Wafgr$ in terms of \(\widetilde{gd}_\lambda(y)\).

\subsubsection*{Acknowledgements}
First of all, T.~I. thanks Ryo Kohira and Yuki Ozaki
who studied  $H_{*}^T(\Gr_G)$ in type $\B{n}$ in their master’s theses.
We thank Satoshi Naito and Kohei Yamaguchi for helpful discussions. This work was supported by JSPS KAKENHI Grant Numbers
24KF0258, 25KF0074, 
23K25772, 22K03239. 
The first author was supported by a JSPS Postdoctoral Fellowship for Research in Japan, during which this research was carried out.

\section{Torus-equivariant $K$-homology of the affine Grassmannian}\label{sec_alg_cons}
The \(K\)-theoretic Peterson isomorphism is most naturally described
using the affine Grassmannian of the adjoint group.  In this form, the
translation elements indexed by coweights are visible inside the same
space, and the component group \(P^\lor/Q^\lor\) given by the coweight lattice
mod the coroot lattice, records the resulting
degree shifts.  This point of view is also important for the connection
with integrable systems.

\subsection{$\Sigma$-grading of the $K$-homology of the adjoint affine Grassmannian}

Consider the affine Grassmannian $\Gr_{\Gad}$ of a complex connected semisimple linear algebraic group $G$ of adjoint type. For \(\gamma\in P^\lor\), let
$
\bar\gamma:=\gamma+Q^\lor\in P^\lor/Q^\lor\cong \Sigma.
$
We denote by
$
\Gr_{\Gad}^{\bar\gamma}
$
the connected component of \(\Gr_{\Gad}\) containing the translation
point \(t_\gamma \Gad(\mathbb{O})\).  Thus
{
\[
\Gr_{\Gad}
= \bigsqcup_{i\in I_\af^s} \Gr_{\Gad}^{\overline{\varpi_i^\vee}}
\]
and the neutral component denoted by
$\Gr_{\Gad}^{\bar 0}$
is isomorphic to $\Gr_{\Gsc}$ (see \cite[Theorem 1.3.11(3)]{Zhu}).

We have an
\(R(\Tad)\)-module decomposition
\[
K_*^{\Tad}(\Gr_{\Gad})
\cong
\bigoplus_{i\in I_\af^s}
K_*^{\Tad}
\bigl(\Gr_{\Gad}^{\overline{\varpi_i^\vee}}\bigr).
\]

For $i\in I_\af^s$ let $\xi_i = [e_{t_{\overline{\varpi_i^\vee}}}]$ be the equivariant point class.
The Pontryagin product with $\xi_i$ gives the equality
\[
K_*^{\Tad}(\Gr_{\Gad}^{\overline{\varpi_i^\vee}})
=
K_*^{\Tad}(\Gr_{\Gad}^{\bar 0})\, \xi_i.
\]
Thus there is an $R(T)$-algebra isomorphism
\[
K_*^{\Tad}(\Gr_{\Gad})
\cong
\bigoplus_{i\in I_\af^s}
K_*^{\Tad}(\Gr_{\Gad}^{\bar 0})
\,\xi_i
\]
with the product structure determined by \eqref{E:pontryagin}.}
In particular the Pontryagin product respects the component grading: for any $\gamma,\gamma'\in P^\vee$
\[
K_*^{\Tad}
\bigl(\Gr_{\Gad}^{\overline{\gamma}}\bigr)
\cdot K_*^{\Tad}
\bigl(\Gr_{\Gad}^{\overline{\gamma'}}\bigr)
\subset
K_*^{\Tad}
\bigl(\Gr_{\Gad}^{\overline{\gamma}+\overline{\gamma'}}\bigr).
\]

\subsection{Algebraic constructions of $K^T_*(\Gr_G)$}
The finite Weyl group $W$ acts naturally on { $R(T)\cong \Z[Q]$}.
The extended affine Weyl group $\Wx$ acts on $R(T)$ via the level-zero action, that is, via the group homomorphism 
$\cl:\Wx\to W$ given by $t_\la u\mapsto u$ for $\la\in P^\lor$ and $u\in W$. Let $R(T)^\Delta$ be the localization of $R(T)$ at the multiplicatively closed subset $\Delta$ generated by $\{c(\alpha) \mid \alpha\in \Phi\}\subset R(T)$, where $\Phi$ is the set of all roots of the finite root system. Since $\Wx$ preserves $\Delta$, $\Wx$ acts on $R(T)^\Delta$.

Let $\tKGd = R(T)^\Delta \otimes_\Z \Z[\Wx]$ be the twisted group algebra. Its multiplication is given by
\[
(f \otimes v)(g \otimes w) = f \, v(g) \otimes vw,
\]
for $f,g \in R(T)^\Delta$ and $v,w \in \Wx$. We regard
\[
 \KGd:= R(T)^\Delta \otimes_{\mathbb{Z}} \mathbb{Z}[\Waf]
\]
as a subalgebra of $\tKGd$. Note that $\tKGd$ acts on $R(T)^\Delta$ where $R(T)^\Delta$ acts on itself by left multiplication and $\Wx$ acts as above. Similarly $\KGd$ acts on $R(T)^\Delta$. For $i \in I_\af$, define $T_i$ and $D_i$ in $\KGd$ by
\begin{align}\label{E:T}
T_i &= c(\alpha_i)^{-1}(s_i-1) \\
\label{E:D}
D_i &= T_i+1 .
\end{align}
Here $\alpha_0 = -\theta$.
They satisfy $T_i^2 = -T_i$ and $D_i^2=D_i$. The elements $T_i$ and $D_i$ satisfy the braid relations. For $w \in W_\af$ one can define
\[
T_w:= T_{i_1} \cdots T_{i_\ell},~ D_w:= D_{i_1} \cdots D_{i_\ell} \text{ for a reduced expression $w = s_{i_1} \cdots s_{i_\ell}$.}
\]
We further define
\begin{equation}\label{eq:A-pi-v}
T_{\pi v} = \pi T_v, \quad D_{\pi v}=\pi D_v \quad \text{for } \pi \in \Sigma, \; v \in W_\af.
\end{equation}

Let $\KG$ be the subring of $\KG^\Delta$ generated by $R(T)$ and $T_i$ for $i\in I_\af$. Let $\tKG$ be the subring of $\tKGd$ generated by $\KG$ and $\pi\in \Fund$. 
We call $\KG$ (resp. $\tKG$) the \emph{level zero} (resp. \emph{extended}) affine nil-Hecke algebra. 
Note that $R(T)$ is not central in $\KG$ nor in $\tKG$.

One may show that 
$\{T_w \mid w\in \Wx\}$ is a left $R(T)$-basis of $\tKG$ and $\{T_w\mid w\in W_\af\}$ is a left $R(T)$-basis of $\KG$. The action of $\tKGd$ on $R(T)^\Delta$ restricts to an action of $\tKG$ on $R(T)$. 
Similarly $\KG$ acts on $R(T)$. We can regard $\KG$ as a subring of $\tKG$. One may show that $W_\af\subset \KG$: solving \eqref{E:T} for $s_i$ one has $s_i = 1+c(\alpha_i) T_i\in \KG$. Similarly $\Wx\subset \tKG$.

We have $\tKG \cong K^T_*(\tilde{\Fl}_G)$; see \cite{KashiwaraNotes}, where
$\tilde{\Fl}_G$ is the affine Kac-Moody flag ind-variety and the right hand side carries a noncommutative convolution product.

Define
$$\tLG=Z_{\tKG}(R(T)),\qquad \LG=Z_{\KG}(R(T)),$$
the centralizers of $R(T)$ in $\tKG$ and $\KG$, respectively.
Both $\tLG$ and $\LG$ are commutative $R(T)$-subalgebras.
Moreover, the elements $\{ t_{\lambda }\mid \lambda \in P^{\lor }\}$ form an $R(T)^{\Delta }$-basis of
$\tLG^{\Delta }=R(T)^{\Delta }\otimes _{R(T)}\tLG$,
and similarly, the elements  $\{t_{\lambda }\mid \lambda \in Q^{\lor }\}$  form an $R(T)^{\Delta }$-basis of
$\LG^{\Delta }=R(T)^{\Delta }\otimes _{R(T)}\LG$.

There is an action of $\tKG^\Delta$ on $\tLG^\Delta$ \cite{Kato} given by:
\begin{align}
\label{E:action W}
  u * (c \otimes t_\la ) &= u(c) \otimes t_{u(\la)} \\
\label{E:action t}
  t_\mu *(c \otimes t_\la) &= c \otimes t_{\mu+\la}  
\end{align}
for $u\in W$, $\la,\mu\in P^\lor$ and $c\in R(T)^\Delta$.
This restricts to an action of $\tKG$ on $\tLG$ and an action of $\KG$ on $\LG$.

Let $\kappa: \tKG\to\tLG$ be the map defined by $\kappa(a) = a*1_{\tLG}$ for $a\in \tKG$.
It satisfies
\begin{align}\label{E:kappa L linear}
    \kappa(ba) &= b \kappa(a)&\qquad&\text{for all $b\in\tLG$, $a\in\tKG$} \\
    \kappa(b) &= b &&\text{for all $b\in\tLG$.}
\end{align}
See \cite[\S 2.3.2]{ISY} for more properties of the map $\kappa.$
For $w\in \Wxgr$ define the elements $k_w,\ell_w\in \tLG$ by
\begin{align}\label{E:kw}
  k_w &= \kappa(T_w)=T_w* 1_{\tLG}, \\
  \label{E:lw}
  \ell_w &= \kappa(D_w)=D_w* 1_{\tLG}.
\end{align}

Then we have $$\tLG \cong \bigoplus_{w\in \Wxgr} R(T)k_w=\bigoplus_{w\in \Wxgr} R(T)\ell_w;\qquad \LG \cong \bigoplus_{w\in \Wafgr} R(T)k_w=\bigoplus_{w\in \Wafgr} R(T)\ell_w.$$

Using arguments analogous to those in \cite{LSS:K} for the simply-connected group $G_{\text{sc}}$ and  \cite{ISY} for adjoint group $G$ we have the following.
\begin{prop}
    The horizontal maps of the following diagram are $R(T)$-algebra isomorphisms:
\begin{equation}\label{E:geom to algebra iso}
\begin{tikzcd}
K^{\Tsc}_*(\Gr_{G_{\text{sc}}}) \arrow[r,"\cong"] \arrow[d,hookrightarrow] & R(\Tsc)\otimes_{R(T)}
\LG \arrow[d,hookrightarrow] \\
K^T_*(\Gr_G) \arrow[r, swap, "\psi"]\arrow[r,"\cong"] & \tLG 
\end{tikzcd}
\end{equation}
Moreover, we have $\psi(\mathscr{I}_w) = k_w$ and
$\psi(\cO{w})=\ell_w$ for all $w\in \Wafexgr$.
\end{prop}



\subsection{Anti-dominant translation elements and the elements $\kappa_i$}

The factors $\kappa_i$ we wish to describe are closely related to 
quantum parameters under the Quantum Equals Affine correspondence.

In light of \eqref{E:extended by fund} for $\la\in P^\lor$ write 
\begin{align}\label{E:translation auto affine}
t_\la = \sigma_\la u_\la\qquad\text{for $\sigma_\la\in\Fund$ and $u_\la\in \Waf$.}
\end{align}
It is known \cite[Lemma 2.33]{ISY} that $u_\la\in\Wafgr$ if $\la$ is antidominant.
For $i\in I$ define 
\begin{align}\label{E:kappa}
  \kappa_i = u_{-\varpi_i^\lor}\qquad\text{for $i\in I$.}
\end{align} 
Since $-\varpi_i^\vee$ is antidominant, $\kappa_i\in\Wafgr$. By convention define $\kappa_0=\id$.

\begin{lem} \label{L:kappa}
Let $i\in I$. 
If $i\in I_\af^s$ then \cite[Prop. 2.34]{ISY}
\begin{align}\label{E:kappa special}
\kappa_i = \pi_i u_i \pi_i^{-1}.
\end{align}
If $i\notin I_\af^s$, let $i^s\in I_\af^s$ be such that $
\varpi_i^\vee \equiv \varpi_{i^s}^\vee\in Q^\vee$. Then \cite[(2.67)]{ISY}
\begin{align}\label{E:kappa i}
\kappa_i =
\pi_{i^s} t_{-\varpi_i^\vee}
\end{align}
In particular if $i^s=0$ then $\kappa_i = t_{-\varpi_i^\vee}$.
\end{lem}

\begin{prop}[{\cite[Proposition 2.37]{ISY}}]\label{P:factor kappa}
Suppose $w\in\Wx$ and $i\in I$ are such that 
$w\kappa_{i}\in \Wxgr$ and $\ell(w\kappa_i)=\ell(w)+\ell(\kappa_i)$.
Then with $i^s$ defined as in Lemma \ref{L:kappa}, $\pi_{i^s}^{-1}(w)\in \Wxgr$ and 
\begin{align*}
  \ell_{w\kappa_i} &= \ell_{\kappa_i} u_{i^s}(\ell_{\pi_{i^s}^{-1}(w)}).
\end{align*}
\end{prop}

\section{Root systems of types $\mathrm{B}_n^{(1)}$, $\D{n}^{(1)}$, and $\D{\infty}$}\label{sec_root system}

We now go into specifics for types B and D. This is used to make symmetric function realizations of the Schubert bases $k_w$ and $\ell_w$ and study their properties, especially the factorization property and the explicit description of these functions for small classes.

\subsection{Dynkin diagrams}
\label{SS:Dynkin}
We use the root systems for the following Dynkin diagrams.
\begin{center}
\begin{tikzpicture}[node distance=1cm and 1cm, every node/.style={draw, circle, inner sep=2pt}]
    \node[label=below:{$2$}] (2) {};
\node[label=left:{$0$}, above left=1cm of 2] (0) {};
\node[label=left:{$1$}, below left=1cm of 2] (1) {};
      \node[label=below:{$3$}, right=of 2] (3) {};
      \begin{scope}[every node/.style={inner sep=2pt}]
        \node[right=of 3] (dots) {$\cdots$};
       \end{scope}
    \node[label={[yshift=2.6mm]below:{$n-1$}},right=of dots](n-1) {};
   \node[label=below:{$n$},right=of n-1](n) {};
   \begin{scope}[every node/.style={inner sep=2pt}]
     \draw[double distance=2pt] (n-1) -- (n) node[midway] {$>$}; 
    \end{scope}
    \draw(dots) -- (n-1);
    \draw (1) -- (2);
    \draw (2) -- (0);
    \draw (3) -- (2);
    \draw (3) -- (dots);
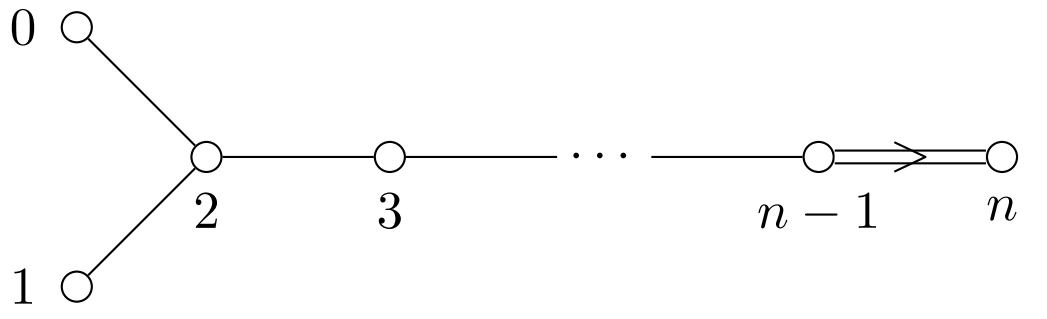
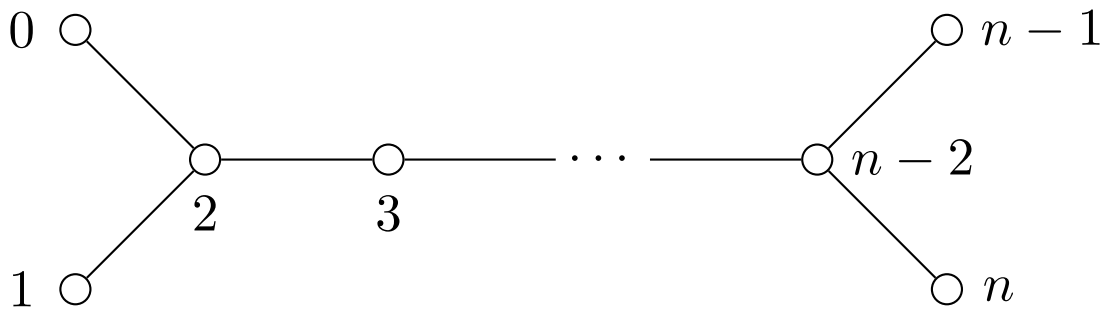
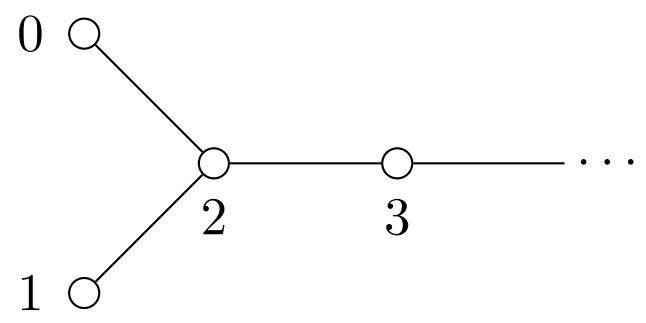
    \begin{scope}[every node/.style={inner sep=2pt}]
     \node[left=2cm of 2]{$\B{n}^{(1)}$};
     \node[right=5mm of n]{};
    \end{scope}
 \end{tikzpicture}
 \end{center}
\begin{center}
\begin{tikzpicture}[node distance=1cm and 1cm, every node/.style={draw, circle, inner sep=2pt}]
    \node[label=below:{$2$}] (2) {};
\node[label=left:{$0$}, above left=1cm of 2] (0) {};
\node[label=left:{$1$}, below left=1cm of 2] (1) {};
      \node[label=below:{$3$}, right=of 2] (3) {};
      \begin{scope}[every node/.style={inner sep=2pt}]
        \node[right=of 3] (dots) {$\cdots$};
       \end{scope}
    \node[label={right:{$n-2$}},right=of dots](n-2) {};
   \node[label=right:{$n-1$},above right=1cm of n-2](n-1) {};
    \node[label=right:{$n$},below right=1cm of n-2](n) {};
    \draw(dots) -- (n-2);
    \draw (1) -- (2);
    \draw (2) -- (0);
    \draw (3) -- (2);
    \draw (3) -- (dots);
    \draw (n-2)--(n-1);
    \draw (n-2)--(n);
    \begin{scope}[every node/.style={inner sep=2pt}]
     \node[left=2cm of 2]{$\D{n}^{(1)}$};
    \end{scope}
 \end{tikzpicture}
 \end{center}
\begin{center}
\begin{tikzpicture}[node distance=1cm and 1cm, every node/.style={draw, circle, inner sep=2pt}]
    \node[label=below:{$2$}] (2) {};
\node[label=left:{$0$}, above left=1cm of 2] (0) {};
\node[label=left:{$1$}, below left=1cm of 2] (1) {};
      \node[label=below:{$3$}, right=of 2] (3) {};
      \begin{scope}[every node/.style={inner sep=2pt}]
        \node[right=of 3] (dots) {$\cdots$};
       \end{scope}
       \begin{scope}[every node/.style={inner sep=2pt}]
        \node[right=2.9cm of dots] (x) {};
       \end{scope}
    \draw (1) -- (2);
    \draw (2) -- (0);
    \draw (3) -- (2);
    \draw (3) -- (dots);
    \begin{scope}[every node/.style={inner sep=2pt}]
     \node[left=2cm of 2]{$\D{\infty}$};
    \end{scope}
 \end{tikzpicture}
 \end{center}
 The Dynkin node sets are designated $I(\B{n}^{(1)})=I(\D{n}^{(1)})=\{0,1,\dots,n\}$ and $I(\D{\infty})=\Z_{\ge0}$. The classical Dynkin diagrams of types $\B{n}$ and $\D{n}$ are given by removing the $0$ nodes from $\B{n}^{(1)}$ and $\D{n}^{(1)}$ respectively. Their Dynkin node sets are denoted $I(\B{n})=I(\D{n})=\{1,2,\dotsc,n\}$. All three diagrams admit an automorphism $\auto$ that exchanges nodes $0$ and $1$ and fixes other nodes. $\auto$ generates $\Aut(\B{n}^{(1)})$ and $\Aut(\D{\infty})$. We will use $\autn\in \Aut(\D{n}^{(1)})$ which exchanges nodes $n-1$ and $n$ and fixes other nodes.
Type $\D{n}^{(1)}$ has an automorphism $\autDrev{n}$ that sends $i\mapsto n-i$ for $i\in I(\D{n}^{(1)})$.

\subsection{Classical root data}
\label{SS:root data}
Let $\mathfrak{g} = \mathfrak{so}_{2n+1}(\mathbb{C})$ for type $\mathrm{B}_n$ and
$\mathfrak{g}=\mathfrak{so}_{2n}(\mathbb{C})$ for type $\mathrm{D}_n$. Here we write $I=\{1,2,\dots,n\}$ for the classical Dynkin node set. Let $P = \bigoplus_{i\in I} \Z \varpi_i\subset \bigoplus_{k=1}^n \mathbb{Q} a_k$ be the weight lattice, with fundamental weights
\begin{align*}
  \varpi_i &= a_1 + \dotsm + a_i\qquad\text{$1\le i\le n-1$ for $\mathrm{B}_n$ and $1\le i\le n-2$ for $\mathrm{D}_n$}, \\
  \varpi_{n-1} &= \half(a_1+a_2+\dotsm+a_{n-1}-a_n) \qquad\text{for $\mathrm{D}_n$}, \\
  \varpi_n &= \half\sum_{k=1}^n a_k \qquad\text{for $\mathrm{B}_n$ and $\mathrm{D}_n$}.
\end{align*}
Let $Q =\bigoplus_{i=1}^n \Z \alpha_i \subset P$ be the root lattice, with simple roots 
\begin{align}
\label{E:simple A root}
  \alpha_i &= a_i - a_{i+1} \qquad\text{for $1\le i\le n-1$}, \\
\label{E:simple BD n root}
  \alpha_n &= \begin{cases}
      a_n & \text{for $\mathrm{B}_n$} \\
      a_{n-1}+a_n &\text{for $\mathrm{D}_n$.}
  \end{cases}
\end{align}
We use the coweight lattice $P^\lor=\bigoplus_{i=1}^n \Z \varpi_i^\lor \subset \bigoplus_{k=1}^n \Q \eps_k$, with
fundamental coweights 
\begin{align*}
  \varpi_i^\lor &= \eps_1 + \dotsm + \eps_i\qquad\text{$1\leq i\le n$ for $\mathrm{B}_n$ and $1\leq i\le n-2$ for $\mathrm{D}_n$}, \\
  \varpi_{n-1}^\lor &= \half(\eps_1+\eps_2+\dotsm+\eps_{n-1}-\eps_n)\qquad\text{for $\mathrm{D}_n$}, \\
\varpi_{n}^\lor &= \half(\eps_1+\eps_2+\dotsm+\eps_{n-1}+\eps_n)\qquad\text{for $\mathrm{D}_n$},
\end{align*}
and coroot lattice $Q^\lor =\bigoplus_{i=1}^n \Z \alpha_i^\lor \subset P^\lor$, with simple coroots
\begin{align}
\label{E:simple A coroot}
  \alpha_i^\lor &= \eps_i - \eps_{i+1} \qquad\text{for $1\le i\le n-1$}, \\
\label{E:simple BD n coroot}
  \alpha_n^\lor &= \begin{cases}
      2\eps_n & \text{for $\B{n}$} \\
      \eps_{n-1}+\eps_n &\text{for $\D{n}$.}
  \end{cases}
\end{align}
The highest root $\theta$ and associated coroot $\theta^\lor$ are given by
\begin{align}
\label{E:highest root and associated coroot}
    \theta &= a_1 + a_2 \qquad \theta^\lor = \eps_1 + \eps_2.
\end{align}
We use the notation 
\begin{align}\label{E:zero root and coroot}
    \alpha_0 = -\theta\qquad \alpha_0^\lor = -\theta^\lor.
\end{align}
There is an evaluation pairing $\pair{\cdot}{\cdot}: (\bigoplus_{i=1}^n \Q \eps_i) \times (\bigoplus_{i=1}^n \Q a_i) \to \Q$ such that $\pair{\eps_i}{a_j}=\delta_{ij}$ for $1\le i,j\le n$.
Under this pairing $Q^\lor$ and $P$ are dual over $\Z$ with respective dual bases
$\{\alpha_i^\lor\mid i\in I\}$ and $\{\varpi_i\mid i\in I\}$.
Similarly $P^\lor$ and $Q$ are dual over $\Z$ with respective dual bases
$\{\varpi^\lor_i\mid i\in I\}$ and $\{\alpha_i\mid i\in I\}$.

The Dynkin node set of type $\mathrm{D}_\infty$ is $I(\mathrm{D}_\infty)=\Z_{\ge0}$.
We use the lattice $\Xinf = \bigoplus_{k\ge1} \Z a_k$ 
and its graded dual $\Xvinf = \bigoplus_{k\ge1} \Z\eps_k$ so that 
$\pair{\cdot}{\cdot}: \Xvinf\otimes \Xinf\to \Z$ is given by
$\pair{\eps_i}{a_j}=\delta_{ij}$ for $i,j\in \Z_{\ge 1}$. The simple roots
$\alpha_i\in \Xinf$ are defined by \eqref{E:zero root and coroot} for $i=0$ 
and \eqref{E:simple A root} for $\alpha_i$ for $i\ge1$. The simple coroots
$\alpha_i^\lor\in \Xvinf$ are defined by \eqref{E:zero root and coroot} for $i=0$
and \eqref{E:simple A coroot} for $i\ge1$. 

\subsection{Classical Weyl group}
\label{SS:classical weyl}
The classical Weyl group $W$ is generated by $\{s_i\mid i\in I\}$ and satisfies braid relations $s_i^2=1$ and $s_i s_j s_i = s_j s_i s_j$ for $i,j\in I$ that are adjacent by a single bond, $s_{n-1}s_ns_{n-1}s_n=s_ns_{n-1}s_ns_{n-1}$ for type $\mathrm{B}_n$, and otherwise $s_i$ and $s_j$ commute. $W$ acts on $P$ and $P^\lor$ by 
\begin{align}
\label{E:simple reflection on weight}
  s_i(\la) &= \la - \pair{\alpha_i^\lor}{\la} \alpha_i&\qquad&\text{for $i\in I$, $\la\in P$} \\
\label{E:simple reflection on coweight}
  s_i(\mu) &= \mu -  \pair{\mu}{\alpha_i}\alpha_i^\lor&&\text{for $\mu\in P^\lor$.}
\end{align}
Explicitly, for $1\le i\le n-1$, $s_i(a_i)=a_{i+1}$, $s_i(a_{i+1})=a_i$, and $s_i(a_k)=a_k$ for other $k$. For $\mathrm{B}_n$, $s_n(a_n)=-a_n$ and $s_n(a_k)=a_k$ for other $k$. For $\mathrm{D}_n$, 
$s_n(a_{n-1})=-a_n$, $s_n(a_n)=-a_{n-1}$ and $s_n(a_k)=a_k$ for other $k$. For the action of $W$ on $P^\lor$, just replace $a_k$ by $\eps_k$ for all $k$. The reflection associated with the highest root $\theta$ is given by 
$s_\theta(a_1)=-a_2$, $s_\theta(a_2)=-a_1$, and $s_\theta(a_k)=a_k$ for other $k$.
The action of $W$ stabilizes $Q\subset \bigoplus_{k=1}^n \Q a_k$
and $Q^\lor \subset \bigoplus_{k=1}^n \Q \eps_k$.

The Weyl group $\WD{\infty}$ is generated by 
$\{s_i\mid i\geq 0\}$ satisfying $s_i^2=1$ for $i\ge 0$ and the braid relations $s_is_{j}s_i=s_{j}s_is_{j}$ if there is an edge between $i$ and $j$ and otherwise
$s_i$ and $s_j$ commute. 
The Weyl group $\WD{\infty}$ acts on $\Xinf$ by 
\eqref{E:simple reflection on weight} and on $\Xvinf$ by  
\eqref{E:simple reflection on coweight} but for $i\in I(\mathrm{D}_\infty)$, $\la\in \Xinf$ and $\mu\in \Xvinf$. The action of $s_0$ is same as the action of $s_\theta$.

\subsection{Affine Dynkin automorphisms and the fundamental group $\Fund$}

The special nodes in type $\B{n}^{(1)}$ and $\D{n}^{(1)}$ are given by
\begin{align*}
  I_\af^s = \begin{cases}
      \{0,1\} & \text{for $\mathrm{B}_n^{(1)}$} \\
      \{0,1,n-1,n\} &\text{for $\mathrm{D}_n^{(1)}$.}
  \end{cases}
\end{align*}
The action of elements of $\Fund$ can be described explicitly by their cycle structure as a permutation of $I_\af^s$ as follows, computing from the general rules in \S \ref{SS:extended} involving addition by simple coroots in $P^\vee/Q^\vee$.
These permutations of $I_\af^s$ extend uniquely to affine Dynkin automorphisms.

We have $\Fund = \langle \pi_1 \rangle \cong \Z/2\Z$ for type $\mathrm{B}_n^{(1)}$: as a permutation of $I_\af$, $\pi_1$ exchanges $0$ and $1$ and fixes other nodes. 
For type $\mathrm{D}_n^{(1)}$ and $n$ even, $\Fund\cong \Z/2\Z\times \Z/2\Z$:  $\pi_0$ is the identity, $\pi_1 = (0,1)(n-1,n)$,
$\pi_{n-1}=(0,n-1)(1,n)$
and $\pi_n=(0,n)(1,n-1)$.
For type $\mathrm{D}_n^{(1)}$ and $n$ odd, $\Fund\cong \Z/4\Z$ with $\pi_0$ the identity,
$\pi_1=(0,1)(n-1,n)$, $\pi_{n-1}=(0,n-1,1,n)$, $\pi_n=(0,n,1,n-1)$. For types $\mathrm{B}_n^{(1)}$ and $\mathrm{D}_n^{(1)}$, for all $i\in I^s_\af$ and $j\in I_\af\setminus I^s_\af$ we have $\pi_i(j)=j$ for $i\in\{0,1\}$ and $\pi_i(j)=n-j$ for $i\in \{n-1,n\}$.

\begin{remark} \label{R:switching Dynkin autos}
Any Dynkin automorphism induces an automorphism of the Weyl group, 
the weight lattice, and the root lattice. 
For any of $\mathrm{B}_n^{(1)}$, $\mathrm{D}_n^{(1)}$, or $\mathrm{D}_\infty$, the Dynkin automorphism $\auto$ (see \S \ref{SS:Dynkin}) acts by
\begin{align}\label{E:auto 0}
\auto(a_1) &= -a_1 \\
\auto(a_k) &= a_k \qquad\text{for $k\ge 2$.}
\end{align}
Similarly the Dynkin automorphism $\autn$ of $\mathrm{D}_n$ (or $\mathrm{D}_n^{(1)}$) acts by 
\begin{align}\label{E:auto n}
\autn(a_n) &= -a_n \\
\autn(a_k) &= a_k \qquad\text{for $1\le k\le n-1$.}
\end{align}
\end{remark}
\begin{remark}\label{R:cosets of coweights}
In type $\mathrm{B}_n$, $\varpi_i^\lor\in Q^\lor$ if and only if $i$ is even,
and $\varpi_i^\lor\equiv \varpi_1^\lor \mod Q^\lor$ if $i$ is odd.
In type $\mathrm{D}_n$ the same holds for $1\le i\le n-2$. 
\end{remark}

\subsection{Pieri factors}
The notion of Pieri factor, introduced in \cite{Pon}, gives a unified definition of ``single row" affine Grassmannian elements $\rho_i$ for untwisted affine types A,B,C and D and is used to define affine Stanley symmetric functions, which are representatives for Schubert classes in the cohomology of the affine Grassmannian (\cite{LM:ktab}, \cite{Lam:as}
for type A, \cite{LSS:C} for type C, and \cite{Pon} for types B and D). Consider the orbit $W\cdot \varpi_1^\lor = \{ \pm \eps_i \mid 1\le i\le n\}$.
We picture the orbits below, first for $\B{n}$ and then for $\D{n}$, where $\eps_i$ and $-\eps_i$ are pictured as $i$ and $\bar{i}$ respectively.
\[
\begin{tikzcd}[column sep=6mm]
 1 \arrow[r,"s_1"] & 2 \arrow[r,"s_2"] & \dotsm \arrow[r,"s_{n-2}"] & n-1 \arrow[r,"s_{n-1}"] & n \arrow[r,"s_n"] & \bar{n} \arrow[r,"s_{n-1}"] & \ba{n-1} \arrow[r,"s_{n-2}"] & \dotsm \arrow[r,"s_2"] & \ba{2} \arrow[r,"s_1"] & \ba{1} 
 \end{tikzcd}
\]
\[
\begin{tikzcd}[column sep = 6mm]
&&&& n \arrow[dr,"s_n"] \\
  1 \arrow[r,"s_1"] & 2 \arrow[r,"s_2"] & \dotsm \arrow[r,"s_{n-2}"] & n-1 \arrow[ur,"s_{n-1}"] \arrow[dr,swap,"s_n"] && \ba{n-1} \arrow[r,"s_{n-2}"] & \dotsm \arrow[r,"s_2"]& \bar{2} \arrow[r,"s_1"] & \bar{1} \\
  &&&& \ba{n} \arrow[ur,swap,"s_{n-1}"]
\end{tikzcd}
\]

For a coweight $\mu$, let $v_\mu\in W$ be the shortest element such that $v_\mu^{-1}(\mu)$ is \emph{dominant}.
For $1\le i\le 2n-1$ in type $\B{n}$ and
$1\le i\le 2n-2$ in type $\D{n}$, let $\rho_i^{(1)}$
denote the unique element of length $i$ among
$\{v_\mu\mid \mu\in W\cdot\varpi_1^\lor\}$,
whenever such an element is unique.
The only failure of uniqueness within these ranges is for type $\mathrm{D}_n$ and $i=n-1$, in which there are two elements we denote by
\begin{align}\label{E:two Pieri factors}
\rho_{n-1}^{+(1)} &= v_{\eps_n} = s_{n-1} s_{n-2}\dotsm s_2 s_1, \\
\rho_{n-1}^{-(1)} &= v_{-\eps_n} = s_n s_{n-2} \dotsm s_2 s_1.
\end{align}
We will call these \emph{twin Pieri factors}. For the other cases (\emph{nontwin factors})
\begin{align}
  \rho^{(1)}_i &= 
  \begin{cases}
      s_i \dotsm s_1&\text{$1\le i\le n$}\,\,(\text{$i\ne n-1$ if in type $\D{n}$}), \\
      s_{2n-i}\dotsm s_{n-2} s_{n-1}\rho^{(1)}_n&\text{for $\B{n}$, $n+1\le i\le 2n-1$},\\
       s_{2n-1-i} s_{2n-i}\dotsm s_{n-2} \rho^{(1)}_n &\text{for $\D{n}$, $n+1\le i\le 2n-2$.}
  \end{cases}
\end{align}

Define the nontwin Pieri factors 
\begin{align}\label{E:Pieri factors}
  \rho_i ^{(0)}= \auto(\rho_i^{(1)})
\end{align}
with $\auto$ as in \eqref{E:auto 0}, $1\le i\le 2n-1$ for type $\B{n}$ and $1\le i\le 2n-2$ ($i\ne n-1$) for type $\D{n}$.
Define the twin Pieri factors (type $\D{n}$) by
\begin{align}\label{E:twin Pieri factors}
\rho_{n-1}^{\mp(0)} &= \auto(\rho_{n-1}^{\pm(1)})
\end{align}
 In other words, $\rho_{n-1}^{\pm(0)} = \pi_1(\rho_{n-1}^{\pm(1)})$.
Since each $\rho_i^{(1)},\,\rho_{n-1}^{\pm (1)}$ is by definition a minimum length coset representative in $W$
with respect to the subgroup generated by $s_i$ for $i\in I\setminus\{1\}$ and $\auto$ is an affine Dynkin automorphism sending $1$ to $0$, all the Pieri factors $\rho_i^{(0)}$ and $\rho_{n-1}^{\pm(0)}$ are in $\Wafgr$. 

\subsection{Affine Grassmannian elements and partitions}
\label{SS:Grass and partitions}
We describe a parametrization of 
$\Wafgr$ in types $\mathrm{B}_n^{(1)}$ and $\mathrm{D}_n^{(1)}$ using partitions. 
Let $\Y^m_n$ be the set of partitions that have first part at most $m$ and each part of size at most $n$ occurs at most once.
For type $\B{n}^{(1)}$, define 
$$\Y^n_{\B{}}:=\Y_{n-1}^{2n-1}.$$

For \(r\geq 0\), define
\begin{equation}\label{eq_rho_i_l}
  \rho_\ell^{(r)}= \begin{cases} \rho_\ell^{(0)}, & \text{if \(r\) is even},\\ \rho_\ell^{(1)}, & \text{if \(r\) is odd}. \end{cases}  
\end{equation}
Define the map $w: \Y^n_{\B{}}\to \WBgr{n}$ by the following: Let $\lambda=(\lambda_1,\dots,\lambda_r)\in \Y^n_{\B{}}$. We define
\begin{align}\label{eq_w_lam}
w_\lambda:=\rho_{\lambda_r}^{(r-1)}\cdots\rho_{\lambda_2}^{(1)}\rho_{\lambda_1}^{(0)}.
\end{align}
This gives a bijection $\Y^n_{\B{}}\to \WBgr{n}$, see \cite[Theorem 28]{BM}.

The partitions corresponding to the affine Grassmannian elements of type \(\mathrm{D}_n^{(1)}\) are defined similarly, with one additional feature. If a partition \(\lambda\) contains parts of size $n-1$, then all such parts must carry the same sign. For such partitions $\lambda$, we define the sign of \(\lambda\) to be the common sign of such parts of size $n-1$. Consequently, whenever \(\lambda\) contains a part of size $n-1$, we distinguish two signed copies, denoted by \(\lambda^{+}\) and \(\lambda^{-}\), corresponding to the two types of twin Pieri factors. Let $\Y^n_{\D{}}$ be the set of partitions $\la$ with $\la_1 \le 2n-2$ 
with distinct parts of size smaller than $n-1$, where either the partition (a) has no parts of size $n-1$, (b) has parts of size $n-1$ and is marked $+$, or (c) has parts of size $n-1$ and is marked $-$.

Let ${\Y'}^{2n-2}_{n-2}$ be the set of those partitions $\lambda$ in $ \Y^{2n-2}_{n-2}$ having at least one part of size $n-1$ and define ${\Y}''^{2n-2}_{n-2}:={\Y}^{2n-2}_{n-2}\setminus {\Y}'^{2n-2}_{n-2}$. Then, 
\begin{align}
    \Y^n_{\D{}}={\Y}''^{2n-2}_{n-2}\sqcup{\Y'}^{2n-2,+}_{n-2}\sqcup{\Y'}^{2n-2,-}_{n-2}.
\end{align}

We define a bijection $w:\Y_{\D{}}^{n}\longrightarrow \WDgr{n}$ as follows. 
If \(\lambda\) contains no part equal to \(n-1\), 
then we define $w_\lambda$ 
just by \eqref{eq_w_lam}. 
We can also define \(\rho_{n-1}^{\pm(r)}\) by \(r\bmod2\), similarly using \eqref{eq_rho_i_l}. Now if \(\lambda\) contains some parts equal to \(n-1\), then the expressions for \(w_{\lambda^\pm}\) are obtained by replacing each factor \(\rho_{\lambda_i}^{(i-1)}\) corresponding to \(\lambda_i=n-1\) with \(\rho_{\lambda_i}^{\pm(i-1)}\), respectively. Note that if $\la$ contains a part of size $n-1$, then $w_{\la^+}$ and $w_{\la^-}$ are related by the affine Dynkin automorphism $\autn$ of \eqref{E:auto n}.

\begin{example}
For $\mathrm{B}_3^{(1)}$, $\rho_i^{(0)}$ is the horizontal row with $i$ boxes from the left in first row and $\rho_i^{(1)}$ is the horizontal row with $i$ boxes from the left in second row. Also in Figure \ref{F:B3poset} the affine Grassmannian left weak order poset is drawn for $\mathrm{B}_3^{(1)}$.

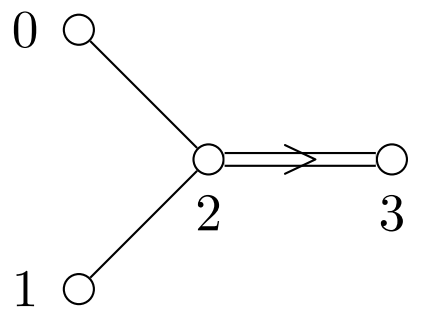
\begin{center}
\begin{tikzpicture}[node distance=1cm and 1cm, every node/.style={draw, circle, inner sep=2pt}]
    \node[label=below:{$2$}] (2) {};
    \node[label=left:{$0$}, above left=1cm of 2] (0) {};
    \node[label=left:{$1$}, below left=1cm of 2] (1) {};
    \node[label=below:{$3$}, right=1cm of 2] (3) {};
    \begin{scope}[every node/.style={inner sep=2pt}]
   \draw[double distance=2pt] (2) -- (3) node[midway, yshift=0ex] {$>$};
   \end{scope}
       \draw (1) -- (2)--(0);
 \end{tikzpicture}
 \end{center}

\[\begin{ytableau}
0&2&3&2&0\\1&2&3&2&1\\0&2&3&2&0\\1&2&3&2&1\\\vdots&\vdots&\vdots&\vdots&\vdots
\end{ytableau}\]

\ytableausetup{smalltableaux}

\begin{figure}[ht]
\centering
\resizebox{0.95\textwidth}{!}{%
\begin{tikzpicture}[
    scale=.8,
    tab/.style={inner sep=1pt, outer sep=2pt, align=center},
    edge/.style={-{Latex[length=2mm]}, thick},
    lab/.style={font=\scriptsize, fill=white, inner sep=1pt},
    rho/.style={font=\small, text=red}
]


\node[tab] (v0)  at (0,0) {$\varnothing$};

\node[tab] (v1)  at (0,-1.4) {
\begin{ytableau}
0
\end{ytableau}};
\node[rho] at (-0.5,-1.35) {$\rho_1^{(0)}$};

\node[tab] (v2)  at (0,-2.8) {
\begin{ytableau}
0&2
\end{ytableau}};
\node[rho] at (-0.8,-2.75) {$\rho_2^{(0)}$};

\node[tab] (v3)  at (-2.4,-4.3) {
\begin{ytableau}
0&2&3
\end{ytableau}};
\node[rho] at (-3.3,-4.25) {$\rho_3^{(0)}$};

\node[tab] (v21) at (2.4,-4.3) {
\begin{ytableau}
0&2\\
1
\end{ytableau}};

\node[tab] (v4)  at (-2.4,-5.9) {
\begin{ytableau}
0&2&3&2
\end{ytableau}};
\node[rho] at (-3.5,-5.85) {$\rho_4^{(0)}$};

\node[tab] (v31) at (2.4,-5.9) {
\begin{ytableau}
0&2&3\\
1
\end{ytableau}};

\node[tab] (v5) at (-5.6,-7.6) {
\begin{ytableau}
0&2&3&2&0
\end{ytableau}};
\node[rho] at (-6.9,-7.55) {$\rho_5^{(0)}$};

\node[tab] (v41) at (-2.4,-7.6) {
\begin{ytableau}
0&2&3&2\\
1
\end{ytableau}};

\node[tab] (v32) at (2.4,-7.6) {
\begin{ytableau}
0&2&3\\
1&2
\end{ytableau}};

\node[tab] (v51) at (-5.6,-9.5) {
\begin{ytableau}
0&2&3&2&0\\
1
\end{ytableau}};

\node[tab] (v421) at (-2.4,-9.5) {
\begin{ytableau}
0&2&3&2\\
1&2
\end{ytableau}};

\node[tab] (v33) at (2.4,-9.5) {
\begin{ytableau}
0&2&3\\
1&2&3
\end{ytableau}};

\node[tab] (v430) at (5.6,-9.5) {
\begin{ytableau}
0&2&3\\
1&2\\
0
\end{ytableau}};

\node[tab] (v52) at (-5.7,-11.5) {
\begin{ytableau}
0&2&3&2&0\\
1&2
\end{ytableau}};

\node[tab] (v420) at (2.4,-11.5) {
\begin{ytableau}
0&2&3&2\\
1&2\\
0
\end{ytableau}};

\node[tab] (v423) at (-2.4,-11.5) {
\begin{ytableau}
0&2&3&2\\
1&2&3
\end{ytableau}};

\node[tab] (v330) at (5.6,-11.5) {
\begin{ytableau}
0&2&3\\
1&2&3\\
0
\end{ytableau}};


\draw[edge] (v0) -- node[lab,right] {$0$} (v1);
\draw[edge] (v1) -- node[lab,right] {$2$} (v2);

\draw[edge] (v2) -- node[lab,left] {$3$} (v3);
\draw[edge] (v2) -- node[lab,right] {$1$} (v21);

\draw[edge] (v3) -- node[lab,left] {$2$} (v4);
\draw[edge] (v3) -- node[lab,above] {$1$} (v31);
\draw[edge] (v21) -- node[lab,right] {$3$} (v31);

\draw[edge] (v4) -- node[lab,left] {$0$} (v5);
\draw[edge] (v4) -- node[lab,right] {$1$} (v41);
\draw[edge] (v31) -- node[lab,right] {$2$} (v32);

\draw[edge] (v5) -- node[lab,left] {$1$} (v51);
\draw[edge] (v41) -- node[lab,left] {$0$} (v51);
\draw[edge] (v41) -- node[lab,right] {$2$} (v421);
\draw[edge] (v32) -- node[lab,left] {$1$} (v421);
\draw[edge] (v32) -- node[lab,right] {$3$} (v33);
\draw[edge] (v32) -- node[lab,above] {$0$} (v430);

\draw[edge] (v51) -- node[lab,left] {$2$} (v52);
\draw[edge] (v421) -- node[lab,pos=.70,right] {$0$} (v420);
\draw[edge] (v421) -- node[lab,right] {$3$} (v423);
\draw[edge] (v33) -- node[lab,pos=.70,right] {$0$} (v330);
\draw[edge] (v33) -- node[lab,pos=.75,right] {$1$} (v423);
\draw[edge] (v430) -- node[lab,above] {$3$} (v330);
\draw[edge] (v430) -- node[lab,pos=.72,below] {$1$} (v420);



\end{tikzpicture}
}
\caption{$\WBgr{3}$ left weak poset}
\label{F:B3poset}
\end{figure}

\end{example}

\begin{example}
Consider type $\mathrm{D}_4^{(1)}$. For $i\neq 3$, $\rho_i^{(0)}$ is the horizontal row with $i$ boxes in first row and $\rho_i^{(1)}$ is the horizontal row with $i$ boxes in second row. $\rho_3^{-(0)}=s_3s_2s_0$, $\rho_3^{+(0)}=s_4s_2s_0$, $\rho_3^{+(1)}=s_3s_2s_1$, $\rho_3^{-(1)}=s_4s_2s_1$.

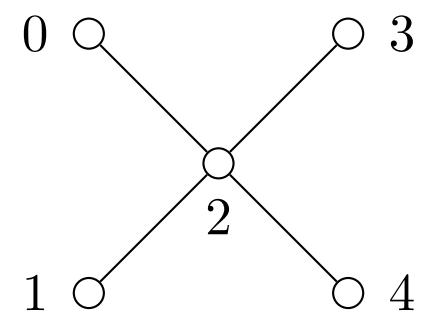
\ytableausetup{nosmalltableaux}
\begin{center}
\begin{tikzpicture}[node distance=1cm and 1cm, every node/.style={draw, circle, inner sep=2pt}]
    \node[label=below:{$2$}] (2) {};
    \node[label=left:{$0$}, above left=1cm of 2] (0) {};
    \node[label=left:{$1$}, below left=1cm of 2] (1) {};
    \node[label=right:{$3$}, above right=1cm of 2] (3) {};
    \node[label=right:{$4$}, below right=1cm of 2] (4) {};
       \draw (1)--(2)--(0);
       \draw (3)--(2)--(4);
 \end{tikzpicture}
 \end{center}

\[\begin{ytableau}
0&2&4&3&2&0\\1&2&3&4&2&1\\0&2&4&3&2&0\\1&2&3&4&2&1\\\vdots&\vdots&\vdots&\vdots&\vdots&\vdots
\end{ytableau}\]

\ytableausetup{smalltableaux}

\newcommand{\YT}[1]{\begin{ytableau}#1\end{ytableau}}

\begin{figure}[ht]
\centering
\resizebox{\textwidth}{!}{%
\begin{tikzpicture}[
    x=0.8cm,
    y=1cm,
    >=Stealth,
    arr/.style={-{Stealth[length=2mm]}, thin},
    tnode/.style={inner sep=1pt, outer sep=2pt},
    every node/.style={font=\scriptsize},
    lab/.style={font=\tiny, fill=white, fill opacity=.95, text opacity=1, inner sep=1.2pt},
    rho/.style={font=\scriptsize, text=red}
]


\node[tnode] (v0) at (0,8.8) {$\varnothing$};

\node[tnode] (v1) at (0,7.9) {$\YT{0}$};

\node[tnode] (v2) at (0,7.0) {$\YT{0&2}$};

\node[tnode] (v3p) at (-3.8,5.95) {$\color{Magenta}\YT{0&2&3}$};
\node[tnode] (v21)  at (0,5.95) {$\YT{0&2\\1}$};
\node[tnode] (v3m) at (3.8,5.95) {$\color{Plum}\YT{0&2&4}$};


\node[tnode] (v31p) at (-3.8,4.75) {$\color{Magenta}\YT{0&2&3\\1}$};
\node[tnode] (v4)   at (0,4.75) {$\YT{0&2&4&3}$};
\node[tnode] (v31m) at (3.8,4.75) {$\color{Plum}\YT{0&2&4\\1}$};


\node[tnode] (v32p) at (-3.8,3.45) {$\color{Magenta}\YT{0&2&3\\1&2}$};
\node[tnode] (v41)  at (0,3.45) {$\YT{0&2&4&3\\1}$};
\node[tnode] (v32m) at (3.8,3.45) {$\color{Plum}\YT{0&2&4\\1&2}$};
\node[tnode] (v5)   at (2.2,3.45) {$\YT{0&2&4&3&2}$};


\node[tnode] (v320p) at (-6,2) {$\color{Magenta}\YT{0&2&3\\1&2\\0}$};
\node[tnode] (v33p)  at (-3.8,2.15) {$\color{Magenta}\YT{0&2&3\\1&2&4}$};
\node[tnode] (v321)  at (0,2.15) {$\YT{0&2&4&3\\1&2}$};

\node[tnode] (v420m) at (6,2.15) {$\color{Plum}\YT{0&2&4\\1&2\\0}$};
\node[tnode] (v43m)  at (3.8,2.15) {$\color{Plum}\YT{0&2&4\\1&2&3}$};
\node[tnode] (v42)   at (-2.2,2.15) {$\YT{0&2&4&3&2\\1}$};
\node[tnode] (v6)    at (2.2,2.25) {$\YT{0&2&4&3&2&0}$};


\node[tnode] (b340)  at (-6,0.5) {$\color{Magenta}\YT{0&2&3\\1&2&4\\0}$};
\node[tnode] (b44)   at (-3.8,0.5) {$\color{Magenta}\YT{0&2&4&3\\1&2&4}$};

\node[tnode] (b4210) at (-2,0.5) {$\YT{0&2&4&3\\1&2\\0}$};
\node[tnode] (b522)  at (2.2,0.5) {$\YT{0&2&4&3&2\\1&2}$};

\node[tnode] (b430)  at (6,0.5) {$\color{Plum}\YT{0&2&4\\1&2&3\\0}$};
\node[tnode] (b433)  at (3.8,0.5) {$\color{Plum}\YT{0&2&4&3\\1&2&3}$};

\node[tnode] (b620)  at (0,0.5) {$\YT{0&2&4&3&2&0\\1}$};


\draw[arr] (v0) -- node[lab,pos=.55,right] {$0$} (v1);
\draw[arr] (v1) -- node[lab,pos=.55,right] {$2$} (v2);

\draw[arr] (v2) -- node[lab,pos=.35,above left] {$3$} (v3p);
\draw[arr] (v2) -- node[lab,pos=.50,right] {$1$} (v21);
\draw[arr] (v2) -- node[lab,pos=.35,above right] {$4$} (v3m);

\draw[arr] (v3p) -- node[lab,pos=.35,left] {$1$} (v31p);
\draw[arr] (v3p) -- node[lab,pos=.25,above] {$4$} (v4);

\draw[arr] (v21) -- node[lab,pos=.70,above left] {$3$} (v31p);
\draw[arr] (v21) -- node[lab,pos=.35,above right] {$4$} (v31m);

\draw[arr] (v3m) -- node[lab,pos=.25,above] {$3$} (v4);
\draw[arr] (v3m) -- node[lab,pos=.35,right] {$1$} (v31m);

\draw[arr] (v31p) -- node[lab,pos=.35,left] {$2$} (v32p);
\draw[arr] (v31p) -- node[lab,pos=.28,above] {$4$} (v41);

\draw[arr] (v4) -- node[lab,pos=.35,above left] {$1$} (v41);
\draw[arr] (v4) -- node[lab,pos=.35,above right] {$2$} (v5);

\draw[arr] (v31m) -- node[lab,pos=.78,below left] {$3$} (v41);
\draw[arr] (v31m) -- node[lab,pos=.35,right] {$2$} (v32m);

\draw[arr] (v32p) -- node[lab,pos=.45,left] {$0$} (v320p);
\draw[arr] (v32p) -- node[lab,pos=.35,right] {$4$} (v33p);

\draw[arr] (v41) -- node[lab,pos=.35,left] {$2$} (v321);

\draw[arr] (v32m) -- node[lab,pos=.35,left] {$0$} (v420m);
\draw[arr] (v32m) -- node[lab,pos=.35,right] {$3$} (v43m);

\draw[arr] (v5) -- node[lab,pos=.32,left] {$1$} (v42);
\draw[arr] (v5) -- node[lab,pos=.40,above] {$0$} (v6);

\draw[arr] (v320p) -- node[lab,pos=.40,left] {$4$} (b340);
\draw[arr] (v33p)  -- node[lab,pos=.40,left] {$0$} (b340);
\draw[arr] (v33p)  -- node[lab,pos=.42,right] {$2$} (b44);
\draw[arr] (v321) to[bend right=8]  node[lab,pos=.6,above left] {$4$} (b44);
\draw[arr] (v321)  to[bend left=8]  node[lab,pos=.39,below left] {$0$} (b4210);
\draw[arr] (v321) to[bend right=8]  node[lab,pos=.4,above left] {$1$} (b522);

\draw[arr] (v321) to[bend left=8]  node[lab,pos=.7,above right] {$3$} (b433);

\draw[arr] (v420m) -- node[lab,pos=.40,right] {$3$} (b430);

\draw[arr] (v43m)  -- node[lab,pos=.40,right] {$0$} (b430);
\draw[arr] (v43m)  -- node[lab,pos=.40,right] {$2$} (b433);


\draw[arr] (v42) -- node[lab,pos=.35,right] {$2$} (b522);
\draw[arr] (v42) -- node[lab,pos=.35,left] {$0$} (b620);

\draw[arr] (v6) -- node[lab,pos=.45,right] {$1$} (b620);

\node[rho, right=-5pt of v1]
    {$\rho_1^{(0)}$};

\node[rho, right=-5pt of v2,yshift=3pt]
    {$\rho_2^{(0)}$};

\node[rho, left=-5pt of v3p]
    {$\rho_3^{-(0)}$};

\node[rho, right=-5pt of v3m]
    {$\rho_3^{+(0)}$};

\node[rho, right=-4pt of v4]
    {$\rho_4^{(0)}$};

\node[rho, above=-4pt of v5, xshift=6pt]
    {$\rho_5^{(0)}$};

\node[rho, above=-4pt of v6, xshift=9pt]
    {$\rho_6^{(0)}$};

\end{tikzpicture}%
}
\caption{$\WDgr{4}$ left weak poset}
\end{figure}

\end{example}

{
\section{Anti-dominant translation elements}
\label{S:translation-factorization}

We determine the affine Grassmannian elements corresponding to the
anti-dominant fundamental translations.  Recall from
\eqref{E:kappa special} that, for a special node $i\in I_\af^s$,
\begin{equation}
\label{E:kappa-special-recalled}
  \kappa_i=\pi_i u_i\pi_i^{-1}
  \qquad\text{and hence}\qquad
  t_{-\varpi_i^\lor}=\pi_i^{-1}\kappa_i.
\end{equation}
For type $\mathrm{B}_n^{(1)}$, we have $I_\af^s=\{0,1\}$ and
$\pi_1=\auto$; for type $\mathrm{D}_n^{(1)}$, we have
$I_\af^s=\{0,1,n-1,n\}$ and $\pi_1=\auto\autn$.
We use the notation for the Pieri factors and colored partitions introduced
in \S \ref{SS:Grass and partitions}.

\begin{lem} 
In types $\mathrm{B}_n^{(1)}$ and $\mathrm{D}_n^{(1)}$, if $i\notin I_\af^s$ then
\begin{align}\label{E:kappa nonspecial}
\kappa_i = 
\begin{cases}
t_{-\varpi_i^\lor} & \text{if $i$ is even} \\
\pi_1 t_{-\varpi_i^\lor} = u_1 t_{-\varpi_1^\vee-\varpi_i^\vee}& \text{if $i$ is odd.}
\end{cases}
\end{align}
\end{lem}
\begin{proof} This is an application of Lemma \ref{L:kappa} using the fact that if $i\in I$ is not special
then the index $i^s$ of Lemma \ref{L:kappa} is the unique $i^s\in \{0,1\}$ satisfying $i\equiv i^s \mod 2$.
For the equality in the $i$ odd case, by \eqref{eq_act_pi} $\pi_1 = t_{\varpi_1^\vee} u_1= u_1 t_{-\varpi_1^\vee}$.
\end{proof}

\subsection{Type $\mathrm{B}_n^{(1)}$}

For $1\leq i\leq n$, set $z_1=e$ and for $i\ge 2$, set
\begin{equation}
\label{E:zi}
  z_i=s_{i-1}\dotsm s_2s_1.
\end{equation}

\begin{lem}
\label{L:B-even-rectangles}
For $1\leq i\leq n$ and $0\leq 2k\leq i$, one has
\begin{equation}
\label{E:B-even-rectangles}
  w_{(2n-i)^{2k}}=z_i^{2k}t_{-\varpi_{2k}^\lor}.
\end{equation}
\end{lem}

\begin{proof}
The assertion is clear for $k=0$.  For $k=1$, using
$s_0=s_\theta t_{-\theta^\lor}$ and
$\theta^\lor=\varpi_2^\lor$, we obtain
\begin{align*}
  w_{(2n-i)^2}
  &=(s_i\dotsm s_ns_{n-1}\dotsm s_1)
    (s_i\dotsm s_ns_{n-1}\dotsm s_2s_\theta)
    t_{-\varpi_2^\lor} \\
  &=z_i^2t_{-\varpi_2^\lor}.
\end{align*}
The second equality may be checked on the coordinate vector
$x=(x_1,\ldots,x_n)$: the finite Weyl group factor sends it to
\[
  (x_3,\ldots,x_i,x_1,x_2,x_{i+1},\ldots,x_n)=z_i^2x.
\]
Assume that \eqref{E:B-even-rectangles} holds for $k-1$.  Since the
partition factorization is length-additive,
\begin{align*}
  w_{(2n-i)^{2k}}
  &=w_{(2n-i)^2}w_{(2n-i)^{2k-2}} \\
  &=z_i^2t_{-\varpi_2^\lor}
    z_i^{2k-2}t_{-\varpi_{2k-2}^\lor} \\
  &=z_i^{2k}
    t_{z_i^{2-2k}(-\varpi_2^\lor)-\varpi_{2k-2}^\lor}
   =z_i^{2k}t_{-\varpi_{2k}^\lor},
\end{align*}
as required.
\end{proof}

\begin{prop}
\label{P:B-kappa}
For $1\leq i\leq n$, one has
\begin{equation}
\label{E:B-antidominant-translation}
  t_{-\varpi_i^\lor}=\pi_1^i\kappa_i,
  \qquad
  \kappa_i=w_{(2n-i)^i}.
\end{equation}
In particular,
\begin{equation}
\label{E:u1-B}
  u_1=\rho_{2n-1}^{(1)},
  \qquad
  \kappa_1=\rho_{2n-1}^{(0)}=w_{(2n-1)}.
\end{equation}
\end{prop}

\begin{proof}
The identity for $u_1$ in \eqref{E:u1-B} follows directly from its reduced
expression; the identity for $\kappa_1$ is the case $i=1$ of
\eqref{E:B-antidominant-translation}.

The first identity follows from the general formula for $\kappa_i$: if $i$
is even, then $\varpi_i^\lor\in Q^\lor$ and
$\kappa_i=t_{-\varpi_i^\lor}$; if $i$ is odd, then
$\varpi_i^\lor\equiv\varpi_1^\lor\pmod{Q^\lor}$ and
$\kappa_i=\pi_1t_{-\varpi_i^\lor}$.  Since $\pi_1^2=1$, both cases are
equivalent to the first identity in \eqref{E:B-antidominant-translation}.

It remains to identify $\kappa_i$.  If $i=2k$ is even, Lemma
\ref{L:B-even-rectangles} and $z_i^i=1$ give
$w_{(2n-i)^i}=t_{-\varpi_i^\lor}=\kappa_i$.
Suppose that $i=2k+1$.  Then $z_i^{2k}=z_i^{-1}$, and Lemma
\ref{L:B-even-rectangles} gives
\begin{align*}
  \pi_1w_{(2n-i)^i}
  &=t_{\varepsilon_1}u_1w_{(2n-i)}z_i^{-1}
    t_{-\varpi_{2k}^\lor}.
\end{align*}
A direct calculation in the finite Weyl group gives
\begin{equation}
\label{E:B-finite-calculation}
  u_1s_i\dotsm s_ns_{n-1}\dotsm s_2s_\theta=z_i.
\end{equation}
Using
$w_{(2n-i)}=s_i\dotsm s_ns_{n-1}\dotsm s_2s_\theta
t_{-\varpi_2^\lor}$, equation \eqref{E:B-finite-calculation} yields
\[
  t_{\varepsilon_1}u_1w_{(2n-i)}z_i^{-1}=t_{-\varepsilon_i}.
\]
Hence
\[
  \pi_1w_{(2n-i)^i}
   =t_{-\varepsilon_i-\varpi_{i-1}^\lor}
   =t_{-\varpi_i^\lor}
   =\pi_1\kappa_i,
\]
and therefore $w_{(2n-i)^i}=\kappa_i$.
\end{proof}

\subsection{Type $\mathrm{D}_n^{(1)}$}

\begin{prop}
\label{P:D-kappa}
In type $\mathrm{D}_n^{(1)}$, one has
\begin{align}
\label{E:D-antidominant-nonspin}
  t_{-\varpi_i^\lor}&=\pi_1^i\kappa_i
  &&\text{for $1\leq i\leq n-2$,}
  \\
\label{E:D-antidominant-spin}
  t_{-\varpi_{n-1}^\lor}&=\pi_{n-1}^{-1}\kappa_{n-1},
  &
  t_{-\varpi_n^\lor}&=\pi_n^{-1}\kappa_n.
\end{align}
Moreover,
\begin{equation}
\label{E:D-kappa}
  \kappa_i=
  \begin{cases}
    w_{(2n-1-i)^i},&1\leq i\leq n-2,\\
    w_{(n-1,n-2,\ldots,2,1)^+},&i=n-1,\\
    w_{(n-1,n-2,\ldots,2,1)^-},&i=n.
  \end{cases}
\end{equation}
In particular,
\begin{equation}
\label{E:u1-D}
  u_1=\rho_{2n-2}^{(1)},
  \qquad
  \kappa_1=\rho_{2n-2}^{(0)}=w_{(2n-2)}.
\end{equation}
\end{prop}
\begin{proof}
The identity for $u_1$ in \eqref{E:u1-D} follows directly from its reduced
expression; the identity for $\kappa_1$ is the case $i=1$ of
\eqref{E:D-kappa}.

Equations \eqref{E:D-antidominant-nonspin} and
\eqref{E:D-antidominant-spin} follow from the general formula for
$\kappa_i$ and from \eqref{E:kappa-special-recalled}.

For $1\leq i\leq n-2$, the same finite-Weyl-group calculation as in Lemma
\ref{L:B-even-rectangles}, with the type $\mathrm{D}_n$ Pieri factor
$\rho_{2n-1-i}$ in place of $\rho_{2n-i}$, gives
\[
  w_{(2n-1-i)^i}=\kappa_i.
\]
More explicitly, for even $i$ the finite Weyl group factor cancels after
$i$ repetitions and the product is $t_{-\varpi_i^\lor}$; for odd $i$ one
obtains $\pi_1t_{-\varpi_i^\lor}$.  These are precisely the two cases of
\eqref{E:D-antidominant-nonspin}.

It remains to prove the two spin cases.  We prove the formula for $\kappa_n$;
the formula for $\kappa_{n-1}$ follows by applying $\autn$.  Suppose first
that $n$ is even.  Put $h=1/2$.  Since
$w_0(\varpi_n^\lor)=-\varpi_n^\lor$, successively moving the negative
coordinates of $(-h)^n$ to the right gives the reduced expression
\begin{equation}
\label{E:un-even}
  u_n=\dotsm
  (s_4\dotsm s_{n-2}s_{n-1})
  (s_3\dotsm s_{n-2}s_n)
  (s_2\dotsm s_{n-2}s_{n-1})
  (s_1s_2\dotsm s_{n-2}s_n).
\end{equation}
Conjugation by $\pi_n$ sends $s_i$ to $s_{n-i}$ for $0\leq i\leq n$;
hence
\begin{align*}
  \kappa_n
  &=\pi_nu_n\pi_n^{-1} \\
  &=\dotsm
    (s_{n-4}\dotsm s_2s_1)
    (s_{n-3}\dotsm s_2s_0)
    (s_{n-2}\dotsm s_2s_1)
    (s_{n-1}\dotsm s_2s_0) \\
  &=\dotsm
    \rho_{n-4}^{(1)}\rho_{n-3}^{(0)}
    \rho_{n-2}^{(1)}\rho_{n-1}^{-(0)}
   =w_{(n-1,n-2,\ldots,2,1)^-}.
\end{align*}

Suppose now that $n$ is odd.  Then
$w_0(\varpi_n^\lor)=-\varpi_{n-1}^\lor$, and the same calculation gives
\begin{equation}
\label{E:un-odd}
  u_n=\dotsm
  (s_5\dotsm s_{n-1})
  (s_4\dotsm s_{n-2}s_n)
  (s_3\dotsm s_{n-1})
  (s_2\dotsm s_{n-2}s_n)
  (s_1\dotsm s_{n-1}).
\end{equation}
Here $\pi_n$ acts on the special nodes by the cycle $(0,n,1,n-1)$ and
sends $i$ to $n-i$ for $2\leq i\leq n-2$.  It follows that
\begin{align*}
  \kappa_n
  &=\pi_nu_n\pi_n^{-1} \\
  &=\dotsm
    (s_{n-5}\dotsm s_2s_0)
    (s_{n-4}\dotsm s_2s_1)
    (s_{n-3}\dotsm s_2s_0)
    (s_{n-2}\dotsm s_2s_1)
    (s_{n-1}\dotsm s_2s_0) \\
  &=\dotsm
    \rho_{n-5}^{(0)}\rho_{n-4}^{(1)}
    \rho_{n-3}^{(0)}\rho_{n-2}^{(1)}\rho_{n-1}^{-(0)}
   =w_{(n-1,n-2,\ldots,2,1)^-}.
\end{align*}
\end{proof}

}

\begin{example}
In type $\mathrm{B}_3^{(1)}$, we abbreviate the product
$s_{i_1}s_{i_2}\cdots s_{i_L}$ by $s_{i_1i_2\cdots i_L}$.
The elements $\kappa_i$ are given by
\begin{align*}
\kappa_1= s_{02320}
     = w_{(5)},~
\kappa_2= s_{2321\cdot 2320}
     = w_{(4,4)},~\kappa_3= s_{320\cdot 321\cdot 320}
     = w_{(3,3,3)}.
\end{align*}
The translations associated with the negative fundamental coweights
factor as
\begin{align*}
t_{-\varpi_1^\lor}= \pi_1\kappa_1,~
t_{-\varpi_2^\lor}= \kappa_2,~
t_{-\varpi_3^\lor} = \pi_1\kappa_3.
\end{align*}
\end{example}

\begin{example}
In type $\mathrm{D}_5^{(1)}$, the elements $\kappa_i$ are given by
\begin{align*}
\kappa_1
    &= s_{02354320}
     = \rho_8^{(0)},\\
\kappa_2
    &= s_{2354321\cdot 2354320}
     =\rho_7^{(1)}
       \rho_7^{(0)}= w_{(7,7)},\\
\kappa_3
    &= s_{354320\cdot 354321\cdot 354320}
     =\rho_6^{(0)}\rho_6^{(1)}
       \rho_6^{(0)}= w_{(6,6,6)},\\
\kappa_4
    &= s_{1\cdot 20\cdot 321\cdot 5320}
     = \rho_1^{(1)}\rho_2^{(0)}
       \rho_3^{(1)}\rho_4^{+(0)}
     = w_{(4,3,2,1)^+},\\
\kappa_5
    &= s_{1\cdot 20\cdot 321\cdot 4320}
     = \rho_1^{(1)}\rho_2^{(0)}
       \rho_3^{(1)}\rho_4^{-(0)}
     = w_{(4,3,2,1)^-}.
\end{align*}
The translations associated with the negative fundamental coweights
factor as
\begin{align*}
t_{-\varpi_1^\lor} &= \pi_1\kappa_1,
&
t_{-\varpi_2^\lor} &= \kappa_2,
&
t_{-\varpi_3^\lor} &= \pi_1\kappa_3,\\
t_{-\varpi_4^\lor} &= \pi_5\kappa_4,
&
t_{-\varpi_5^\lor} &= \pi_4\kappa_5.
\end{align*}
\end{example}

\section{Symmetric function bases for the equivariant $K$-homology of affine Grassmannians in types $\B{}$ and $\D{}$}\label{S:sym_fun_rea}
We introduce a pair of Schubert bases for each of the types $\mathrm{B}_n^{(1)}$ and $\mathrm{D}_n^{(1)}$, for the equivariant $K$-homology of the affine Grassmannians, namely, the basis of structure sheaves of Schubert varieties and the basis of ideal sheaves of boundaries of Schubert varieties in the ind affine Grassmannians.
This is achieved by transporting the equivariant Schubert bases $k_w$ and $l_w$ from $\LG$ to symmetric functions.


\subsection{Realization of extended Peterson algebra by equivariant symmetric functions}
We define notation for the $K$-theoretic group law: $$x\oplus y = x + y - xy,\quad
x\ominus y = (x-y)/(1-y), \quad \bar{y} = \ominus(y) = -y/(1-y).$$ 
Recalling the notation $b_i$ from \S \ref{SS:introBD}, we have
$\bar{b}_i  = 1 - e^{a_i}$ for $1\le i\le n$.

For types $\mathrm{B}_n^{(1)}$ and $\mathrm{D}_n^{(1)}$ we have (see \eqref{E:TDc})
\begin{align*}
    c(\alpha_0) &= b_1 \oplus b_2 \\
    c(\alpha_i) &= b_{i+1} \ominus b_i \qquad\text{for $1\le i\le n-1$} \\
    c(\alpha_n) &= \begin{cases}
        \bar{b}_n & \text{for $\mathrm{B}_n^{(1)}$} \\
        \bar{b}_{n-1} \oplus \bar{b}_n &  \text{for $\mathrm{D}_n^{(1)}$.}
    \end{cases}
\end{align*}
Let $\hL_{\sR}$ be the completion of the algebra of symmetric functions in variables $(y_1,y_2,\dotsc)$ with coefficients in $\sR$.



In both types $\mathrm{B}_n^{(1)}$ and $\mathrm{D}_n^{(1)}$, the natural $W$-action on $\sR$ extends to an action of the affine Weyl group $\Waf$ on $\hL_{\sR}$ by letting $W$ act on coefficients and translation elements $t_\la$ act by multiplication by $\tf{\la}$:
\begin{align}\label{eq_t_lam_act}
   t_\la f(y| b) = \tf{\la} f(y| b), \quad \text{for $f(y| b) \in \hL_{\sR}$}. 
\end{align}

In particular, for both types $\mathrm{B}_n^{(1)}$ and $\mathrm{D}_n^{(1)}$ using \eqref{E:highest root and associated coroot} the action of the simple reflection $s_0=t_{\theta^\lor}s_\theta$ is given by
\begin{align*}
(s_0 f)(y| b) &= (t_{\theta^\lor}s_\theta )f(y| b) 
=\Omega(b_1 |  y)\Omega(b_2 |  y) f(y| s_\theta b)\\
&= \Omega(b_1 |  y)\Omega(b_2 |  y) f(y|  \bar{b}_2, \bar{b}_1,b_3, \dotsc, b_n).
\end{align*}
The action of the simple reflection $s_i$ for $1\leq i\leq n$ is given by
\begin{align*}
(s_i f)(y| b_1,\dots,b_n)& = f(y|  b_1,b_2,\dots,b_{i+1},b_i, \dots, b_n)\quad (1\le i\le n-1),\\
(s_n f)(y| b_1,\dots,b_n) &=\begin{cases}
    f(y|  b_1,b_2, \dots,b_{n-1}, \bar{b}_n) &  \text{for $\mathrm{B}_n^{(1)}$}\\
     f(y|  b_1,b_2, \dots,b_{n-2}, \bar{b}_n,\bar{b}_{n-1}) & \text{for $\mathrm{D}_n^{(1)}$}.
\end{cases} 
\end{align*}
This \(\Waf\)-action extends further to an action of the extended affine Weyl group $\Wafex$ on $\hL_{\sR}$. Indeed, using the decomposition $\Wx\cong \Sigma {\ltimes} \Waf$ every element of $\Wx$ can be written as $\sigma w$ for $\sigma\in\Sigma$ and $w\in\Waf$.
The action of $\sigma\in \Sigma$ on $\hL_{\sR}$ is defined using \eqref{eq_act_pi} together with the translation action as in \eqref{eq_t_lam_act}.

We recall that the fundamental group is $\Sigma=\{\pi_0,\pi_1\}$ for type $\B{n}$ and is $\Sigma=\{\pi_0,\pi_1,\pi_{n-1},\pi_n\}$ for type $\D{n}$. The action of the fundamental group $\Sigma$ on {$\hL_{\sR}$} is given explicitly as follows:

 Let $f(y | b) \in \hL_{\sR}$. In type $\B{n}$,
\[
\pi_{1} f(y | b)
   = \tf{\varpi_1^\vee}(y | b)\, f\bigl(y | u_{1}(b)\bigr).\]
In type $\D{n}$, for $i\in \{1,n-1,n\}$ we have
\begin{align*}
\pi_{i} f(y | b)
   &= \tf{\varpi_i^\vee}(y | b)\, f\bigl(y | u_{i}(b)\bigr),
\end{align*}

\begin{prop}[{see \cite[Proposition 3.1]{ISY}}]\label{prop_kg_mod}
The $\Wafex$-action together with the 
left multiplication of $\RT$, extends to 
a left $\tKG$-module structure on $\hL_\sR$ .    
\end{prop}

 We recall the operator $\hat{T}_i$ and $\hat{D}_i$ on $\hL_{\sR}$ for $i\in I_\af$ from \eqref{E:TDc}. For $w\in W_\af$, the operators $\hat{T}_w$ and $\hat{D}_w$ are well-defined as $\hat{T}_i$ and $\hat{D}_i$ satisfy the braid relations. Now consider an element $\sigma w\in\Wx$ for $\sigma\in\Sigma$ and $w\in\Waf$, the operators  $\hat{T}_{\sigma w} = \sigma \hat{T}_w$ are well-defined on $\hL_{\sR}$.



For $w\in \WBexgr{n}$ define
$$gq^{(n)}_w(y| b) = \hat{T}_w(1_{\hL_{\sR}}), \quad \widetilde{gq}^{(n)}_w(y| b) = \hat{D}_w(1_{\hL_{\sR}}).$$
 We call $gq^{(n)}_w(y| b)$ the affine double $gq$-function and $\widetilde{gq}^{(n)}_w(y| b)$ the closed affine double $gq$-function.

For $w\in \WDexgr{n}$ define
$$\gd{w}{n}(y| b) = \hat{T}_w(1_{\hL_{\sR}}), \quad \tgd{w}{n}(y| b) = \hat{D}_w(1_{\hL_{\sR}}).$$

 We call $gd^{(n)}_w(y| b)$ the \emph{affine double $gd$-function} and $\widetilde{gd}^{(n)}_w(y| b)$ the \emph{closed affine double $gd$-function}.

Define the map $$\Psi:\tLG\to \hL_{\sR} \text{ by }\Psi(c) = c\cdot 1_{\hL_{\sR}}.$$ Let $\gGtil_{G}$ denote the $\RT$-span of $\Psi(k_w)$ for $w\in \Wafexgr$ of the corresponding types. Also consider the restriction map $\Psi|_{\LG}: \LG \to \hL_{\sR}$ and let $\gGhat_{G}$ denote the $\RT$-span of $\Psi(k_w)$ for $w\in \Wafgr$ of the corresponding types.

\begin{thm} \label{T:Peterson to symm funcs}
The map $\Psi$ restricts to an isomorphism $\tLG \cong \gGtil_{G}$ of $R(T)$-algebras and left $\tKG$-modules, 
such that
\begin{align}
    \Psi(k_w) &= \gq{w}{n}(y| b), & \Psi(\ell_w) &= \tgq{w}{n}(y| b)&&\text{for $w\in \WBexgr{n}$,} \label{eq:def_gq}\\
    \Psi(k_w) &= \gd{w}{n}(y| b), &    \Psi(\ell_w)&= \tgd{w}{n}(y| b) &&\text{for $w\in \WDexgr{n}$.}
   \label{eq:def_gd}
\end{align}
After extending scalars, it induces an isomorphism
\[
\Rsc\otimes_{\RT}\LG
\;\cong\;
\Rsc\otimes_{\RT}\gGhat_G
\]
of $\Rsc$-algebras.\end{thm}
\begin{proof}
We only prove this for type $\D{n}$. Type $\B{n}$ follows similarly. 

We show that the morphism $\Psi$ is injective. 
Let $F=\mathrm{Frac}(\sR)$ be the fraction field of $\sR$.
Let $F[P^\lor]$ be the group ring of $P^\lor$ considered as 
a subalgebra of $F[t_{\pm \frac{1}{2}\eps_1},\ldots,t_{\pm \frac{1}{2}\eps_n}]$.
Consider the homomorphism of $F$-algebras defined by 
\begin{equation}
F[t_{\pm \frac{1}{2}\eps_1},\ldots,t_{\pm \frac{1}{2}\eps_n}]
\rightarrow F\otimes_{\sR}\hL_{\sR},\quad 
t_{\frac{1}{2}\eps_i} \mapsto 
\Omega(b_i| y)^{\frac{1}{2}},\label{eq:Laurent_map}
\end{equation}
which is well-defined because  
$\Omega(b_i| y)^{\frac{1}{2}}$ is an invertible
element in $\hL_{\sR}$. 

Put $z_i=2-2\cosh(a_i)$ and
$t_j=y_j(1-y_j)$. Since
\[
\frac{\partial}{\partial y_j}
\Omega(b_i| y)^{1/2}
=
\frac{\Omega(b_i| y)^{1/2}\sinh(a_i)}
{1-z_it_j},
\]
the Cauchy determinant identity gives
\[
\det\left(
\frac{\partial\Omega(b_i| y)^{1/2}}
{\partial y_j}
\right)_{1\le i,j\le n}
=
\frac{
\displaystyle
\prod_{i=1}^n
\Omega(b_i| y)^{1/2}\sinh(a_i)
\prod_{1\le i<j\le n}
(z_i-z_j)(t_i-t_j)
}{
\displaystyle
\prod_{i,j=1}^n(1-z_it_j)
}.
\]
Since this determinant is nonzero, the Jacobian criterion for algebraic independence (see, e.g., \cite{Matsumura}) implies that the functions \(\Omega(b_i| y)^{1/2}\), \(1\leq i\leq n\), are algebraically independent over \(F\). Hence the map (5.4) is injective.

Restricting it to
\(F[P^\lor]\) gives an embedding
\[
F[P^\lor]\hookrightarrow F\otimes_{\sR}\hL_{\sR}.
\]
By \cite[Lemma~5.2]{LSS:K} for \(\LG\)
(cf.\ \cite[Lemma~2.14]{ISY} for \(\tLG\)), there is an
\(F\)-algebra isomorphism
\[
F\otimes_{\RT}\tLG
=
\bigoplus_{\lambda\in P^\lor}Ft_\lambda
\xrightarrow{\sim}
F[P^\lor].
\]
Composing this isomorphism with the embedding above, we obtain an
injective \(F\)-algebra homomorphism
\[
\Psi^{\mathrm{rat}}\colon
F\otimes_{\RT}\tLG
\xrightarrow{\sim}
F[P^\lor]
\hookrightarrow
F\otimes_{\sR}\hL_{\sR}.
\]
Indeed, both maps send $t_\lambda$ to
$\eta_\lambda$ and are compatible with
the action of the extended affine nil-Hecke
algebra. Hence their restrictions to
$\tLG$ agree.

Since $t_\lambda\cdot 1=\eta_\lambda$ by~\eqref{eq_t_lam_act}, we have $\Psi_{\mathrm{rat}}(1\otimes c)=1\otimes(c\cdot 1)=1\otimes\Psi(c)$ for every $c\in\widetilde L_G$, so $\Psi_{\mathrm{rat}}$ restricts to $\Psi$ on $\tLG$. Thus $\Psi$ is injective.

By \eqref{E:kw} and \eqref{E:lw}, we have
\eqref{eq:def_gq} and \eqref{eq:def_gd}. 
The proofs of the remaining statements
are straightforward (cf. \cite[Theorem 1.1]{ISY}).
\end{proof}

\begin{prop}
If $i$ is odd then
   $$\widetilde{gq}^{(n)}_{t_{-\varpi_i^\lor}}(y| b)=\tf{-\varpi_1^\lor}(y| b)\widetilde{gq}^{(n)}_{\kappa_i}(y| b)=\tf{\varpi_1^\lor}(y| b)\widetilde{gq}^{(n)}_{\kappa_i}(y| u_1(b)).$$ 
\end{prop}

\begin{prop}
If $i$ is odd and $1\leq i\leq n-2$ then
   $$\widetilde{gd}^{(n)}_{t_{-\varpi_i^\lor}}(y| b)=\tf{-\varpi_1^\lor}(y| b)\widetilde{gd}^{(n)}_{\kappa_i}(y| b)=\tf{\varpi_1^\lor}(y| b)\widetilde{gd}^{(n)}_{\kappa_i}(y| u_1(b)).$$
   For $i\in \{n-1, n\}$ we have 
   $$\widetilde{gd}^{(n)}_{t_{-\varpi_i^\lor}}(y| b)=\tf{-\varpi_i^\lor}(y| b)\widetilde{gd}^{(n)}_{\kappa_i}(y| b)=\tf{\varpi_{i^*}^\lor}(y| b)\widetilde{gd}^{(n)}_{\kappa_{i}}(y| u_{i}^{-1}(b))
   $$
   where $i^*\in \{n-1,n\}$ is such that $\pi_{i^*}=\pi_i^{-1}$.
\end{prop}

\subsection{A type $\mathrm{D}_n^{(1)}$ Dynkin symmetry}\label{subsec_dyn_sym}

Recall the affine Dynkin automorphism $\autn$ of $\mathrm{D}_n^{(1)}$ from \S \ref{SS:Dynkin}.
We have $\autn(b_i)=b_i$ for $1\le i\le n-1$ and $\autn(b_n)=\bar{b}_n$.

\begin{remark} If $\la\in \Y^n_{\D{}}$ has a part of size $n-1$ then $\autn(\la)$ changes the marking $\pm \to \mp$. Otherwise $\autn(\la)=\la$.
\end{remark}

The following Theorem holds by applying $\autn$.

\begin{thm}\label{T:gd auto}
Let $\lambda\in \Y^n_{\D{}}$ be a signed partition. Then 
\begin{align}
\gd{w_{\la}}{n}(y| b_1,\dotsc,b_{n-1},b_n) &= 
  \gd{\autn(w_{\la})}{n}(y| b_1,\dotsc,b_{n-1},\bar{b}_n), \\
  \tgd{w_{\la}}{n}(y| b_1,\dotsc,b_{n-1},b_n) &= 
  \tgd{\autn(w_{\la})}{n}(y| b_1,\dotsc,b_{n-1},\bar{b}_n).
\end{align}
In particular, we have
\begin{align}
\gd{w_{\la}}{n}(y| 0) &= 
  \gd{\autn(w_{\la})}{n}(y| 0), \label{eq:non-eq:gd1}\\
  \tgd{w_{\la}}{n}(y| 0) &= 
  \tgd{\autn(w_{\la})}{n}(y| 0).\label{eq:non-eq:gd2}
\end{align}
\end{thm}

\subsection{The non-equivariant $K$-homology ring $K_*(\Gr_{\Gsc})$}
\label{sec:non-eq}
Let $\epsilon\colon\RT=\Z[Q]\to\Z$ be the augmentation
homomorphism defined by $e^\lambda\mapsto 1$ for $\lambda\in Q$.
Its kernel is the augmentation ideal
\[
J:=\Ker(\epsilon)
  =(1-e^{-\alpha_i}\mid i\in I).
\]
Let \(\cO{w}^0\) denote the non-equivariant structure-sheaf Schubert classes. The Schubert-cell decomposition gives
\[
K_*^{T_{\mathrm{sc}}}(\Gr_{\Gsc})
=\bigoplus_{w\in W_{\mathrm{af}}^0}R(\Tsc)\,\mathscr O_w ~;
\qquad
K_*(\Gr_{\Gsc})
=\bigoplus_{w\in W_{\mathrm{af}}^0}\mathbb Z\,\mathscr O_w^0.
\]

The non-equivariant $K$-homology ring can be obtained by forgetting
\[
K_*^{\Tsc}(\Gr_{\Gsc})
\longrightarrow
\Z\otimes_{\Rsc}K_*^{\Tsc}(\Gr_{\Gsc})
\cong K_*(\Gr_{\Gsc}),
\]
where $\Z$ is regarded as an $\Rsc$-algebra via the augmentation
$e^\lambda\mapsto 1$ for $\lambda\in P$.
\begin{remark}
The specialization at $b=0$ is induced by the homomorphism
$\epsilon_{\R}\colon\R\to\Q$ defined by $e^{a_i}\mapsto 1$.
It gives a map $\hL_{\R}\to\hL_{\Q}$, whose restriction
to $\RT$ agrees with the augmentation
$\epsilon\colon\RT\to\Z$ followed by $\Z\hookrightarrow\Q$.
Equations \eqref{eq:non-eq:gd1} and \eqref{eq:non-eq:gd2} result from specializing the ambient symmetric functions at $b=0$. {This specialization should not be confused} with the base change
$\Z\otimes_{R(T)}\gGhat_{G},
$ which realizes the non-equivariant Schubert basis
$\{\,1\otimes\tgx{w}{n}(y| b) \mid w\in W_{\af}^0\,\}.$
 \end{remark}

Following the symmetric function realizations as in Theorem \ref{T:Peterson to symm funcs},
a natural expectation is that $gx_w^{(n)}(y):=gx_w^{(n)}(y|0)$ (resp. $\widetilde{gx}_w^{(n)}(y):=\widetilde{gx}_w^{(n)}(y|0)$)
give Schubert bases. However, in type $\mathrm{D}$, using \eqref{eq:non-eq:gd1} and \eqref{eq:non-eq:gd2} we saw
that the functions $\{gd_w(y)\mid w\in \Wafgr\}$ are not 
linearly independent.

Let $g\Gamma_{\Gsc}$ be the $\Z$-submodule of $\hL_{\Q}$
spanned by the specialized functions:
\[
g\Gamma_{\Gsc}
:=\sum_{w\in\Wafgr}\Z\,\tgx{w}{n}(y| 0).
\]
The specialization at $b=0$ induces a surjective
$\Z$-algebra homomorphism
\[
\Phi_{\Gsc}\colon
K_*(\Gr_{\Gsc})
\cong \Z\otimes_{\RT}\gGhat_{G}
\longrightarrow g\Gamma_{\Gsc},
\qquad
1\otimes\tgx{w}{n}(y| b)
\longmapsto \tgx{w}{n}(y| 0).
\]
Under this identification, $\Phi_{\Gsc}$ sends the structure-sheaf
class $\mathscr O_w^0$ to $\tgx{w}{n}(y| 0)$.
Thus we have an isomorphism 
$$g\Gamma_{\Gsc}\cong {K_*(\Gr_{\Gsc})}/{\Ker(\Phi_{\Gsc})}.$$
Next we prove that $\Ker(\Phi_{\Gsc})=\{0\}$ if $\Gsc=\Spin_{2n+1}(\C)$.

\begin{thm} The family of symmetric functions
$\tgq{w}{n}(y)\;(w\in \Wafgr)$ is linearly independent
over $\Z$. The $\Z$-span $g\Gamma_{\Gsc}$ of 
these functions is isomorphic to the
$K$-homology ring $K_*(\Gr_{\Spin_{2n+1}(\C)})$
such that the Schubert structure sheaf $\cO{w}^0$
corresponds to $\tgq{w}{n}(y).$
 \end{thm}
 \begin{proof}
Let $\tilde{G}_w^{B_n}(y)$ be the homology Schubert basis introduced by Pon, see \cite[Section 2.2]{Pon} and ${\widehat{B}}_w^{(n)}(y)$ be the
homology version of $\tgq{w}{n}(y)$ 
constructed in Appendix \ref{S:eq_hom}.
Then 
\begin{equation}
{\widehat{B}}_w^{(n)}(y)=\tilde{G}_w^{B_n}(y),
\end{equation}
see Theorem \ref{thm_hom_eq}. 
We consider the filtration on the ring of symmetric functions
given by the total degree in the variables $y$.
With respect to this filtration, the highest-degree
homogeneous component of $\widetilde{gq}_w^{(n)}(y)$
is the homology Schubert function $\widehat B_w^{(n)}(y)$.
Since $\{\tilde{G}_w^{B_n}(y)\mid w\in \Wafgr\}$ is linearly 
independent, we deduce the linear independent-ness of 
$\{\tgq{w}{n}(y)\}.$
 \end{proof}

If $\Gsc=\Spin_{2n}(\C)$ then $\Ker(\Phi_{\Gsc})$ is not a zero ideal. Following the equalities \eqref{eq:non-eq:gd1} and \eqref{eq:non-eq:gd2} define $\delta_{\Spin}:=\cO{\rho_{n-1}^{+(0)}}^0-\cO{\rho_{n-1}^{-(0)}}^0$ and denote its principal ideal by $I_{\Spin}:=(\delta_{\Spin})\subset K_*(\Gr_{\Spin_{2n}(\C)})$. Then it follows that $I_{\Spin}\subseteq \Ker(\Phi_{\Spin_{2n}(\C)})$.
Note that $I_{\Spin}$ also contains relations beyond the direct identifications of the two colors.  
For example in $\D{4}$, we have the following equality: 
$$\cO{213420}^0 + \cO{023420}^0 + \cO{13420}^0-\cO{123420}^0-\cO{321420}^0-\cO{23420}^0=\delta_{\Spin}\cO{420}^0.$$
We propose the following conjecture.
\begin{conj}\label{conj:D}
    We have the following $\Z$-algebra  isomorphism:
    $$g\Gamma_{\Spin_{2n}(\C)}\cong {K_*(\Gr_{\Spin_{2n}(\C)})}/{I_{\Spin}}.$$
\end{conj}

\subsection{Factorization of symmetric function bases}
Using the map $\Psi$ of Theorem \ref{T:Peterson to symm funcs} the factorization of the
$\ell$-basis of $\LG$ is transported to symmetric function bases.

\subsubsection*{Type $\B{}$} For $\mathrm{B}_n^{(1)}$, $u_1$ acts on $\sR$ by 
\begin{equation*}u_1(b_1)=\bar{b}_1\quad\text{and $u_1(b_k)=b_k$ for $2\le k\le n$}.
\end{equation*}This action preserves \(R(T)\).

\begin{thm}\label{thm_fac_type_B}
    Let \(v \in \widetilde{W}_{\B{n}^{(1)}}\). For $i\in I(\mathrm{B}_n)$, suppose that $v\kappa_i \in \WBexgr{n}$ and 
$\ell(v\kappa_i) = \ell(v) + \ell(\kappa_i)$. 
If $i$ is odd, then 
\[
\bar{v} := \pi_1^{-1} \, v \, \pi_1 \in \WBexgr{n},
\]
and otherwise \(v \in \WBexgr{n}\). Moreover, we have
\begin{equation}
\tgq{v\kappa_i}{n}(y| b) =
\begin{cases}
\tgq{\kappa_i}{n}(y| b)\tgq{v}{n} (y| b) & \text{if $i$ is even} \\
\tgq{\kappa_i}{n}(y| b) \tgq{\bar{v}}{n}(y| u_1(b)) & \text{if $i$ is odd}.
\end{cases}
\end{equation}
\end{thm}
\begin{proof}
The proof follows using Proposition \ref{P:factor kappa} and Theorem \ref{T:Peterson to symm funcs}.    
\end{proof}

\begin{remark} This factorization can be iterated: if $v$ or $\bar{v}$ has some $\kappa_j$ as length-additive right factor, the factorization can be applied again, regardless of whether the parameters are $b$ or $u_1(b)$: in the latter case one just applies $u_1$ to the usual factorization.
\end{remark}

\begin{example}
    In $\mathrm{B}_3^{(1)}$ we 
    have the following factorization:
\begin{align*}
   \widetilde{gq}^{(3)}_{s_1\kappa_1}(y| b)&=\widetilde{gq}^{(3)}_{\kappa_1}(y| b)\widetilde{gq}^{(3)}_{s_0}(y| u_1(b)),\\
   \widetilde{gq}^{(3)}_{s_0\kappa_2}(y| b)&=\widetilde{gq}^{(3)}_{\kappa_2}(y| b)\widetilde{gq}^{(3)}_{s_0}(y| b),\\ 
\widetilde{gq}^{(3)}_{s_2s_0\rho_3^{(1)}\rho_3^{(0)}\rho_3^{(1)}\rho_5^{(0)}}(y| b)&=\widetilde{gq}^{(3)}_{\kappa_1}(y| b)\widetilde{gq}^{(3)}_{s_2s_1\rho_3^{(0)}\rho_3^{(1)}\rho_3^{(0)}}(y| u_1(b))\\
   &=\widetilde{gq}^{(3)}_{\kappa_1}(y| b)\widetilde{gq}^{(3)}_{\kappa_3}(y| u_1(b))\widetilde{gq}^{(3)}_{s_2s_0}(y| b).
\end{align*}
\end{example}

\subsubsection*{Type $\D{}$} For type $\D{n}^{(1)}$, the actions of $u_1,u_{n-1}$ and $u_n$ on $\sR$ are given by:
\begin{equation*}u_1(b_1)=\bar{b}_1, ~u_1(b_n)=\bar{b}_n,~\text{and $u_1(b_k)=b_k$ for $2\le k\le n-1$}.
\end{equation*}
\[u_n(b_j)=
\begin{cases}
\bar b_{n+1-j},&1\le j\le n-1,\\
\bar b_1,&j=n,\ n\text{ even},\\
b_1,&j=n,\ n\text{ odd},
\end{cases}
\qquad
u_{n-1}(b_j)=
\begin{cases}
b_n,&j=1,\\
\bar b_{n+1-j},&2\le j\le n-1,\\
b_1,&j=n,\ n\text{ even},\\
\bar b_1,&j=n,\ n\text{ odd.}
\end{cases}
\]
This action preserves \(R(T)\).

 \begin{thm}\label{thm_fac_type_D}
    Let \(v \in \widetilde{W}_{\D{n}^{(1)}}\). For $i\in I(\mathrm{D}_n)$, assume \(v\kappa_i \in \WDexgr{n}\) and 
\(\ell(v\kappa_i) = \ell(v) + \ell(\kappa_i)\).
\begin{enumerate}
    \item If $i$ is even and $1\leq i\leq n-2$ then  \(v \in \WDexgr{n}\).\\
    \item If $i$ is odd and $1\leq i\leq n-2$ then 
\[
\bar{v} := \pi_1^{-1} \, v \, \pi_1 \in \WDexgr{n}.
\]
\item If $i=n-1$, then 
\[
\bar{v}_{(n-1)} := \pi_{n-1}^{-1} \, v \, \pi_{n-1} \in \WDexgr{n}.
\]
\item If $i=n$, then 
\[
\bar{v}_{(n)} := \pi_{n}^{-1} \, v \, \pi_{n} \in \WDexgr{n}.
\]
\end{enumerate}

Moreover, we have
\begin{equation}
\widetilde{gd}^{(n)}_{v\kappa_i}(y| b) =
\begin{cases}
\widetilde{gd}^{(n)}_{\kappa_i}(y| b)\,
\widetilde{gd}^{(n)}_{v}(y| b), & 1 \le i \le n-2, i \text{ is even} \\
\widetilde{gd}^{(n)}_{\kappa_i}(y| b)
\widetilde{gd}^{(n)}_{\bar{v}}(y| u_1(b)), & 1 \le i \le n-2, i \text{ is odd}\\
\widetilde{gd}^{(n)}_{\kappa_i}(y| b)
\widetilde{gd}^{(n)}_{\bar{v}_{(n-1)}}(y| u_{n-1}(b)), & i= n-1\\
\widetilde{gd}^{(n)}_{\kappa_i}(y| b)
\widetilde{gd}^{(n)}_{\bar{v}_{(n)}}(y| u_{n}(b)), & i= n.
\end{cases}
\end{equation}
\end{thm}
\begin{proof}
The proof follows using Proposition \ref{P:factor kappa} and Theorem \ref{T:Peterson to symm funcs}.    
\end{proof}

\section{
{Type \(D_\infty\) symmetric functions and parameter changes}
}
\label{subsec_D_infty}
In this section, we recall the definition of the function $gd_\la(y| b)$ of Nakagawa and Naruse, characterized as the dual of the factorial $GD$-function of Ikeda and Naruse. We also give an explicit formula for the expansion of $gd_\la(y| b)$ in terms of $gd_\mu(y| c)$ for two different sets of equivariant parameters $b=(b_1,b_2,\dots)$, $c=(c_1,c_2,\dots)$ and describe the transition coefficients in terms of double Grothendieck polynomials.

A strict partition $\lambda\in \SP$ is a strictly decreasing sequence of positive integers
$\lambda=(\lambda_1>\lambda_2>\dotsm>\lambda_\ell> 0)$. Its length $\ell$ is denoted $\ell(\lambda)$ and size of $\lambda$ is given by $|\lambda|=\lambda_1+\dots+\lambda_{\ell}$. We introduce a coefficient ring \[
R_\infty
:=\Z[b_i,(1-b_i)^{-1}\mid i\ge1]
=\Z[b_i,\bar b_i\mid i\ge1],
\qquad
\bar b_i=-\frac{b_i}{1-b_i}.
\]

Let $\lambda\in \SP$ and $ r=\ell(\lambda)\leq n.$ We recall the definition of the factorial $GP$-function following \cite{IN}: 
\begin{align*}
GP_\lambda(x_1,\ldots,x_n| b)&=
    \frac{1}{(n-r)!}
    \sum_{w\in \mathfrak{S}_n}
    w\left(
\prod_{i=1}^r[x_i| b]^{\lambda_i}
    \prod_{i=1}^r \prod_{j=i+1}^n 
    \frac{x_i\oplus x_j}{x_i\ominus x_j}
    \right),
\end{align*}    
where $[x| b]^k=(x\ominus b_1)\cdots (x\ominus b_{k})$ for $k\geq 1$ and $w$ permutes $x_1,\ldots,x_n$. The polynomials $GP_\lambda(x_1,x_2,\dots | b)$ are symmetric in $(x_1,\ldots,x_n)$ with coefficients in $R_\infty$. Let $\SP(r)$ denote the set of all strict partitions of length $r$. Let $\lambda\in \SP(r)$ and $n$ be an even integer such that $r\leq n$. We define 
\[GD_\lambda(x_1,x_2,\dots,x_{n}| b_1,b_2,\dots)=
    GP_\lambda(x_1,x_2,\dots, x_{n}| b_1,b_2,\dots) 
\]

Also, we note that $$GD_\lambda(x_1,\dots,x_{2m-2}| b)=GD_\lambda(x_1,\dots,x_{2m-2},0,0| b).$$ 
We denote the inverse limit of $GD_\lambda(x_1,\dots,x_{2m}| b)$ by $GD_\lambda(x_1,x_2,\dots| b)$.

\begin{remark}
The $GP_\lambda(x_1,x_2,\dots| b_1,b_2,\dots )$ defined above is the same as $GP_\lambda(x_1,x_2,\dots| \ominus b_1,\ominus b_2, \dots)$ in the notation of \cite{IN}. Moreover, the family of functions $GD_\lambda(x_1,x_2,\dots| b)$ defined above represents the Schubert classes of the maximal isotropic Grassmannian for type $\D{}$. More precisely, the Schubert class corresponding to $\lambda$ in the torus-equivariant $K$-theory of the maximal isotropic Grassmannian of type $\mathrm{D}_{n+1}$ corresponds to $GD_\lambda(x_1,x_2,\dots,x_n| b)$ or $GD_\lambda(x_1,x_2,\dots,x_{n},0| b)$ depending on whether $n$ is even or odd respectively. 
\end{remark}

The `dual factorial $gd$-function' ${gd}_\lambda(y| b)$ is defined to be the dual of $GD_\lambda(x| b)$ following \cite{NN:Uni}:
\begin{equation}\label{E:gdGD Cauchy}
\prod_{i,j} \frac{1 - \bar{x}_i y_j}{1 - x_i y_j}
= \sum_{\lambda\in \SP} gd_\lambda(y| b) GD_\lambda(x| b).
\end{equation}
In other words, there is a bilinear map between the two $R_\infty$-algebras generated by $\{GD_\lambda(x| b)\}$ and $\{gd_\lambda(y| b)\}$ respectively satisfying
\begin{equation}\label{E:gdGD pairing}
\langle {gd}_\lambda(y| b), GD_\mu(x| b) \rangle = \delta_{\lambda\mu}.
\end{equation}

Given $\lambda \in \SP$, let $D^*(\lambda)$ denote its shifted diagram, where
$$D^*(\lambda) := \{(i, j) | 1 \le i \le j < i+\lambda_i\}.$$
$\SP$ is a graded poset under containment of shifted diagrams. For each box $(i,j)$ in the $i$-th row and $j$-th column in $D^*(\lambda)$, we define its content $c(i,j)$ by
\[c(i,j)=\begin{cases}
    j-i+1 &\text{ if }i<j\\
    0 & \text{ if }i=j \text{ is odd}\\
    1 & \text{ if }i=j \text{ is even}
\end{cases}\\
\]
The contents of the boxes in the shifted diagram are shown in the following diagram. 
\[\begin{ytableau}
0&2&3&4&5&6&\dots\\
\none &1&2&3&4&5&\dots\\
\none&\none&0&2&3&4&\dots\\
\none&\none&\none&1&2&3&\dots\\
\none&\none&\none&\none&\vdots&\vdots&\vdots\\
\end{ytableau}\]

The infinite Grassmannian elements $\Winfgr$ are given by the minimal length coset representatives $$\Winfgr=\{w\in W_{\mathrm{D}_\infty}\mid ws_i>w \text{ for }i\ge1\}$$

For type $\mathrm{D}_\infty$, we define $\rho_i^{(1)}= s_i \dotsm s_2s_1$ for $i\geq 1$ and  $\rho_i^{(0)} = \auto(\rho_i^{(1)})$.
The infinite Grassmannian elements of 
$\Winfgr$ can be parametrized in terms of strict partitions using a bijection $w:\SP\to \Winfgr$ as in \eqref{eq_w_lam}.  Thus, the size $|\lambda|$ coincides with the Coxeter length $\ell(w_{\lambda})$.

Let $\hL_{R_\infty}$ denote the completion, with respect
to the $y$-degree filtration, of the ring of symmetric functions
with coefficients in $R_\infty$.
We define an action of $\Winf$ on the ring $\widehat{\Lambda}_{R_\infty}({y})$ by setting
\begin{align}
s_0 f(y| b) &:= \Omega(b_1| y)\Omega(b_2| y) \cdot f(y | \bar{b}_2, \bar{b}_1, b_3, \ldots),  \\
s_i f(y| b) &:= f(y | \ldots, b_{i+1}, b_i, \ldots) \quad \text{for } i > 0,
\end{align}
for $f\in\hL_{R_\infty}$.
The action of $s_0$ is well-defined since
$\Omega(b_1|y)\Omega(b_2|y)$ belongs to $\hL_{R_\infty}$.
The operator $T_i$ for $i\geq 0$ is defined by
\begin{equation}\label{eq_T_inf}
\hat{T}_i={c(\alpha_i)}^{-1}(s_i-1).
\end{equation}

For type $\mathrm{D}_\infty$, we have
\begin{align*}
    c(\alpha_0)= b_1 \oplus b_2, \quad
    c(\alpha_i) = b_{i+1} \ominus b_i \text{ for $ i\ge 1$}
\end{align*}

These operators satisfy the braid relations as in $W_{\mathrm{D}_\infty}$.
Thus the operators $\{\hat{T}_w~|~w\in \Winf\}$ are well defined. 

\begin{prop}[\cite{NN:Uni}, see also \cite{ISY2}]\label{prop_T_i}
For $i\in \Z_{\ge 0}$ we have
\[
\hat{T}_i gd_\lambda(y| b)=
\begin{cases}
gd_{s_i \lambda}(y| b)
& \text{if $s_i w_\lambda>w_\lambda$ and $s_i w_\lambda\in \Winfgr$,}\\
0
& \text{if $s_i w_\lambda>w_\lambda$ and $s_i w_\lambda\notin \Winfgr$,}\\
- gd_{\lambda}(y| b)
& \text{if $s_i w_\lambda<w_\lambda$.}
\end{cases}
\]
Moreover, for $\lambda\in\SP$, we have
 $$gd_{\lambda}(y| b)=\hat{T}_{w_\lambda}(1).$$
\end{prop}

For each \(w\in \mathfrak{S}_\infty=\bigcup_{n\geq 1}\mathfrak{S}_n\), Lascoux and Sch\"utzenberger \cite{LaSc} defined double Grothendieck polynomials \(\G_w(b;c)\). For every \(i\geq 1\), define an operator $\hat{D}_i$ on $\mathbb Z[b_i,c_i,\bar{b}_i,\bar{c}_i:i\ge1]$ by \(\hat{D}_i:=\hat{T}_i+1\), where \(\hat{T}_i\) is defined in \eqref{eq_T_inf}. Note that $\hat{D}_i(c_j)=c_j$ for all $i,j\geq 1$. Let $w_0^{(n)}$ denote the longest element of $\mathfrak{S}_n$. The family $ \{\G_w(b;c)\mid w\in \mathfrak{S}_\infty\}$
is uniquely determined by
\[
\G_{w_0^{(n)}}(b;c)
=\prod_{i+j\le n}
b_i\ominus c_j,
\]
\begin{equation}\label{eq_p_i_equi}
\hat{D}_i(\G_w(b;c))=
\begin{cases}
\G_{ws_i}(b;c), & \text{if } \ell(ws_i)<\ell(w),\\
\G_w(b;c), & \text{if } \ell(ws_i)>\ell(w).
\end{cases}
\end{equation}

\begin{remark}
 A permutation \(w\in \mathfrak{S}_\infty\) is called \(k\)-\emph{Grassmannian} if
\[
w(1)<\cdots<w(k)
\qquad\text{and}\qquad
w(k+1)<w(k+2)<\cdots.
\]
Such a permutation is uniquely determined by \(w(1),\dots,w(k)\). Let $\mathscr{P}_k$ denote the partitions of length less than or equal to $k$. For \(\lambda=(\lambda_1,\dots,\lambda_k)\in\mathscr{P}_k\), there is a unique
\(k\)-Grassmannian permutation \(w_\lambda^{(k)}\in \mathfrak{S}_\infty\) satisfying $w_\lambda^{(k)}(i)=\lambda_{k-i+1}+i$ for $1\leq i\leq k$. Ikeda and Naruse give another description of double Grothendieck polynomials for $k$-Grassmann permutations in \cite[Section 2.4]{IN}.
\end{remark}

Let $\lambda \in \mathscr{P}_{k+1}$ and $\delta_k=(k,k-1,\dots,1)$, the factorization property as in \cite[Proposition 2.3]{IN} implies
\begin{align*}
GP_{\delta_{k}}(b_1,\ldots,b_{k+1}| c)
&=GP_{\delta_{k}}(b_1,\ldots,b_{k+1}),\\
GP_{\lambda+\delta_{k}}(b_1,\ldots,b_{k+1}| c)
&=GP_{\delta_{k}}(b_1,\ldots,b_{k+1})\G_{w_\lambda^{(k+1)}}(b;c).
\end{align*}

 Moreover, if $\mu:=\lambda+\delta_k$ then
$$\G_{w_\lambda^{(k+1)}}(b;c)=\G_{w_\mu w_{\delta_k}^{-1}}(b;c).$$

Note that if $\lambda\in \mathscr{P}_k\subset \mathscr{P}_{k+1}$ then $\lambda+\delta_k\in \SP(k)$. For $\lambda \in \mathscr{P}_{k+1}\setminus\mathscr{P}_{k}$, we have $\lambda+\delta_k\in \SP(k+1)$. Then we have the following expansion formula:

 \begin{lem}\label{lem_odd_k}
 If $k$ is odd, then 
 \begin{equation}\label{eq_odd_k}
gd_{\delta_k}(y| b)=\sum_{\mu\in \SP(k)\cup\SP(k+1)} \big(\G _{w_{{\mu}}w_{\delta_k}^{-1}}(b;c)\big)gd_{\mu}(y| c). 
 \end{equation}
 \end{lem}
\begin{proof}



If $\lambda\in \mathscr{P}_{k+1}$ then using the similar argument as in \cite[Proposition 3.6]{IN} we have 
\begin{equation}\label{eq_exp_GP}
    GP_{\lambda+\delta_{k}}(x_1,x_2,\dots,x_{k+1}| c)=\sum_{\substack{\eta\in \mathscr{P}_{k+1}\\ \eta\subset \lambda}}C_{\lambda\eta}(b,c)GP_{\eta+\delta_{k}}(x_1,x_2,\dots,x_{k+1}| b).
\end{equation}
Thus $C_{\lambda\phi}(b,c)=\langle gd_{\delta_k}(y| b),GD_{\lambda+\delta_k}(x| c)\rangle$ and we have the expansion
$$gd_{\delta_k}(y| b)=\sum_{\lambda\in \mathscr{P}_{k+1}}C_{\lambda\phi}(b,c)gd_{\delta_k+\lambda}(y| c).$$

Since $k$ is odd, substituting $y=b_{\delta_{k}}=( b_{k+1},\dots,b_2, b_1,0,\dots)$ in \eqref{eq_exp_GP} and using the vanishing theorem as in \cite[Proposition 7.1]{IN} we have
$$GP_{\lambda+\delta_{k}}(b_{\delta_{k}}| c)=C_{\lambda\phi}(b,c)GP_{\delta_{k}}(b_{\delta_{k}}| b).$$

Then $$C_{\lambda\phi}(b,c)=\frac{GP_{\delta_{k}+\lambda}(b_{\delta_{k}}| c)}{GP_{\delta_{k}}(b_{\delta_{k}}| b)}=\frac{GP_{\delta_{k}}(b_{\delta_{k}})\G_{w_\lambda^{(k+1)}}(b;c)}{GP_{\delta_{k}}(b_{\delta_k})}=\G_{w_\lambda^{(k+1)}}(b;c).$$
\end{proof}

\begin{example}\label{ex_exp}
\begin{align*}
gd_1(y| b)=&gd_1(y| c)+\G_{s_2}(b;c)gd_2(y| c)+\G_{s_3s_2}(b;c)gd_3(y| c)+\G_{s_4s_3s_2}(b;c)gd_4(y| c)\\
&+\dots+\G_{s_1s_2}(b;c)gd_{(2,1)}(y| c)+\G_{s_3s_1s_2}(b;c)gd_{(3,1)}(y| c)+\dots\\
gd_2(y| b)=&(1-\G_{s_2}(b;c))gd_2(y| c)+(\G_{s_3}(b;c)-\G_{s_3s_2}(b;c))gd_3(y| c)+(\G_{s_4s_3}(b;c)\\ &-\G_{s_4s_3s_2}(b;c))gd_4(y| c)+\dots+(\G_{s_1}(b;c)-\G_{s_1s_2}(b;c))gd_{(2,1)}(y| c)\\
&+(\G_{s_3s_1}(b;c)-\G_{s_3s_1s_2}(b;c))gd_{(3,1)}(y| c)+\dots\\
gd_{(2,1)}(y| b)=&(1-\G_{s_1}(b;c))gd_{(2,1)}(y| c)+(\G_{s_3}(b;c)-\G_{s_3s_1}(b;c))gd_{(3,1)}(y| c)\\
&+(\G_{s_4s_3}(b;c)-\G_{s_4s_3s_1}(b;c))gd_{(4,1)}(y| c)+\dots
\end{align*}
\end{example}

\begin{lem}\label{lem_even_k}
    If $k$ is even, then 
$$gd_{\delta_k}(y| b)=\sum_{\mu\in \SP(k)} \big(\G _{w_{{\mu}}w_{\delta_k}^{-1}}(b;c)- \G _{w_{{\mu}}w_{\delta_k}^{-1}s_1}(b;c)\big)gd_{\mu}(y| c),$$ 
\end{lem}
 \begin{proof}
The expansion of $gd_{\delta_{k-1}}(y| b)$ follows from Lemma \ref{lem_odd_k}. We apply finitely many $\hat{T}_i$ in both sides of \eqref{eq_odd_k}. We use $\hat{T}_i=\hat{D}_i-1$ for $i\geq 1$ and Proposition \ref{prop_T_i} and the relation in \eqref{eq_p_i_equi}. Note that $$w_{\delta_k}=s_1\dots s_{k-1}s_kw_{\delta_{k-1}}, ~gd_{\delta_k}(y| b) = \hat{T}_1\cdots \hat{T}_{k} gd_{\delta_{k-1}}(y| b).$$  Now the proof follows using the induction.  
 \end{proof}

 If $\lambda\in \SP(k)$ then  $\lambda\setminus \delta_k$ lies in $\mathscr{P}_k$.
Let $R(\lambda\setminus\delta_k)$ denote the set of the contents of the removable boxes of $\lambda\setminus \delta_k$. A box is called a removable box of $\lambda\setminus \delta_k$ if it can be removed from the Young diagram of $\lambda\setminus \delta_k$ so that the remaining diagram is still a Young diagram of some partition. Note that the simple reflections corresponding  to two elements in $R(\lambda\setminus\delta_k)$ pairwise commute.
So for any subset $J\subset R(\lambda\setminus\delta_k)$ we denote $w_J=\prod_{j\in J}s_j$ which is independent of the order of the elements of $J$. 
If $\lambda\in \SP(k)$ then $\lambda\setminus \delta_{k-1}$ also lies in $\mathscr{P}_k$ and we can define $R(\lambda\setminus\delta_{k-1})$ similarly as above. The proof of the following Theorem follows similarly as in \cite[Theorem 4.10]{BIIY}.

\begin{thm}\label{thm:gp expansion}
Let $\lambda=(\lambda_1,\dots,\lambda_k)\in \SP(k)$ and $\delta_k=(k,k-1,\dots,1)$.  Then
 \begin{align*}
gd_\lambda(y| b)&=\sum_{\substack{\mu:\,\mu\supset \lambda\\ \ell(\mu)\leq k_2}}\Big(\sum_{J\subset R(\lambda\setminus \delta_{k_1})}
(-1)^{|J|} \G_{w_{\mu}w_\lambda^{-1}w_J}(b;c)\Big)
 gd_{\mu}(y| c),
 \end{align*}
 where $k_1$ is the largest odd integer and $k_2$ is the smallest even integer such that $k_1\leq k\le k_2$.
\end{thm}

\begin{proof}
    The proof follows by induction on $|\lambda|-|\delta_k|$. We start with the expansion of $gd_{\delta_k}(y| b)$ following Lemma \ref{lem_odd_k} or Lemma \ref{lem_even_k} according as $k$ is odd or even. Using Proposition \ref{prop_T_i}, $gd_\lambda(y| b)$ can be obtained from $gd_{\delta_k}(y| b)$ by applying finitely many $\hat{T}_i$ in the left hand side. On the right hand side we use $\hat{T}_i=\hat{D}_i-1$ for $i\geq 1$ and \eqref{eq_p_i_equi}. 
\end{proof}

\begin{remark}
Let $\mathfrak{S}_\infty$ be the subgroup of $W_{D_\infty}$ generated by $s_i$ for $i>0$. Suppose \(\lambda\subset\mu\) are strict partitions and \(\ell(\mu)\le k_2\). Then one can directly check the following equality:
\[
\{u\in \mathfrak{S}_\infty:u*w_\lambda=w_\mu\}
=
\{w_\mu w_\lambda^{-1}w_J:J\subseteq R(\lambda\setminus \delta_{k_1})\},
\]
where $*$ is the Demazure product. These subsets \(J\) give distinct elements, and
\[
\ell(w_\mu w_\lambda^{-1}w_J)
=
|\mu|-|\lambda|+|J|.
\] 
Therefore the expansion in Theorem \ref{thm:gp expansion} can also be expressed by the following as a generalization of \cite[Proposition 5.36]{LLS}.
\begin{align*}
    gd_\la(y| b) &= \sum_{\substack{\mu:~ \mu\supset\lambda\\\ell(\mu)\leq k_2}} \Big( \sum_{\substack{u\in {\mathfrak{S}_\infty} \\ u*w_\la=w_\mu}} (-1)^{|\mu/\la|-\ell(u)} \G_u(b;c) \Big)gd_\mu(y| c).
\end{align*}
\end{remark}




\section{Small classes}\label{S: small_cl}

In this section, we express certain affine double $gq$- and $gd$-functions in terms of the type $\mathrm{D}_\infty$ functions $gd_\lambda(y| b)$, possibly after a suitable specialization of the equivariant parameters.

Let $\lambda(w)=(\lambda_1,\dots,\lambda_\ell)$ be nonempty for $w\in \Wafgr$. Whenever \((\lambda_1-1,\lambda_2,\ldots,\lambda_\ell)\) is a strict partition, denote it by \(\lambda'(w)\), with the convention \(\lambda'(w)=\varnothing\) if \(\lambda(w)=(1)\). Let \(w'\) be the corresponding affine Grassmannian element. In type \(\D{}\), retain the existing sign if a part of size \(n-1\) remains, assign the sign \(-\) if such a part is newly created, and omit the sign if no such part remains.



\begin{thm}\label{small class type BD}
\begin{enumerate}
\item Let $w\in \WBgr{n}$ with $\la(w)$ strict and $\la_1 \le n-1$. 
Then 
\begin{align}\label{E:gqgd}
    gq_w^{(n)}(y| b)= gd_{\lambda(w)}(y| b).
\end{align}
\item Let $w\in \WDgr{n}$ with $\la(w)$ strict and $\la_1\le n-1$.
If either $\la_1<n-1$ or $\la_1=n-1$ and $\la(w)$ is marked $-$, then
\begin{align}\label{E:gd+}
 gd_w^{(n)}(y| b)= gd_{\lambda(w)}(y| b).
\end{align}
If $\la(w)_1=n-1$ and $\la(w)$ is marked $+$ then
\begin{align}\label{E:gd-}
 gd_w^{(n)}(y| b)= gd_{\lambda(w)}(y| b_1,\dotsc,b_{n-1},\bar{b}_n).
\end{align}
\end{enumerate}
\end{thm}
\begin{proof}
Equations \eqref{E:gqgd} and \eqref{E:gd+}
follow from the fact that the subdiagrams of $\D{\infty}$, $\mathrm{B}_n^{(1)}$, and $\D{n}^{(1)}$, when restricted to the nodes $\{0,1,\dotsc,n-1\}$, all agree.
Equation \eqref{E:gd-} follows from \eqref{E:gd+} by Theorem \ref{T:gd auto}.
\end{proof}

\begin{example} For $n=5$ the relevant partitions in Theorem \ref{small class type BD}
are the strict partitions contained in $(4,3,2,1)$.
In type $D_5^{(1)}$, if $\la$ contains the part $n-1=4$ and is marked $-$ (resp. $+$)
then $s_{n-1}=s_4$ (resp. $s_n=s_5$) is used in $w_\la\in \WDgr{5}$.
\[
\begin{ytableau}
0 & 2 & 3 & 4 \\
\none&1 & 2 & 3\\
\none&\none&0&2 \\
\none&\none&\none&1
\end{ytableau}\qquad \begin{ytableau}
0 & 2 & 3 & 5 \\
\none&1 & 2 & 3\\
\none&\none&0&2 \\
\none&\none&\none&1
\end{ytableau}
\]
\end{example}

Next, we extend Theorem \ref{small class type BD} for more affine Grassmannian elements for both types $\B{}$ and $\D{}$ after some appropriate reparameterization of equivariant parameters.

\subsection{Type B}
Let $\pi:\Z[b_i,\bar{b}_i: 1\le i\le 2n]\to \Z[b_1,\bar{b}_1,\dots, b_n,\bar{b}_n]$ be the map defined by 
\[\pi(b_i)=
\begin{cases}
   b_i& \text{for } 1\leq i\leq n\\
   \bar{b}_{2n+1-i} & \text{for } n+1\leq i\leq 2n-1\\
   {b}_1 & \text{for }  i=2n\\
\end{cases}.
\]

Let us write \(\alpha_i^\infty,s_i^\infty,T_i^\infty\) for the type \(\D{\infty}\) simple roots, operators for $i\geq 1$ and \(\alpha_j^{\B{}},s_j^{\B{}},T_j^{\B{}}\) for the type \(\B{n}^{(1)}\) simple roots, operators for $1\le j\le n$.

\begin{lem}\label{lem_1}
    Let $\lambda=(\lambda_1,\lambda_2,\dots,\lambda_\ell)$ be a strict partition. The equivariant parameters in the function \(gd_\lambda(y| b)\) depend only on
$b_1,\ldots,b_{\lambda_1+1}.$ 
\end{lem}
\begin{proof}
Using Proposition \ref{prop_T_i}, 
\[
gd_\lambda(y| b)=T^\infty_{w_\lambda}(1).
\]
The canonical reduced expression for \(w_\lambda\) uses only the generators $s_0,s_1,\ldots,s_{\lambda_1}.$
The operator \(T_0^\infty\) involves only \(b_1,b_2\), and \(T_j^\infty\), for \(j>0\), involves only \(b_j,b_{j+1}\). Thus no parameter beyond \(b_{\lambda_1+1}\) occurs in the function \(gd_\lambda(y| b)\). 
\end{proof}

\begin{lem}\label{lem_2}
    Suppose $
2\le r\le \lambda_1-1$ and 
$\lambda_2\le r-2$.
Then
\[
T_r^\infty gd_\lambda(y| b)=0,
\qquad
s_r^\infty gd_\lambda(y| b)=gd_\lambda(y| b).
\]
\end{lem}
\begin{proof}
    Write the canonical reduced factorization as $
w_\lambda=v\,\rho_{\lambda_1}^{(0)}.$
Every generator in \(v\) has index at most \(\lambda_2\). Therefore \(s_r\) commutes with \(v\), because \(\lambda_2\le r-2\). The braid relations give
\[
s_r\rho_m^{(0)}
=
\rho_m^{(0)}s_{r+1} \text{ and consequently }
s_rw_\lambda=w_\lambda s_{r+1}.
\]
Now the proof follows using Proposition \ref{prop_T_i}.
\end{proof}

\begin{lem}\label{intertwine_n}
    If \(f\) depends only on $b_1,\ldots,b_{n+1}$ then
\[
T_n^B\pi(f)
=
(2-\bar b_n)\,\pi(T_n^\infty f).
\] 
\end{lem}
\begin{proof}
The specialization $\pi$ sends the pair \((b_n,b_{n+1})\) to $(b_n,\bar b_n)$.
Thus, if \(f\) depends only on $b_1,\ldots,b_{n+1}$ then
\[
s_n^B\pi(f)=\pi(s_n^\infty f).
\]
Furthermore,
\[
\pi\bigl(c(\alpha_n^{\infty})\bigr)
=
\bar b_n\ominus b_n
=
(2-\bar b_n)\bar b_n
=
(2-\bar b_n)c(\alpha_n^{B}).
\]
Dividing the intertwining identity for the numerators completes the proof.
\end{proof}

\begin{lem}\label{intertwine_i}
    Let $n+1\le i\le 2n-2$ and
$r=2n-i$ and \(f\) depends only on \(b_1,\ldots,b_{i+1}\) and satisfies
$s_r^\infty f=f$.
Then
\[
T_r^B\pi(f)=\pi(T_i^\infty f).
\]
\end{lem}
\begin{proof}
    Observe that \(s_r^{\B{}}\) interchanges \(b_r,b_{r+1}\), so it also interchanges their formal inverses. Under \(\pi\), those inverses occur at positions \(i+1,i\). Therefore
$s_r^B\pi(f)=\pi(s_r^\infty s_i^\infty f).$ 
Note that $s_r^\infty s_i^\infty=s_i^\infty s_r^\infty$ and the assumed symmetry gives
$s_r^B\pi(f)
=
\pi(s_i^\infty f).$ Furthermore, $$\pi(c(\alpha_{i}^{\infty}))=c(\alpha_{2n-i}^{\mathrm{B}}) \text{ for } n+1\leq i\leq 2n-2.$$
This completes the proof.
\end{proof}

\begin{thm}\label{thm_sm_cl_b}
Let $w\in \WBgr{n}$ and write $\lambda=\lambda(w)$.
\begin{enumerate}
    \item If $\lambda_1=n$ and $\lambda'(w)$ is strict, then
    \begin{align}\label{sm_cl_n}
        gq_w^{(n)}(y| b)
        =
        (2-\bar{b}_n)\,gd_{\lambda(w)}(y| \pi(b)).
    \end{align}

    \item If $\lambda_1\geq n+1$ and $\ell(w)\leq 2n-2$, then
    \begin{align}\label{sm_cl_geq_np1}
        gq_w^{(n)}(y| b)
        =
        (2-\bar{b}_{2n-\lambda_1})\,
        gd_{\lambda(w)}(y| \pi(b))
        -
        \frac{1}{1-b_{2n+1-\lambda_1}}\,
        gd_{\lambda'(w)}(y| \pi(b)).
    \end{align}
\end{enumerate}
\end{thm}

\begin{proof}
Suppose $\lambda_1=n$ and $\lambda'(w)$ is strict. By Theorem~\ref{small class type BD}, we have
\[
    gq_{w'}^{(n)}(y| b)
    =
    gd_{\lambda'(w)}(y| b).
\]
By Lemma \ref{lem_1}, the right-hand side depends only on \(b_1,\ldots,b_n\), so it is unchanged by \(\pi\). Applying $T_n^{\mathrm{B}}$ to both sides and using Lemma \ref{intertwine_n}, we get
\[
    gq_w^{(n)}(y| b)
    =
    (2-\bar{b}_n)\,
    gd_{\lambda(w)}(y| \pi(b)),
\]
which proves \eqref{sm_cl_n}.

Next we prove \eqref{sm_cl_geq_np1}. Let $\lambda_1=a$, where $a\geq n+1$. Since $\ell(w)\leq 2n-2$, we have
\[\lambda_j\leq 2n-a-2
    \qquad \text{for all } j\geq 2.
\]
Since, $a\geq n+1$, this also gives $\lambda_j\leq n-3, \text{ for all } j\geq 2.$  Set
\[\lambda^{(k)}=(n+k,\lambda_2,\ldots,\lambda_r),
    \qquad 0\leq k\leq a-n.
\]
Thus $\lambda^{(0)}=(n,\lambda_2,\ldots,\lambda_\ell)$ and
$\lambda^{(a-n)}=\lambda$. The case for $k=0$ follows by using \eqref{sm_cl_n} 
\begin{equation}\label{eq_sm_gq}
     gq_{w_{\lambda^{(0)}}}^{(n)}(y| b)
    =
    (2-\bar{b}_n)\,
    gd_{\lambda^{(0)}}(y| \pi(b)).
\end{equation}

Now fix \(k\) with \(1\le k\le a-n\). Let $i=n+k$ and $
r=2n-i.$
Thus, we have
\begin{equation}\label{eq_lam_j}
  \lambda_j\le r-2, \text{ for all }j\ge 2\text{ and } 2\le r\le i-2. 
\end{equation}

Lemma \ref{lem_2}, applied to the partition \((i-1,\lambda_2,\dots,\lambda_\ell)\), therefore gives
\[
s_r^\infty gd_{(i-1,\lambda_2,\dots,\lambda_\ell)}(y| b)
=gd_{(i-1,\lambda_2,\dots,\lambda_\ell)}(y| b).
\]
By Lemma \ref{lem_1}, \(gd_{(i-1,\lambda_2,\dots,\lambda_\ell)}(y| b)\) depends only on \(b_1,\ldots,b_{i}\). Thus Lemma \ref{intertwine_i} applies, giving
\[
T_r^B \pi(gd_{(i-1,\lambda_2,\dots,\lambda_\ell)}(y| b))
=
\pi(T_i^\infty gd_{(i-1,\lambda_2,\dots,\lambda_\ell)}(y| b))
=
\pi(gd_{(i,\lambda_2,\dots,\lambda_\ell)}(y| b)).
\]
The corresponding type \(\B{}\) canonical words give
\[T_r^B gq^{(n)}_{w_{(i-1,\lambda_2,\dots,\lambda_\ell)}}(y| b)=gq^{(n)}_{w_{(i,\lambda_2,\dots,\lambda_\ell)}}(y| b).
\]
Indeed, \(s_r\) commutes past the tail factors by \eqref{eq_lam_j}, and it extends the first-row Pieri factor by one simple reflection.
 If \(i\ge n+2\), Lemma \ref{lem_2} applies to \((i-2,\lambda_2,\dots,\lambda_r)\), so
\[
s_r^\infty gd_{(i-2,\lambda_2,\dots,\lambda_r)}(y| b)=gd_{(i-2,\lambda_2,\dots,\lambda_r)}(y| b).
\]

Thus, applying $T_{n-1}^{\mathrm{B}}$ to both sides of \eqref{eq_sm_gq} gives
\[
    gq_{w_{\lambda^{(1)}}}^{(n)}(y| b)
    =
    (2-\bar{b}_{n-1})\,
    gd_{\lambda^{(1)}}(y| \pi(b))
    -
    \frac{1}{1-b_n}\,
    gd_{\lambda^{(0)}}(y| \pi(b)).
\]
This is precisely \eqref{sm_cl_geq_np1} for $\lambda_1=n+1$.

We proceed by induction on $k$. Suppose that the formula holds for
$\lambda^{(k-1)}$, where $1\leq k-1<a-n$. Applying
$T_{n-k}^{\mathrm{B}}$ and using Lemma \ref{intertwine_i}, we obtain
\[
    T_{n-k}^{\mathrm{B}}
    \left(
        (2-\bar{b}_{n-k+1})\,
        gd_{\lambda^{(k-1)}}(y| \pi(b))
    \right)
    =
    (2-\bar{b}_{n-k})\,
    gd_{\lambda^{(k)}}(y| \pi(b))
    -
    \frac{1}{1-b_{n-k+1}}\,
    gd_{\lambda^{(k-1)}}(y| \pi(b)).
\]
Moreover, the remaining term in the induction hypothesis vanishes as
\[
    T_{n-k}^{\mathrm{B}}
    \left(
        \frac{1}{1-b_{n-k+2}}\,
        gd_{\lambda^{(k-2)}}(y| \pi(b))
    \right)
    =0.
\]
Therefore
\[
    gq_{w_{\lambda^{(k)}}}^{(n)}(y| b)
    =
    (2-\bar{b}_{n-k})\,
    gd_{\lambda^{(k)}}(y| \pi(b))
    -
    \frac{1}{1-b_{n-k+1}}\,
    gd_{\lambda^{(k-1)}}(y| \pi(b)),
\]
using Lemma \ref{intertwine_i}.
Taking $k=a-n$ gives
\[
    gq_w^{(n)}(y| b)
    =
    (2-\bar{b}_{2n-a})\,
    gd_{\lambda(w)}(y| \pi(b))
    -
    \frac{1}{1-b_{2n+1-a}}\,
    gd_{\lambda'(w)}(y| \pi(b)),
\]
where $a=\lambda_1$. This proves \eqref{sm_cl_geq_np1}.
\end{proof}

Let us denote the functions $gq_{\rho_i^{(0)}}^{(n)}(y| b)$ and $gd_{\rho_i^{(0)}}(y| b)$ by $gq_{i}^{(n)}(y| b)$ and $gd_{i}(y| b)$ respectively.

\begin{cor}\label{cor_per_cl_k}
\[ gq_{i}^{(n)}(y| b)=
\begin{cases}
   (2-\bar{b}_n) gd_n(y| \pi(b))& \text{ for } i=n\\
     (2-\bar{b}_{2n-i}) gd_i(y| \pi(b))-\frac{1}{1-b_{2n+1-i}}gd_{i-1}(y| \pi(b))& \text{ for } n+1\leq i\leq 2n-2\\
      (2-{b}_{1}) gd_{2n-1}(y| \pi(b))-\frac{1}{1-b_{2}}gd_{2n-2}(y| \pi(b))& \text{ for } i= 2n-1.
\end{cases}\]
\end{cor}
\begin{proof}
For $n\leq i\leq 2n-2$ this follows from Theorem \ref{thm_sm_cl_b}. 
To prove the case for $i=2n-1$, note that for any function $f$ whose equivariant parameters depend only on \(b_1,\ldots,b_{2n}\), then
\[
s_0^{\B{}
}\pi(f)=\pi(s_0^\infty s_{2n-1}^\infty f),
\]
where both actions $s_0^{\B{}}$ and $s_0^\infty$ include their
\(\Omega(b_1| y)\Omega(b_2| y)\) multiplier.
Consequently, if $s_0^\infty f=f$,
then, since \(s_0^\infty\) and \(s_{2n-1}^\infty\) commute, we have
\begin{equation}\label{eq_com}
    T_0^{\B{}}\pi(f)=\pi(T_{2n-1}^\infty f).
\end{equation}

The denominators agree because
\[
\pi(b_{2n}\ominus b_{2n-1})
=
b_1\ominus\bar b_2
=
b_1\oplus b_2.
\]

For \(m\ge2\), writing \(w_{(m)}=s_m\cdots s_2s_0\), the
braid relations give \(s_0w_{(m)}=w_{(m)}s_2\). Hence
\[
T_0^\infty gd_m(y| b)
=T_{w_{(m)}}^\infty T_2^\infty(1)=0,
\qquad
s_0^\infty gd_m(y| b)=gd_m(y| b).
\]

Thus using \eqref{eq_com},
\[
T_0^{\B{}} gd_{2n-2}(y| \pi(b))=gd_{2n-1}(y| \pi(b)),
\]
and, when \(n\ge3\),
\[
T_0^{\B{}} gd_{2n-3}(y| b)=0.
\]
Now the proof follows by applying \(T_0^B\) to the preceding one-row formula of $gd_{2n-2}(y| b)$ and using the Leibniz rule:
\[
T_0^{\B{}}\Big((2-\bar{b}_2)gd_{2n-2}(y| \pi(b)\Big)
=s_\theta(2-\bar{b}_2)\,T_0^{\B{}}(gd_{2n-2}(y| \pi(b)))
+\frac{s_\theta(2-\bar{b}_2)-(2-\bar{b}_2)}{b_1\oplus b_2}\,gd_{2n-2}(y| \pi(b)).
\]
For the coefficient, note that \(s_\theta\) satisfies
\[
s_\theta(2-\bar b_2)=2-b_1 \text{ and }
\frac{s_\theta(2-\bar b_2)-(2-\bar b_2)}
     {b_1\oplus b_2}
=
-\frac1{1-b_2}.
\]
\end{proof}

\begin{cor}
Let $w\in \WBgr{n}$ with $\ell(w)\le 2n-2$ and $\lambda=\lambda(w)\in \Y^n_{\B{}}$ be the corresponding partition. If $\lambda_1\ge n$ and $\lambda'$ is a strict partition then
\[ gq_{w}^{(n)}(y)=2gd_{\lambda(w)}(y)-\chi(\lambda_1> n)gd_{\lambda'(w)}(y)
\]
\end{cor}
\begin{proof}
    The proof follows from Theorem \ref{thm_sm_cl_b} and using the augmentation homomorphism  $\epsilon_{\sR}$ sending $b_i$ to 0.
\end{proof}

\begin{remark}  
For all $\lambda$ satisfying the assumption of Theorem \ref{thm_sm_cl_b}, and Corollary \ref{cor_per_cl_k}, we have $\lambda_1 \le 2n-1$. Thus using Lemma \ref{lem_1}, equivariant parameters of $gd_\lambda(y\mid b)$ depends only on $b_1,\dots,b_{2n}$. Thus $gd_{\lambda}(y\mid \pi(b))$ is well defined. Similar argument works for the map $\pi'$ defined for type $\D{}$ in the Section \ref{sec_D}. 
\end{remark}

\subsection{Type D}\label{sec_D}
Let $\pi':\Z[b_i,\bar{b}_i:i\geq 1]\to \Z[b_1,\bar{b}_1,\dots, b_{n-1},\bar{b}_{n-1}]$ be the map defined by  

\[\pi'(b_i)=
\begin{cases}
   b_i& \text{for } 1\leq i\leq n-1\\
   0& \text{ if } i=n,2n\\
   \bar{b}_{2n-i} & \text{for } n+1\leq i\leq 2n-2\\
   {b}_1 & \text{for }  i=2n-1\\
\end{cases}
\]

We write  \(\alpha_j^{\D{}},s_j^{\D{}},T_j^{\D{}}\) for the type \(\D{n}^{(1)}\) simple roots and operators for $1\le j\le n$.

\begin{lem}\label{lem_sm_cl_type_D}
For a strict partition
$\la=(n-1,\la_2,\ldots,\la_\ell)$ with $\lambda_2\leq n-3$, we have
\begin{align}\label{eq_gd_bar}
&(2-\bar{b}_{n-1})
gd_{\la}(y| b_1,\ldots,b_{n-2},b_{n-1},\bar{b}_{n-1})\nonumber\\
&=gd_{\la}(y| b_1,\ldots,b_{n-2},0,\bar{b}_{n-1})
+gd_{\la}(y| b_1,\ldots,b_{n-2},b_{n-1},0).
\end{align}
\end{lem}
\begin{proof}
We prove for $\lambda=(n-1).$
Let $F(u)$ be arbitrary function. We define the sequence of functions  
$\Phi_n^F(u_1,\ldots,u_n)$ inductively by
$\Phi_1^F(u)=F(u)$ and 
    $$\Phi_{n+1}^F(u_1,\ldots,u_n,u_{n+1})=
    \frac{\Phi_n^F(u_1,\ldots,u_{n-1},u_{n+1})-\Phi_n^F(u_1,\ldots,u_{n-1},u_{n})}{u_{n+1}-u_{n}}.$$

If we take
    $F(u)=\dfrac{\Omega(b_1| y)\Omega(u| y)-1}{b_1\oplus u}$
    then
    $$gd_{(1)}(y| b_1,b_2)=F(b_2)=\Phi_1^F(b_2), $$ 
\begin{equation}\label{eq_div_rel}
     gd_{(i)}(y| b_1,\dots,b_{i+1})=(1-b_2)\cdots (1-b_i)\Phi^{F}_{i}(b_2,\ldots,b_{i+1}) \text{ for } i\ge 2.
\end{equation}
   
Moreover, we have \begin{align*}
\bar{b}_{n-1} \Phi_{n}(b_2,\ldots,b_{n-1},0,\bar{b}_{n-1})&=\Phi_{n-1}(b_2,\ldots,b_{n-1},\bar{b}_{n-1})-\Phi_{n-1}(b_2,\ldots,b_{n-1},0).\\
{b}_{n-1} \Phi_{n}(b_2,\ldots,b_{n-1},0,\bar{b}_{n-1})&=\Phi_{n-1}(b_2,\ldots,b_{n-1},\bar{b}_{n-1})-\Phi_{n-1}(b_2,\ldots,b_{n-2},0,\bar{b}_{n-1}).
\end{align*}
We compare the two equalities above and use \eqref{eq_div_rel} to obtain our desired identity.
To extend it to arbitrary $(n-1,\lambda_2,\dots,\lambda_\ell)$, write the length-additive factorization
\[
w_{(n-1,\lambda_2,\dots,\lambda_\ell)}=v\rho_{n-1}^{(0)},
\]
where \(v=\rho_{\lambda_\ell}^{(\ell-1)}\dots\rho_{\lambda_2}^{(1)}\) uses only nodes at most \(\lambda_2\le n-3\), and apply \(T_{v}^\infty\) to the one-row identity.
\end{proof}

\begin{remark}
   Using \eqref{sm_cl_n}, the left hand side of \eqref{eq_gd_bar} can be written as $gq_{\lambda(w)}^{(n-1)}(y| \pi(b))$. Thus if $\la=(n-1,\la_2,\ldots,\la_\ell)$ with $\lambda_2\leq n-3$ then 
$$gq_{\lambda(w)}^{(n-1)}(y| \pi(b))=gd_{\la(w)}^{(n)}(y| b_1,\ldots,b_{n-2},0,\bar{b}_{n-1})
+gd_{\la(w)}^{(n)}(y| \pi'(b)).$$
\end{remark}

\begin{thm}\label{thm_sm_cl_type_D}
Let $w\in \WDgr{n}$ with $\ell(w)\leq 2n-3$ and write $\lambda=\lambda(w)$.
\begin{enumerate}
    \item If $\lambda_1=n$ then
    \begin{align}\label{eq_sm_cl_D1}
        gd_w^{(n)}(y| \pi'(b))
        =
        (2-\bar{b}_{n-1})\,gd_{\lambda(w)}(y|\pi'(b) )-gd_{\lambda'(w)}(y| \pi'(b)).
    \end{align}

    \item If $\lambda_1\geq n+1$ then
    \begin{align}\label{eq_sm_cl_D2}
        \gd{w}{n}(y| \pi'(b))
        =
        (2-\bar{b}_{2n-1-\lambda_1})\,
        gd_{\lambda(w)}(y| \pi'(b))
        -
        \frac{1}{1-b_{2n-\lambda_1}}\,
        gd_{\lambda'(w)}(y| \pi'(b)).
    \end{align}
\end{enumerate}
\end{thm}
\begin{proof}
    We first prove
    \eqref{eq_sm_cl_D1}. If $\lambda_1=n$ and $\ell(w)\le 2n-3$ then $\lambda'$ is a strict partition. Thus by Theorem~\ref{small class type BD}, we have
\[
    gd_{w'}^{(n)}(y| b_1,\dots,b_n)
    =
    gd_{\lambda'(w)}(y| b_1,\dots,b_n).
\]
Applying $T_n^{\D{}}$ to both sides of the above equality and then applying $\pi'$ we get
$$gd_{w}^{(n)}(y| \pi'(b))=\frac{gd_{\la'(w)}(y| b_1,\ldots,b_{n-2},0,\bar{b}_{n-1})
-gd_{\la'(w)}(y| b_1,\ldots,b_{n-2},b_{n-1},0)}{\bar{b}_{n-1}}.$$
  The equality in \eqref{eq_gd_bar} for the strict partition $\lambda'$ can be written as the following:
    \begin{align*}
    &gd_{\la'(w)}(y| b_1,\ldots,b_{n-2},0,\bar{b}_{n-1})
-gd_{\la'(w)}(y| b_1,\ldots,b_{n-2},b_{n-1},0)\\
&=(2-\bar{b}_{n-1})(
gd_{\la'(w)}(y| b_1,\ldots,b_{n-2},b_{n-1},\bar{b}_{n-1})-gd_{\la'(w)}(y| b_1,\dots,b_{n-2},b_{n-1},0))\\
&-\bar{b}_{n-1}gd_{\la'(w)}(y| b_1,\dots,b_{n-1},0).
    \end{align*}
    Dividing both sides by $\bar{b}_{n-1}$ we get \eqref{eq_sm_cl_D1}. Now \eqref{eq_sm_cl_D2} will follow from the same argument as in Theorem \ref{thm_sm_cl_b}, and using the fact 
    $$T_{2n-i-1}^{\mathrm{D}}\pi'(f)=\pi'T_i^{\infty}(f)\text{ for }n+1\leq i\leq 2n-3,$$
 for any function $f$ such that the equivariant parameters of $f$ depend only on $b_1,\dots, b_{i+1}$ and $s_{2n-i-1}^\infty(f)=f$.
   One needs to observe that for $n+1\le i\le 2n-3$, $\pi'(c(\alpha_{i}^{\infty}))=c(\alpha_{2n-i-1}^{\mathrm{D}})$. 
\end{proof}

\begin{cor}
Let $w\in \WDgr{n}$ with $\ell(w)\le 2n-3$ and $\lambda=\lambda(w)\in \Y^n_{\D{}}$ be the corresponding element. If $\lambda_1\ge n$ then $\lambda'$ is a strict partition and
\[ gd_{w}^{(n)}(y)=2gd_{\lambda(w)}(y)-gd_{\lambda'(w)}(y).
\]
\end{cor}

\section{$\kappa$-irreducible affine Grassmannian elements of type $\B{}$ and $\D{}$}\label{S:ka_irr_elt}
In this section, we characterize the \(\kappa\)-irreducible affine Grassmannian elements combinatorially in terms of their associated partitions. More precisely, we derive explicit criteria classifying the \(\kappa\)-irreducible elements in types \(\B{}\) and \(\D{}\).

\begin{defn} 
An element \(w\in W_{\af}^{0}\) is called \(\kappa\)-reducible if there exist \(v\in W_{\af}\) and \(i\in I\) such that \[ w=v\kappa_i \text{ and } \ell(w)=\ell(v)+\ell(\kappa_i). \] Otherwise, \(w\) is called \(\kappa\)-irreducible. \end{defn}

\subsection{$\kappa$-irreducible elements of type $\mathrm{B}_n$}

In type \(\mathrm{B}_n\), we recall \[ \kappa_i=w_{(2n-i)^i}, \qquad 1\leq i\leq n. \]

\begin{thm}\label{thm:kappa-irr-type-B}
Let \(\lambda\in\Y^n_{\B{}}\) and \(m_r=m_r(\lambda)\) denote the
multiplicity of \(r\) in \(\lambda\). Then \(w_\lambda\) is
\(\kappa\)-irreducible if and only if
\[
m_{2n-i}\leq i-1
\qquad\text{for every }1\leq i\leq n.
\]
Consequently, the maximal (with respect to Bruhat order) \(\kappa\)-irreducible element
\(\lambda^{(0)}\in\Y^n_{\B{}}\) is determined by
\[
m_{2n-i}\bigl(\lambda^{(0)}\bigr)=
\begin{cases}
i-1, & 1\leq i\leq n,\\
1,   & n+1\leq i\leq 2n-1.
\end{cases}
\]
\end{thm}

We recall the root lattice for type $\mathrm{B}_n$ is
$Q=\bigoplus_{i=1}^n \mathbb{Z}a_i$. The action of the elements $\rho_j^{(1)}$ and $\rho_j^{(0)}$ on the root lattice $Q$ is given by
$$\rho_j^{(0)}(a_1)=-\rho_j^{(1)}(a_1)=a_i;~ \rho_j^{(1)}(a_2)=-\rho_j^{(0)}(a_2)=a_1;~\rho_j^{(0)}(a_\ell)=\rho_j^{(1)}(a_\ell)=a_{\ell-1} \text{ if } 3\leq \ell\le i.$$

\begin{lem}\label{lem:rectification-B}
Let $2\le i\le n$ and \(j=2n-i\). For every \(j\le k\leq 2n-1\), one has 
\begin{align}\label{eq_kap_irr}
\rho_k^{(i)}\kappa_i=\rho_j^{(i)}\dots\rho_j^{(2)}\rho_j^{(1)}\rho_k^{(0)}, \quad
\kappa_i\rho_k^{(1)}=\rho_k^{(i+1)}\rho_j^{(i)}\dots\rho_j^{(2)}\rho_j^{(1)}.
\end{align}
\end{lem} 

\begin{proof} Set \[ R_i=\rho_j^{(i)}\rho_j^{(i-1)} \cdots\rho_j^{(2)}\rho_j^{(1)}. \] Using $ws_r=s_{w(\alpha_r)}w $ together with the above action on the root lattice, we obtain \[ R_i s_r=s_rR_i \quad (2\leq r\leq i-1) ,\quad R_i s_0=s_{\bar{i}}R_i, \qquad s_{\bar{i}}= \begin{cases} s_0, & \text{if \(i\) is even},\\ s_1, & \text{if \(i\) is odd}. \end{cases} \] The identities in \eqref{eq_kap_irr} are immediate when \(k=j\). Suppose that \(k>j\) and \(k=2n-r\) for some \(2\leq r\leq i-1\). Using \[ \rho_k^{(0)}=s_rs_{r+1}\cdots s_{i-1}\rho_j^{(0)}, \] we obtain 
\begin{align*} R_i\rho_k^{(0)} = R_i(s_rs_{r+1}\cdots s_{i-1}\rho_j^{(0)})= (s_rs_{r+1}\cdots s_{i-1})R_i\rho_j^{(0)}= \rho_k^{(i)} \rho_j^{(i-1)}\cdots\rho_j^{(1)}\rho_j^{(0)}= \rho_k^{(i)}\kappa_i. \end{align*} For \(k=2n-1\), we similarly use \[ \rho_k^{(0)}=s_0s_2\cdots s_{i-1}\rho_j^{(0)} \] to obtain 
\begin{align*} R_i\rho_k^{(0)} = R_i(s_0s_2\cdots s_{i-1}\rho_j^{(0)})= (s_{\bar{i}}s_2\cdots s_{i-1})R_i\rho_j^{(0)}= \rho_k^{(i)}\kappa_i. \end{align*} 
This proves the first identity in \eqref{eq_kap_irr}. The second identity follows by applying the Dynkin diagram automorphism \(\tau_{01}\) to the first identity. 
\end{proof}

\begin{proof}[Proof of Theorem \ref{thm:kappa-irr-type-B}] 
Suppose first that \(m_{2n-i}\geq i\). Then \(\lambda\) contains at least \(i\) consecutive Pieri-factors of size \(2n-i\), which form a \((2n-i)^i\)-block. In the reduced expression for \(w_\lambda\), the Pieri-factors of size larger than $2n-i$ occur to the right of this block. By Lemma~\ref{lem:rectification-B}, these larger Pieri-factors can be moved to the left of the block. Thus $w_\lambda=v\kappa_i$ for some \(v\in W_{\af}\), and the factorization is length-additive.

Conversely, suppose that $ w_\lambda=v\kappa_i \text{ and } \ell(w_\lambda)=\ell(v)+\ell(\kappa_i).$  According to the parity of \(i\), either \(v\) or \(\tau_{01}(v)\) belongs to \(\WBgr{n}\). Again, Lemma \ref{lem:rectification-B} allows the larger segments to pass across the \((2n-i)^i\)-block without removing any of its \(i\) consecutive Pieri-factors. Therefore, the affine partition corresponding to \(w_\lambda\) contains at least \(i\) parts of size \(2n-i\). Thus $ m_{2n-i}\geq i.$ 

The formula for \(\lambda^{(0)}\) follows by saturating these bounds and using the fact that all elements of $\Y^n_{\B{}}$ have at most one part of each size less or equal to $n-1$.
\end{proof}

\begin{example}
For $n=3$, the elements $\lambda\in \Y^3_{\B{}}$ such that $w_\lambda$ is $\kappa$ irreducible are 
\begin{align*}
&\emptyset, (1), (2), (2,1), (3), (3,1), (4), (3,2), (4,1), (3,2,1),(3,3), (4,2), (4,3), (3,3,1), (4,2,1),\\ & (3,3,2),(4,3,1),  (3,3,2,1), (4,3,2), (4,3,2,1), (4,3,3), (4,3,3,1), (4,3,3,2), (4,3,3,2,1).\end{align*}
\end{example}

\subsection{$\kappa$-irreducible elements of type $\mathrm{D}_n$}

In type $\mathrm{D}_n$, we recall 
$$\kappa_i=w_{(2n-i-1)^i} \text{ for } 1\leq i\leq n-2,~~~ \kappa_{n-1}=w_{\delta_{n-1}^+} \text{ and } \kappa_{n}=w_{\delta_{n-1}^-},$$ where $\delta_{n-1}=(n-1,n-2,\dots,2,1)$. For $1\leq i\leq n-2$, set
\[
 S_i(\lambda)
 :=\sum_{r=i}^{n-1}m_r(\lambda)+\min\{m_n(\lambda),i-1\}.
\]

The following Theorem describes all $\kappa$-irreducible affine Grassmannian elements of type $\D{n}$.
\begin{thm}\label{thm:kappa-irr-type-D}
Let \(\lambda\in\Y^n_{\D{}}\) and \(m_r=m_r(\lambda)\) denote the
multiplicity of \(r\) in \(\lambda\). Then \(w_\lambda\) is
\(\kappa\)-irreducible if and only if all the following conditions hold.
\begin{enumerate}
    \item[(A1)] $m_{2n-i-1}\le i-1 \text{ for all } 1\le i\le n-2$, \\
    \item[(A2)] $m_n+m_{n-1}\leq n-2$, \\
    \item[(A3)] For all $1\leq i\le n-2$, if $m_i=1$ then $S_i(\lambda)\le n-2$.
\end{enumerate}

Consequently, the maximal \(\kappa\)-irreducible element
\(\lambda^{(0)}\in\Y^n_{\D{}}\) is determined by
\[m_{2n-i-1}=
  \begin{cases}i-1 & \text{ for } 1\leq i\leq n-1\\0 & \text{ for } i=n\\ 1 & \text{ for } n+1\leq i \leq 2n-2.
  \end{cases}\]   
\end{thm}

We begin with the commutation relations for elements
$\kappa_i$, $1\leq i\leq n-2$.

\begin{lem}\label{lem:rectification-D}
Let $1\le i\le n-2$ and \(j=2n-i-1\). For every \(k\geq j\), one has 
\begin{align}\label{eq_kap_irr_D}
\rho_k^{(i)}\kappa_i=\rho_j^{(i)}\dots\rho_j^{(2)}\rho_j^{(1)}\rho_k^{(0)}, \quad
\kappa_i\rho_k^{(1)}=\rho_k^{(i+1)}\rho_j^{(i)}\dots\rho_j^{(2)}\rho_j^{(1)}.
\end{align}
\end{lem} 
\begin{proof}
    The proof follows using the same argument as in Lemma \ref{lem:rectification-B}.
\end{proof}
 
We next establish the corresponding commutation relations for the two spin
factors. Define

\[
 \Delta^\pm:=
 \rho_1^{(n-1)}\rho_2^{(n-2)}\cdots
 \rho_{n-2}^{(2)}\rho_{n-1}^{\pm(1)}.
\]

\begin{lem}\label{lem_com_ru}
The following are the spin-rectification moves.  For
$n-1\leq k\leq2n-2$,
\begin{align}
 (\pi_n\rho_k^{(0)}\pi_n^{-1})\kappa_n
   &=\Delta^-\rho_k^{(0)},
 \label{eq_rec_id_D_n-2}\\
 (\pi_n\rho_k^{(1)}\pi_n^{-1})\Delta^-
   &=\kappa_n\rho_k^{(1)}.
 \label{eq_rec_id_D_n-3}
\end{align}
When $k=n-1$, then $\rho_k^{(0)}$ and $\rho_k^{(1)}$ in \eqref{eq_rec_id_D_n-2} and \eqref{eq_rec_id_D_n-3} must be
$\rho_{n-1}^{-(0)}$ and $\rho_{n-1}^{-(1)}$, respectively.  For the
opposite color,
\begin{align}
 (\pi_n\rho_{n-1}^{+(0)}\pi_n^{-1})\kappa_n
 &=\rho_2^{(n-2)}\cdots\rho_{n-2}^{(2)}
  \rho_{n-1}^{-(1)}\rho_n^{(0)}.
 \label{eq:spin-opposite}\\
 (\pi_n\rho_{n-1}^{+(1)}\pi_n^{-1})\Delta^-
 &=\rho_2^{(n-3)}\dots\rho_{n-2}^{(1)}\rho_{n-1}^{-(0)}\rho_n^{(1)}.
 \label{eq:spin-opposite_1}
\end{align}
Finally, for $1\leq k\leq n-2$,
\begin{equation}
 (\pi_n\rho_k^{(0)}\pi_n^{-1})\kappa_n
 =\rho_1^{(n-2)}\cdots\rho_{n-k-2}^{(k+1)}
  \rho_{n-k}^{(k)}\cdots\rho_{n-2}^{(2)}
  \rho_{n-1}^{-(1)}\rho_{n-1}^{-(0)},
 \label{eq:spin-small}
\end{equation}
The corresponding identities for
$\kappa_{n-1}$ involving $\Delta^+$ also hold.
\end{lem}

\begin{proof}
\eqref{eq_rec_id_D_n-3} is obtained after evaluating the Dynkin automorphism $\pi_1$ in \eqref{eq_rec_id_D_n-2}. To verify \eqref{eq_rec_id_D_n-2} 
\begin{equation}\label{eq_com_k}
w_{(k,n-1,\dots,1)^-}= 
\begin{cases} 
(s_{\epsilon} s_2 \dots s_{n-2} s_{n})\kappa_{n} & k=n-1, \\ 
(s_{k-n+1} \dots s_2 s_1 s_0 s_2 \dots s_{n-2} s_{n})\kappa_{n} & n \le k \le 2n-3, \\ 
(s_{n-1} s_{n-2}  \dots s_2s_1 s_0 s_2 \dots s_{n-2} s_{n})\kappa_{n} & k=2n-2, 
\end{cases}   
\end{equation}
where $\epsilon$ is $0$ or, 1 depending on $n$ is odd or, even. 

 The identity in \eqref{eq:spin-opposite_1} follows from \eqref{eq:spin-opposite} after applying $\pi_1$. The equality in \eqref{eq:spin-opposite} follows similarly using \eqref{eq_com_k} for $k=n-1$. The only change is we need to replace $\epsilon$ by $1$ or, 0  according as $n$ is odd or, even. Then the left hand side of \eqref{eq:spin-opposite}  becomes $s_1s_2s_1\rho_3^{(n-3)}\dots\rho_{n-2}^{(2)}\rho_{n-1}^{-(1)}\rho_{n-1}^{-(0)}$ if $n$ is odd or, $s_0s_2s_0\rho_3^{(n-3)}\dots\rho_{n-2}^{(2)}\rho_{n-1}^{-(1)}\rho_{n-1}^{-(0)}$ if $n$ is even. Next we use $s_is_2s_i=s_2s_is_2$ for $i=0,1$ and $s_{i-1}\rho_i^{(j)}=\rho_i^{(j)}s_i$ for $3\leq i\le n-2$. Also  $s_{n-2}\rho_{n-1}^{-(1)}=\rho_{n-1}^{-(1)}s_n$ and $s_n\rho_{n-1}^{-(0)}=\rho_n^{(0)}$.

Finally,
$\pi_n\rho_k^{(0)}\pi_n^{-1}=s_{n-k}\cdots s_{n-2}s_n$ for
$k\leq n-2$; commuting this word through the staircase gives
\eqref{eq:spin-small}.  Applying the Dynkin automorphism $\tau_{n-1,n}$ gives the identities for $\kappa_{n-1}$.
\end{proof}

\begin{lem}\label{lem_m_n_m_n1}
  Let $\lambda\in\Y^n_{\D{}}$. Suppose that either  $m_n+m_{n-1}\ge n-1$ or, there exists $i$, with $1\leq i\leq n-2$,   such that $m_i=1$ and  $S_i(\lambda)\geq n-1,$ then $w_\lambda$ admits a reduced expression such that one of $\kappa_n$, $\kappa_{n-1}$ appears as a right factor.
\end{lem}
\begin{proof}
 A direct computation gives, for all $2\le k\le n-1 $
\begin{align}
    \rho_{n-1}^{\pm(n-1)}\dots \rho_{n-1}^{\pm(k)}\rho_n^{(k-1)}\cdots \rho_n^{(1)}&=\rho_{k-1}^{(n-1)}\rho_{n-1}^{\pm(n-1)}\rho_{n-1}^{\pm(n-2)}\cdots \rho_{n-1}^{\pm(1)}\label{eq_n_n1}\\ 
  \rho_{n}^{(n-1)}\rho_n^{(n-2)}\cdots \rho_n^{(1)}&=\rho_{n-1}^{\mp(n-1)}\rho_{n-1}^{\pm(n-1)}\rho_{n-1}^{\pm(n-2)}\cdots \rho_{n-1}^{\pm(1)} \label{eq_n_n2}.\\  
\rho_{n-1}^{+(n-1)}\cdots\rho_{n-1}^{+(1)}&=(s_{n-\epsilon}s_{n-2}\cdots s_3s_2)\dots(s_{n-1}s_{n-2})s_{n}\rho_1^{(n-1)}\cdots\rho_{n-2}^{(2)}\rho_{n-1}^{+(1)}\label{eq_n_n3}.
\end{align}
where $\epsilon$ is $1$ or $0$ according as $n$ is even or odd. Thus if $m_n+m_{n-1}\geq n-1$ the desired result follows using \eqref{eq_n_n1}, \eqref{eq_n_n2} and \eqref{eq_n_n3}.

Now assume that $S_i(\lambda)\geq n-1$ with $m_i=1$ for some $1\leq i\le n-2$. If $i=1$ then the condition reduces to $\sum_{j=1}^{n-1}m_j\geq n-1$. If $m_i=1$ for all $1\leq i\le n-1$ then the results follows directly. Otherwise, there are exactly $k$ elements $m_\ell$ with $1\leq \ell\le n-2$ such that $m_\ell=0$. Then $\sum_{j=1}^{n-2}m_j= n-2-k$ and $m_{n-1}\ge k+1$. In this case using \eqref{eq_n_n3}, $k+1$ many Pieri-factors of length $n-1$ can be written as $w{\cdot}\rho_{n-k-1}^{(k+1)}\cdots\rho_{n-2}^{(2)}\rho_{n-1}^{+(1)}$. The remaining $n-k-2$ Pieri factors of length $j$ with $1\leq j \le n-2$ such that $m_j=1$ can be written as $w'{\cdot}\rho_1^{(n-1)}\dots\rho_{n-k-2}^{(k+2)}$. Moreover, element $w$ can be moved to the left of $\rho_1^{(n-1)}\dots\rho_{n-k-2}^{(k+2)}$ using the Coxeter relation and one of $\kappa_n$, $\kappa_{n-1}$, $\Delta^+$ and $\Delta^-$ arise as a factor of $w_\lambda$.
If $i>1$ then we have $\sum_{j=i}^{n-2}m_j=n-i-1-k$ for some $0\le k\le n-i-2$. Now $S_i(\lambda)\ge n-1$ implies $m_{n-1}\ge k+1$ and $m_{n-1}+m_n\ge i+k$. The same argument as in the previous case applies here as well. 

Now using Lemma \ref{lem_com_ru}, we can move the larger Pieri-factors in $w_\lambda$ to the left, leaving the $\kappa_n$ or, $\kappa_{n-1}$ as a right factor.
\end{proof}

\begin{prop}\label{prop:count}
The number of signed partitions $\lambda\in \mathscr{P}_n^{\mathrm D}$
satisfying conditions {\rm(A1)--(A3)} of Theorem~\ref{thm:kappa-irr-type-D}
is $2^{n-3}n!$.
\end{prop}

\begin{proof}
Put $N=n-2$ and fix $t=m_n(\lambda)$ and
$J=\{i\le N:m_i(\lambda)=1\}$. Set
$k=|J\cap\{t+1,\ldots,N\}|$.
Conditions (A2) and (A3) are equivalent to
$
m_{n-1}(\lambda)+t+k\le N.
$
For fixed $t$ and $J$, there are therefore
$1+2(N-t-k)$ choices for $m_{n-1}$ and its sign.
Since
$
\sum_{J\subseteq\{1,\ldots,N\}}k
=(N-t)2^{N-1},
$
the number of choices for the parts of sizes
at most $n$, with $t$ fixed, is
\[
\sum_J\bigl(1+2(N-t-k)\bigr)
=2^N(N-t+1).
\]
Finally, (A1) independently contributes
$(n-2)!$ choices. Hence the total number is
\[
(n-2)!\,2^N\sum_{t=0}^N(N-t+1)
=2^{n-3}n!.
\]
\end{proof}

\begin{proof}[Proof of Theorem \ref{thm:kappa-irr-type-D}]
For \(1\leq i\leq n-2\), reducibility by \(\kappa_i\) is equivalent to the presence of at least \(i\) Pieri-factors of size \(2n-i-1\). Indeed, the argument is identical to that of Theorem~\ref{thm:kappa-irr-type-B}. The commutation relation as in Lemma \ref{lem:rectification-D} move all larger Pieri-factors to the left of the \((2n-i-1)^i\)-block. Thus $ w_\lambda=v\kappa_i $ is a length-additive factorization if and only if \(m_{2n-i-1}\geq i\). Thus $w_\lambda$ is not reducible by any of $\kappa_1,\ldots,\kappa_{n-2}$ precisely when condition~$(A1)$ holds.

{
Next suppose that either (A2) or (A3) fails.
By Lemma~\ref{lem_m_n_m_n1}, one of $\kappa_n$ or $\kappa_{n-1}$ is a reduced right factor
of $w_\lambda$.
Hence $w_\lambda$ is $\kappa$-reducible.
}

{
It remains to prove the converse. By Proposition~\ref{prop:eta-kappa-irreducible},
the number of $\kappa$-irreducible elements in
$W_{\af}^0$ is
\[
\frac{|W(D_n)|}{|\Sigma|}=2^{n-3}n!.
\]
On the other hand, Proposition \ref{prop:count} shows
that exactly $2^{n-3}n!$ signed partitions satisfy
(A1)--(A3). Since every $\kappa$-irreducible element
satisfies these conditions, the two sets coincide.
}
\end{proof}

\begin{example}
 Consider the partition $\lambda=(5,5,4,3)\in \Y^5_{\D{}}$.  Here $n=5$, $m_3=1$ and $S_3(\lambda)=m_3+m_4+m_5=4>n-2$. Thus $\lambda$ is $\kappa$ reducible. But $\mu=(5,5,5,3)$ is $\kappa$-irreducible as one can check $S_3(\mu)=m_3+m_4+2=3\leq n-2$. Also one can check $(A1)$ and $(A2)$ are true for $\mu$.  
\end{example}

\begin{example}
For $n=4$, the elements $\lambda\in \Y^4_{\D{}}$ such that $w_\lambda\in ({{W}}_\af^0)_{\kappa-\text{irr}}$ 
are 
\begin{align*}
&\emptyset, (1), (2), (2,1), (3)^+, (3)^-, (3,1)^{+}, (3,1)^{-}, (4), (3,2)^+, (3,2)^-, (4,1), (5), (5,1),(3,3)^+, (3,3)^-,\\& (4,2), (4,3)^+,(4,3)^-, (5,2), (4,2,1),(4,3,1)^+, (4,3,1)^-,(5,3)^+,(5,3)^-,(4,4),(5,2,1),\\ & (5,3,1)^+,(5,4),(5,3,1)^-,(4,4,1), (5,3,2)^+,(5,3,2)^-,(5,4,1),(4,4,2), (4,4,2,1), (5,3,3)^+, \\&(5,3,3)^-, (5,4,2), (5,4,3)^+, (5,4,3)^-, (5,4,2,1), (5,4,3,1)^+, (5,4,3,1)^-,(5,4,4), (5,4,4,1),\\& (5,4,4,2),(5,4,4,2,1).\end{align*}
\end{example}

\appendix

\section{The $K$-theoretic Peterson isomorphism and applications to the quantum $K$-theory of flag manifolds}\label{Ap_qu_K_th}

Let \(I\) be the set of nodes of the finite Dynkin diagram and let \(W\) be the corresponding finite Weyl group. For \(w\in W\), define its descent set by
\[
\operatorname{Des}(w)
:=
\left\{ i\in I \mid  ws_i<w \right\},
\]
where \(<\) denotes the Bruhat order on \(W\). Associated with \(w\), we define the coweight
\[
\gamma_w
:=
-\sum_{i\in\operatorname{Des}(w)}
\varpi_i^\lor,
\]
where \(\varpi_i^\lor\) denotes the \(i\)-th fundamental coweight. Consider the map
\[
\eta\colon W\longrightarrow \Wafexgr,
\qquad
w\longmapsto wt_{\gamma_w}.
\]

Recall that every element of the extended affine Grassmannian can be expressed in the form
\[
\Wafexgr=\Sigma\cdot\Wafgr.
\]
Let $(\Wafgr)_{\kappa\text{-}\mathrm{irr}}$
denote the set of \(\kappa\)-irreducible elements of \(\Wafgr\).

\begin{prop}\label{prop:eta-kappa-irreducible}
The map \(\eta\) induces a bijection
\[
W
\overset{\sim}{\longrightarrow}
\Sigma\,{\cdot}\,
(\Wafgr)_{\kappa\text{-}\mathrm{irr}}.
\]
\end{prop}
{
\begin{proof}
By \cite[Lemma 2.33]{ISY}, every
$x\in\Wafexgr$ has a unique expression
\[
x=wt_\gamma,\qquad
\gamma=\gamma_w-\sum_{i\in I}a_i\varpi_i^\vee,
\quad a_i\in\Z_{\ge0},
\]
where $w\in W$.

We claim that $\kappa_i$ is a reduced right factor
of $x$ if and only if $a_i>0$. 
Indeed, writing $t_{-\varpi_i^\vee}=\pi_{i^s}^{-1}\kappa_i$,
the length-additivity criterion
\cite[Proposition~2.32]{ISY} gives
\[
\kappa_i\text{ is a reduced right factor of }x
\iff xt_{\varpi_i^\vee}\in\Wafexgr.
\]
By \cite[Lemma 2.33]{ISY}, the latter condition holds exactly when \(\gamma+\varpi_i^\vee-\gamma_w\) is antidominant. This is equivalent to \(a_i>0\), since the coefficient of \(\varpi_i^\vee\) is \(1-a_i\) and all other coefficients are nonpositive.

Writing $x=\sigma u$ with $\sigma\in\Sigma$
and $u\in\Wafgr$, we conclude that
$u$ is $\kappa$-irreducible if and only if
$a_i=0$ for all $i\in I$.
Thus
\[
\Sigma\cdot(\Wafgr)_{\kappa\text{-irr}}
=\{wt_{\gamma_w}\mid w\in W\}.
\]
The assertion follows from the uniqueness
of the decomposition $\Wafex=P^\vee\rtimes W$.
\end{proof}
}

\begin{remark}
Let $i\in I$ be the special node. Then for all $w\in W$ \begin{equation}\eta(u_i^{-1}w)=\pi_i^{-1}\eta(w),\end{equation} see \cite[Lemma 3.3]{IKN}.
\end{remark}

We follow the conventions of \cite{MNS1} for $QK_T(G/B)$, including
the indexing of Schubert classes and the Novikov variables.
For $w\in W$, let $\cOO{w}$ denote the Schubert structure-sheaf class
in $QK_T(G/B)$ indexed by $w$. We define
\[
QK_T^{\mathrm{pol}}(G/B)
:=
\bigoplus_{w\in W} R(T)[Q]\,\cOO{w}.
\]
By the finiteness theorem of Anderson, Chen and Tseng \cite{ACT}, $QK_T^{\mathrm{pol}}(G/B)$ is an
$R(T)[Q]$-subalgebra of $QK_T(G/B)$.

For $i\in I$ define $\sigma_i:=\tgx{t_{-\varpi_i^\lor}}{n}$. The \(K\)-theoretic Peterson map \cite{Kato} together with Theorem \ref{T:Peterson to symm funcs} gives the injective homomorphism of $\RT$-algebras
\[
\Phi\colon
QK_T^{\mathrm{pol}}(G/B) 
\longrightarrow
(\gGtil_{G})_{[
\sigma_i^{-1}]}
\]
satisfying
\[
\Phi\bigl(\cOO{w}\bigr)
=
\frac{\tgx{\eta(w)}{n}}
     {\prod_{i\in \operatorname{Des}(w)} \sigma_i};\qquad \Phi(Q_i)=\prod_{j\in I}\sigma_j^{-\langle\alpha_i^\lor, \alpha_j\rangle}.
\]

Thus, the quantum \(K\)-theory Schubert class indexed by \(w\) is identified with the quotient of the affine Grassmannian Schubert class indexed by
\(\eta(w)=wt_{\gamma_w}\) by the invertible translation class
\(t_{\gamma_w}\).

\begin{remark}
     The map $\eta$ corresponds the longest element $w_0\in W$ to an element in the $\Sigma$ orbit of the maximal  $\kappa$-irreducible element  $\lambda^{(0)}\in (\Wafgr)_{\kappa-\text{irr}}$. 
\end{remark}

\subsection{Computation of the structure constants in $QK(\Spin_5(\C)/B)$.}

Let $\WB{2}=\langle s_1,s_2\rangle$ be the Weyl group of type $B_2$, and $w_0=s_1s_2s_1s_2$
be its longest element.  We write $\cOO{w}$ for the Schubert structure-sheaf
class indexed by $w\in W$ and use $\star$ for the quantum $K$-product. We consider the Schubert basis of
\(QK(\Spin_5/B)\) by
\[
 \bigl\{1,\cOO{s_1},\cOO{s_2},\cOO{s_1s_2},\cOO{s_2s_1},
 \cOO{s_1s_2s_1},\cOO{s_2s_1s_2},\cOO{w_0}\bigr\}.
\]

In the non-equivariant specialization, $\tgq{\pi_1}{2}(y)=1$, and the factor $\pi_1$ may be suppressed. Thus in the non-equivariant setting the two translation classes used in the $K$-theoretic Peterson isomorphism are
\begin{align*}
\sigma_1=\widetilde{gq}_{s_0s_2s_0}^{(2)}(y), \quad
 \sigma_2=\widetilde{gq}_{s_2s_1s_2s_0}^{(2)}(y). 
\end{align*}

The localized affine representatives for the quantum Schubert classes are
\begin{align*}
\cOO{s_1}&=\frac{\widetilde{gq}_{s_2s_0}^{(2)}(y)}{\sigma_1},\quad
\cOO{s_2}=\frac{\widetilde{gq}_{s_1s_2s_0}^{(2)}(y)}{\sigma_2},\quad
\cOO{s_1s_2}=\frac{\widetilde{gq}_{s_2s_0}^{(2)}(y)}{\sigma_2},\quad
\cOO{s_2s_1}=\frac{\widetilde{gq}_{s_0}^{(2)}(y)}{\sigma_1},\\
\cOO{s_1s_2s_1}&=\frac{1}{\sigma_1},\quad
\cOO{s_2s_1s_2}=\frac{\widetilde{gq}_{s_0}^{(2)}(y)}{\sigma_2},\quad
\cOO{w_0}=\frac{\widetilde{gq}^{(2)}_{ s_1s_2s_0}(y)}{\sigma_1\sigma_2}.
\end{align*}
The Novikov variables are given by
\[Q_1=\frac{\sigma_2}{\sigma_1^2},\quad  Q_2=\frac{\sigma_1^2}{\sigma_2^2}.\]

Consequently,
$$\frac{1}{\sigma_2}=Q_1Q_2, \quad \frac{1}{\sigma_1^2}=Q_1^2Q_2. $$


In the affine side using the product between symmetric polynomials we get:
\begin{align*}
\tgq{s_0}{2}(y)\tgq{s_0}{2}(y)&=\tgq{s_2s_0}{2}(y).\\
\tgq{s_0}{2}(y)\tgq{s_2s_0}{2}(y)&=\tgq{s_0s_2s_0}{2}(y)+\tgq{s_1s_2s_0}{2}(y)-\tgq{s_2s_0}{2}(y).\label{eq_P3}\\
\tgq{s_0}{2}(y)\tgq{s_1s_2s_0}{2}(y)&=\sigma_1\tgq{s_0}{2}(y) +\sigma_2- \tgq{s_1s_2s_0}{2}(y).\\
\tgq{s_2s_0}{2}(y)\tgq{s_2s_0}{2}(y)&=2\tgq{s_1s_0s_2s_0}{2}(y)+\sigma_2-2\tgq{s_1s_2s_0}{2}(y)-\sigma_1+\tgq{s_2s_0}{2}(y).\\
  \tgq{s_1s_2s_0}{2}(y)\tgq{s_1s_2s_0}{2}(y)
    &=\sigma_1^2+\sigma_2\tgq{s_2s_0}{2}(y)-\sigma_1\tgq{s_2s_0}{2}(y).\\
\tgq{s_1s_2s_0}{2}(y)\tgq{s_2s_0}{2}(y)
    &= \sigma_1 \tgq{s_2s_0}{2}(y)+\sigma_2\tgq{s_0}{2}(y)-\sigma_1 \tgq{s_0}{2}(y) - \sigma_2+\tgq{s_1s_2s_0}{2}(y).
\end{align*}

Dividing the affine identities by the appropriate translation classes gives the quantum products.  Some
representative calculations are
\begin{align*}
\cOO{s_1}\star \cOO{s_1}
 &=\frac{(\tgq{s_2s_0}{2}(y))^2}{\sigma_1^2}\\
 &=2\cOO{s_2s_1}+Q_1-2Q_1\cOO{s_2}-\cOO{s_1s_2s_1}+Q_1\cOO{s_1s_2},\\[2mm]
\cOO{s_1}\star \cOO{s_2s_1}
 &=\frac{\tgq{s_2s_0}{2}(y)\tgq{s_0}{2}(y)}{\sigma_1^2}
 =\cOO{s_1s_2s_1}+Q_1\cOO{s_2}-Q_1\cOO{s_1s_2},\\[2mm]
\cOO{s_2}\star \cOO{s_2s_1s_2}
 &=\frac{\tgq{s_1s_2s_0}{2}(y)\tgq{s_0}{2}(y)}{\sigma_2^2}
 =Q_2\cOO{s_2s_1}+Q_1Q_2-Q_1Q_2\cOO{s_2},\\[2mm]
\cOO{w_0}\star \cOO{w_0}
 &=\frac{(\tgq{s_1s_2s_0}{2}(y))^2}{\sigma_1^2\sigma_2^2}
 =Q_1^2Q_2^2+Q_1^2Q_2\cOO{s_1s_2}-Q_1^2Q_2^2\cOO{s_1}.
\end{align*}

The quantum $K$-product is commutative, and $1\star X=X$ for every Schubert
class $X$.  Therefore, Table~\ref{tab:B2-multiplication} lists each of the
remaining $28$ unordered products exactly once and hence determines all
structure constants of $QK(\Spin_5(\C)/B)$.

\begingroup
\small
\renewcommand{\arraystretch}{1.22}
\begin{longtable}{@{}>{\raggedright\arraybackslash}p{0.24\textwidth}
                       >{\raggedright\arraybackslash}p{0.70\textwidth}@{}}
\caption{Multiplication table for $QK(\Spin_5(\C)/B)$.}
\label{tab:B2-multiplication}\\
\toprule
Product & Expansion in the Schubert basis \\
\midrule
\endfirsthead
\multicolumn{2}{c}{\tablename~\thetable\ (continued)}\\
\toprule
Product & Expansion in the Schubert basis \\
\midrule
\endhead
\midrule
\multicolumn{2}{r}{Continued on the next page}\\
\endfoot
\bottomrule
\endlastfoot

$\cOO{s_1}\star \cOO{s_1}$
& $Q_1-2Q_1\cOO{s_2}+Q_1\cOO{s_1s_2}+2\cOO{s_2s_1}-\cOO{s_1s_2s_1}$ \\
$\cOO{s_1}\star \cOO{s_2}$
& $\cOO{s_1s_2}+\cOO{s_2s_1}-\cOO{s_1s_2s_1}-\cOO{s_2s_1s_2}+\cOO{w_0}$ \\
$\cOO{s_1}\star \cOO{s_1s_2}$
& $-Q_1Q_2+Q_1Q_2\cOO{s_1}+\cOO{s_1s_2s_1}+2\cOO{s_2s_1s_2}-2\cOO{w_0}$ \\
$\cOO{s_1}\star \cOO{s_2s_1}$
& $\cOO{s_1s_2s_1}+Q_1\cOO{s_2}-Q_1\cOO{s_1s_2}$ \\
$\cOO{s_1}\star \cOO{s_1s_2s_1}$
& $Q_1\cOO{s_1s_2}$ \\
$\cOO{s_1}\star \cOO{s_2s_1s_2}$
& $Q_1Q_2-Q_1Q_2\cOO{s_1}+\cOO{w_0}$ \\
$\cOO{s_1}\star \cOO{w_0}$
& $-Q_1^2Q_2+Q_1Q_2\cOO{s_1}+Q_1^2Q_2\cOO{s_2}-Q_1Q_2\cOO{s_2s_1}+Q_1\cOO{s_2s_1s_2}$ \\
\addlinespace

$\cOO{s_2}\star \cOO{s_2}$
& $Q_2-Q_2\cOO{s_1}+\cOO{s_1s_2}$ \\
$\cOO{s_2}\star \cOO{s_1s_2}$
& $-Q_1Q_2+Q_2\cOO{s_1}+Q_1Q_2\cOO{s_2}-Q_2\cOO{s_2s_1}+\cOO{s_2s_1s_2}$ \\
$\cOO{s_2}\star \cOO{s_2s_1}$
& $\cOO{s_1s_2s_1}+\cOO{s_2s_1s_2}-\cOO{w_0}$ \\
$\cOO{s_2}\star \cOO{s_1s_2s_1}$
& $\cOO{w_0}$ \\
$\cOO{s_2}\star \cOO{s_2s_1s_2}$
& $Q_1Q_2-Q_1Q_2\cOO{s_2}+Q_2\cOO{s_2s_1}$ \\
$\cOO{s_2}\star \cOO{w_0}$
& $Q_1Q_2\cOO{s_1}-Q_1Q_2\cOO{s_1s_2}+Q_2\cOO{s_1s_2s_1}$ \\
\addlinespace

$\cOO{s_1s_2}\star \cOO{s_1s_2}$
& $Q_1Q_2-2Q_1Q_2\cOO{s_2}+Q_1Q_2\cOO{s_1s_2}+2Q_2\cOO{s_2s_1}-Q_2\cOO{s_1s_2s_1}$ \\
$\cOO{s_1s_2}\star \cOO{s_2s_1}$
& $Q_1Q_2-Q_1Q_2\cOO{s_1}+\cOO{w_0}$ \\
$\cOO{s_1s_2}\star \cOO{s_1s_2s_1}$
& $Q_1Q_2\cOO{s_1}$ \\
$\cOO{s_1s_2}\star \cOO{s_2s_1s_2}$
& $Q_1Q_2\cOO{s_2}-Q_1Q_2\cOO{s_1s_2}+Q_2\cOO{s_1s_2s_1}$ \\
$\cOO{s_1s_2}\star \cOO{w_0}$
& $Q_1Q_2\bigl(\cOO{s_1s_2}+\cOO{s_2s_1}-\cOO{s_1s_2s_1}-\cOO{s_2s_1s_2}+\cOO{w_0}\bigr)$ \\
\addlinespace

$\cOO{s_2s_1}\star \cOO{s_2s_1}$
& $Q_1\cOO{s_1s_2}$ \\
$\cOO{s_2s_1}\star \cOO{s_1s_2s_1}$
& $Q_1\cOO{s_2s_1s_2}$ \\
$\cOO{s_2s_1}\star \cOO{s_2s_1s_2}$
& $Q_1Q_2\cOO{s_1}$ \\
$\cOO{s_2s_1}\star \cOO{w_0}$
& $Q_1^2Q_2-Q_1^2Q_2\cOO{s_2}+Q_1Q_2\cOO{s_2s_1}$ \\
\addlinespace

$\cOO{s_1s_2s_1}\star \cOO{s_1s_2s_1}$
& $Q_1^2Q_2$ \\
$\cOO{s_1s_2s_1}\star \cOO{s_2s_1s_2}$
& $Q_1Q_2\cOO{s_2s_1}$ \\
$\cOO{s_1s_2s_1}\star \cOO{w_0}$
& $Q_1^2Q_2\cOO{s_2}$ \\
\addlinespace

$\cOO{s_2s_1s_2}\star \cOO{s_2s_1s_2}$
& $Q_1Q_2\cOO{s_1s_2}$ \\
$\cOO{s_2s_1s_2}\star \cOO{w_0}$
& $Q_1Q_2\bigl(\cOO{s_1s_2s_1}+\cOO{s_2s_1s_2}-\cOO{w_0}\bigr)$ \\
\addlinespace

$\cOO{w_0}\star \cOO{w_0}$
& $Q_1^2Q_2^2-Q_1^2Q_2^2\cOO{s_1}+Q_1^2Q_2\cOO{s_1s_2}$ \\

\end{longtable}
\endgroup

\section{Torus-equivariant homology of affine Grassmannian and relation to Pon's results}\label{S:eq_hom}

In this section, we discuss the torus-equivariant homology of the affine Grassmannian of types B and D. We construct a symmetric-function realization of the aforementioned ring and identify the Schubert classes with explicit symmetric functions. This generalizes Pon's non-equivariant construction \cite{Pon}. The equivariant homology of affine Grassmannian is a Hopf algebra over
\[
S:=H_T^*(\pt)\cong\operatorname{Sym}_{\Z}(Q).
\]
We use the inclusion
$S\hookrightarrow
{S}_{\Q}=\Q[a_1,\ldots,a_n]$
to construct its symmetric-function realization in an ambient
ring with coefficients in $\Sc$.

Let $\hL_{\Sc}$ denote the completion, with respect
to the $y$-degree filtration, of the ring of symmetric functions
in $y=(y_1,y_2,\ldots)$ with coefficients in $\Sc$.
The group $W$ acts on $\hat{\Lambda}_{\Sc}$ through its action on the coefficient ring $\Sc$, leaving the variables $y$ fixed. Define
\begin{equation}
\Omega(a_i| y) =\prod_{j=1}^\infty \frac{1+a_iy_j}{
1-a_i y_j
}.
\end{equation}

For an element $\la=\sum_{i=1}^n \la_i \eps_i$ of the coweight lattice $P^\lor$, define 
\begin{align}
  \eta_{\la} = \prod_{i=1}^n \Omega(a_i| y)^{\la_i} 
\end{align}
By abuse of notation, we use the same symbols as in the $K$-theoretic setting for the corresponding objects in homology.
Here $\lambda_i\in\frac12\Z$.
The fractional powers are defined by the binomial
expansion with constant term $1$.
Thus $\eta_\lambda$ is an invertible element
of $\hat{\Lambda}_{\Sc}$.
Moreover, $\eta_{\lambda+\mu}=\eta_\lambda\eta_\mu$
for all $\lambda,\mu\in P^\vee$.

In both types $\B{n}^{(1)}$ and $\D{n}^{(1)}$,
the natural $W$-action on $\Sc$
extends to an action of the extended affine Weyl group
$\Wafex$ on $\hat{\Lambda}_{\Sc}$
by letting $W$ act on coefficients and
$t_\lambda$ act by multiplication by $\eta_\lambda$:
\begin{equation}\label{eq_act_of_Pc}
t_\lambda f(y| a)=\eta_\lambda f(y| a),
\qquad
f(y| a)\in\hat{\Lambda}_{\Sc}.
\end{equation}

In particular, by \eqref{E:highest root and associated coroot} the action of the simple reflection $s_0=t_{\theta^\lor}s_\theta$ is given by
\begin{align*}
(s_0 f)(y| a) &= (t_{\theta^\lor}s_\theta )f(y| a) 
=\Omega(a_1 |  y)\Omega(a_2 |  y) f(y| s_\theta a)\\
&= \Omega(a_1 |  y)\Omega(a_2 |  y) f(y|  -a_2, -a_1,a_3 \dotsc, a_n),
\end{align*}
for both types $\mathrm{B}_n$ and $\D{n}.$

Define the affine Demazure operators $A_i$ 
for $i\in I_\af$ by 
\begin{equation}
A_i=\alpha_i^{-1}(1-s_i)
\end{equation}
The operators $A_i$ satisfy $A_i^2=0$ and the same Braid relation as the $s_i$. Thus $A_w$ is well defined for $w\in \Waf$.
For an element $\sigma w\in\Wx\cong \Sigma {\ltimes} \Waf$ for $\sigma\in\Sigma$ and $w\in\Waf$, the operators  
$A_{\sigma w} = \sigma A_w$ are well-defined on $\hL_{\Sc}$. Note that the action of $\sigma\in \Sigma$ on $\hL_{\Sc}$ is defined using \eqref{eq_act_pi} and \eqref{eq_act_of_Pc}.

For $w\in \Wxgr$, we define
\begin{equation}
\widehat{x}_{w}^{(n)}(y| a)=A_w (1)
\end{equation}
where we denote $\widehat{x}_{w}^{(n)}(y| a)$ by
${\widehat{B}}_{w}^{(n)}(y| a)$ for type $\B{n}$
and $\widehat{D}_{w}^{(n)}(y| a)$ for type $\D{n}.$ 

Set 
$S_{\sco}:=\operatorname{Sym}_{\Z}(P)$.
Let $\htL{G}$ and $\htK{G}$ denote the {$S$}-submodule
of $\hL_{\Sc}$ defined by
\[
\htL{G}
:=\sum_{w\in\Wafexgr}
S\,\widehat{x}_w^{(n)}(y| a),
\qquad
\htK{G}
:=\sum_{w\in\Wafgr}
S\,\widehat{x}_w^{(n)}(y| a).
\]
\begin{thm}
There is an $S$-algebra isomorphism
\[
H_*^T(\Gr_G)\cong\htL{G},
\]
identifying the Schubert class indexed by $w\in\Wafexgr$
with $\widehat{x}_w^{(n)}(y| a)$.
There is also an $S_{\sco}$-algebra isomorphism
\[
H_*^{\Tsc}(\Gr_{\Gsc})
\cong S_{\sco}\otimes_{S}\htK{G},
\]
identifying the Schubert class indexed by $w\in\Wafgr$
with $1\otimes\widehat{x}_w^{(n)}(y| a)$.
\end{thm}

We also have factorization property of $\widehat{x}_{w}^{(n)}(y| a)$ in terms of $\widehat{x}_{\kappa_i}^{(n)}(y| a)$ similar to Theorem \ref{thm_fac_type_B} and Theorem \ref{thm_fac_type_D}.

\begin{thm}
For $i\in I(\mathrm{B}_n)$, if $v\kappa_i \in \WBexgr{n}$ and 
$\ell(v\kappa_i) = \ell(v) + \ell(\kappa_i)$ then
    \begin{equation}
{\widehat{B}}_{v\kappa_i}^{(n)}(y| a) =
\begin{cases}
{\widehat{B}}_{\kappa_i}^{(n)}(y| a)
 {\widehat{B}}_{v}^{(n)}(y| a) & \text{if $i$ is even} \\
{\widehat{B}}_{\kappa_i}^{(n)}(y| a)
{\widehat{B}}_{\bar{v}}^{(n)}(y| u_1(a)) & \text{if $i$ is odd},
\end{cases}
\end{equation}
where $\bar{v}$ is defined in Theorem \ref{thm_fac_type_B}.
\end{thm}

\begin{thm}
For $i\in I(\mathrm{D}_n)$, if $v\kappa_i \in \WDexgr{n}$ and 
$\ell(v\kappa_i) = \ell(v) + \ell(\kappa_i)$ then
    \begin{equation}
\widehat{D}_{v\kappa_i}^{(n)}(y| a) = 
\begin{cases}
\widehat{D}_{\kappa_i}^{(n)}(y| a)
\,\widehat{D}_{v}^{(n)}(y| a) & 1\le i\le n-2~ \text{if $i$ is even}, \\
\widehat{D}_{\kappa_i}^{(n)}(y| a)
\,\widehat{D}_{\bar{v}}^{(n)}(y| u_1(a)) & 1\le i\le n-2~\text{if $i$ is odd},\\
\widehat{D}_{\kappa_i}^{(n)}(y| a)
\,\widehat{D}_{\bar{v}_{(n-1)}}^{(n)}(y| u_{n-1}(a)) & i=n-1,\\
\widehat{D}_{\kappa_i}^{(n)}(y| a)
\,\widehat{D}_{\bar{v}_{(n)}}^{(n)}(y| u_{n}(a)) & i=n.
\end{cases}
\end{equation}
 where $\bar{v}, \bar{v}_{(n-1)}$ and $\bar{v}_{(n)}$ are defined in Theorem \ref{thm_fac_type_D}.
\end{thm}

\begin{thm} If $\lambda$ contains a part of size $n-1$ then
\begin{align*}
\widehat{D}_{w_{\la^+}}^{(n)}(y| a_1,\dotsc,a_{n-1},a_n)= 
  \widehat{D}_{w_{\la^-}}^{(n)}(y| a_1,\dotsc,a_{n-1},-a_n).
  \end{align*}
Consequently, 
we have
\begin{align*}
\widehat{D}_{w_{\la^+}}^{(n)}(y)= 
  \widehat{D}_{w_{\la^-}}^{(n)}(y).
  \end{align*}
\end{thm}
\begin{proof}
The proof is analogous to Theorem \ref{T:gd auto} and follows by applying $\autn$.
\end{proof}

Next, we show that for small $w\in\Wafgr$ the symmetric functions $\widehat{B}_w^{(n)}(y| a)$ and $\DD^{(n)}(y| a)$ can be written in terms of type $\D{\infty}$ equivariant homology symmetric functions. 
For every $\lambda\in \SP(r)$ and a positive integer $n\geq r$, we recall factorial Schur-$P$ functions $P_\lambda(x_1,\dots,x_n| a_1,a_2,\dots)$ due to Ivanov \cite{Iv}. Let $\lambda\in \SP(r)$ and $n$ be an even integer such that $r\leq n$. Define 
\[D_\lambda(x_1,x_2,\dots,x_{n}| a_1,a_2,\dots)=
    P_\lambda(x_1,x_2,\dots, x_{n}| a_1,a_2,\dots).
\]

We denote the inverse limit of $D_\lambda(x_1,\dots,x_{2m}| a_1,a_2,\dots)$ by $D_\lambda(x| a)$ as $D_\lambda(x_1,\dots,x_{2m-2}| a)=D_\lambda(x_1,\dots,x_{2m-2},0,0| a).$ 
Set \(S_\infty:=\mathbb Z[a_i\mid i\ge1]\).
The type $\D{\infty}$ dual factorial functions
$\widehat{D}_\lambda(y| a)\in\hL_{S_\infty}$
are defined by
\begin{equation}
\prod_{i,j} \frac{1 + x_i y_j}{1 - x_i y_j}
= \sum_{\lambda\in \SP} \widehat{D}_\lambda(y| a) D_\lambda(x|  a).
\end{equation}

\begin{remark}\label{rem_eq}
    We denote ${D}_\lambda(x)={D}_\lambda(x| 0)$ and $\widehat{D}_\lambda(y)=\widehat{D}_\lambda(y| 0)$. Then we can identify the functions $D_\lambda(x)$ and  $\widehat{D}_\lambda(y)$ by $P_{\lambda}(x)$ and $Q_\lambda(y)$ respectively for all $\lambda\in \SP$.
\end{remark}

Let \(\pi\) and \(\pi'\) denote the additive parameter
specializations defined by
\[
\begin{aligned}
\pi(a_1,\ldots,a_{2n})
&=(a_1,\ldots,a_n,-a_n,\ldots,-a_2,a_1),\\
\pi'(a_1,\ldots,a_{2n})
&=(a_1,\ldots,a_{n-1},0,-a_{n-1},\ldots,-a_2,a_1,0).
\end{aligned}
\]
These are the additive analogues of the specializations in
Section~\ref{S: small_cl}, with each multiplicative inverse \(\bar b_i\)
replaced by \(-a_i\).

\begin{thm}\label{small class hom type BD}
Let $w\in \Wafgr$ and write $\lambda(w)=(\lambda_1,\dots,\lambda_r)$.
\begin{enumerate}
\item If $w\in \WBgr{n}$ with $\ell(w)\le 2n-2$ then
\begin{align}\label{E:QB}
    \widehat{B}_w^{(n)}(y| a)=2^{\chi(\lambda_1\ge n)} \widehat{D}_{\lambda(w)}(y| \pi(a)).
\end{align}
\item If $w\in \WDgr{n}$ with $\ell(w)\le 2n-3$ then 
\begin{align}\label{E:QD}
    \DD^{(n)}(y| \pi'(a))=2^{\chi(\lambda_1\ge n)} \widehat{D}_{\lambda(w)}(y| \pi'(a)).
\end{align}
\end{enumerate}
\end{thm}

\begin{cor}\label{cor_per_cl}
The functions $\widehat{B}_w^{(n)}(y)$ and $\DD^{(n)}(y)$ for  Pieri classes are given by:
 \begin{align}\label{E:QB_n}
    {\widehat{B}}_{\rho_i^{(0)}}(y)&=2^{\chi(i\ge n)} Q_{i}(y),\quad \text{ for } 1\le i\le 2n-1\\
    \widehat{D}_{\rho_i^{(0)}}(y)&=2^{\chi(i\ge n)} Q_{i}(y) \quad \text{ for } 1\le i\le 2n-2, i\neq n-1.
\end{align}   
Moreover, $ \widehat{D}_{\rho_{n-1}^{-(0)}}(y)= \widehat{D}_{\rho_{n-1}^{+(0)}}(y)=Q_{n-1}(y)$
\end{cor}
\begin{proof}
 If $i=2n-1$ in type $\B{}$ or, $i=2n-2$ in type $\D{}$ follows from same argument as in Corollary \ref{cor_per_cl_k}. The remaining cases follows from Theorem \ref{small class hom type BD}. 
\end{proof}

Let \(\sigma_{w}^0\) denote the non-equivariant Schubert classes. The Schubert-cell decomposition gives
\[
H_*(\Gr_{\Gsc})
=\bigoplus_{w\in W_{\mathrm{af}}^0}\mathbb Z\sigma_w^0.
\]
Define the ring of $\Z$-span of non-equivariant functions by
\[\Gamma_{\Gsc}=\sum_{w\in \Wafgr}\Z \, \widehat{x}_{w}(y).\]
Then $\Gamma_{\Gsc}$ is closed under multiplication and the specialization map gives a natural surjection 
$$\psi_{\Gsc}: H_*(\Gr_{\Gsc})=\bigoplus_{w\in \Wafgr}
\Z \sigma_{w}^0
\twoheadrightarrow  \Gamma_{\Gsc}$$
sending $\sigma_{w}^0$ to $\widehat{x}_{w}^{(n)}(y)$. Thus we have an isomorphism 
$$\Gamma_{\Gsc}\cong {H_*(\Gr_{\Gsc})}/{\Ker(\psi_{\Gsc})}.$$

\begin{thm}\label{thm_hom_eq}
The family of symmetric functions
${\widehat{B}}_{w}^{(n)}(y)\;(w\in \WBgr{n})$ is linearly independent
over $\Z$. Moreover, $\Gamma_{\Spin_{2n+1}(\C)}$ is isomorphic to the
homology ring $H_*(\Gr_{Spin_{2n+1}(\C)})$
such that the Schubert class $\sigma_{w}^0$
corresponds to ${\widehat{B}}_{w}^{(n)}(y)$.
 \end{thm}
\begin{proof}
The homology Schubert basis 
 is identified with the family of symmetric functions
 $\tilde{G}_w^{B_n}(y)$ introduced by Pon \cite[Section 2.2]{Pon}. Following \cite[Theorem 2.5]{Pon} there exists an algebra isomorphism 
$$\psi_G': H_*(\Gr_{\Spin_{2n+1}(\C)})\to \Z[Q_1,\dots,Q_{n-1},2Q_n,\dots,2Q_{2n-1}],$$
that sends the Schubert class $\sigma_w^0$ to $\tilde{G}_w^{B_n}(y)$. In particular, $\psi_G'(\sigma_{\rho_i^{(0)}})=2^{\chi(i\ge n)}Q_i$. 

Using Remark \ref{rem_eq} and Corollary \ref{cor_per_cl} we also have $\psi_G(\sigma_{\rho_i^{(0)}})=2^{\chi(i\ge n)}Q_i$. Moreover, using \cite[Theorem 2.5]{Pon}, the Pieri classes $\{\sigma_{\rho_i^{(0)}}: 1\le i\le 2n-1\}$ generate the ring $H_*(\Gr_{\Spin_{2n+1}(\C)})$. Thus for an arbitrary $w\in \WBgr{n}$, we have ${\widehat{B}}_w^{(n)}(y)=\tilde{G}_w^{B_n}(y)$. Since the family $\{\tilde{G}_w^{B_n}(y)| w\in \Wafgr\}$ is linearly 
independent, we deduce the linear independence of the family
$\{{\widehat{B}}_{w}^{(n)}(y)| w\in \Wafgr\}.$ This implies that $\psi_{\Gsc}$ is injective and the proof follows.
\end{proof}

In type $\D{n}$, the homomorphism $\psi_{\Gsc}$ is not injective. Let $J_{\Spin}$ be the ideal generated by $(\sigma^0_{\rho_{n-1}^{+(0)}}-\sigma^0_{\rho_{n-1}^{-(0)}})$. Then we have the following conjecture:
\begin{conj}
    We have the following $\Z$-algebra  isomorphism:
    $$\Gamma_{\Spin_{2n}(\C)}\cong {H_*(\Gr_{\Spin_{2n}(\C)})}/{J_{\Spin}}.$$
\end{conj}

\section{Expansion of $\tgq{w}{n}(y)$ and $\tgd{w}{n}(y)$-functions in terms of $\widetilde{gd}_\lambda(y)$}\label{S: expansion}
Consider $\lambda\in \SP$ and let us define $$\widetilde{gd}_{\lambda}(y)=\sum_{\mu\in \SP:~\mu\subseteq\lambda}gd_{\mu}(y).$$

\begin{conj}
 For every affine Grassmannian element \(w\in \Wafgr\), there exist nonnegative integers \(a_{w,\lambda}\) and \(b_{w,\lambda}\) such that
  \begin{align}
      \tgd{w}{n}(y)&=\sum_{\lambda\in \SP:~|\lambda|\leq\ell(w)}(-1)^{\ell(w)-|\lambda|}a_{w,\lambda}\widetilde{gd}_{\lambda}(y),\\
      \tgq{w}{n}(y)&=\sum_{\lambda\in \SP:~|\lambda|\leq \ell(w)}(-1)^{\ell(w)-|\lambda|}b_{w,\lambda}\widetilde{gd}_{\lambda}(y).
  \end{align}
\end{conj}


\subsection{Sign-normalized coefficient tables for type $\mathrm{D}_4$}

The following table records the positive integers $a_{w,\lambda}$ for $\kappa$ irreducible $w$ of type $\D{4}$. 
The table contains one representative modulo the involution interchanging
$s_3$ and $s_4$. A blank entry means $a_{w,\lambda}=0$. For $\ell(w)\leq 13$, the rows are labeled by $\lambda(w)$. From length $14$ onward, the Weyl group element
is displayed in the heading and the table contains only its coefficient row.
The columns are ordered by decreasing
$|\lambda|$ and, for fixed size, by decreasing lexicographic order. Horizontal
rules separate distinct values of $\ell(w)$, and vertical rules separate
distinct values of $|\lambda|$. Columns that vanish throughout a displayed
table are omitted.
\par\medskip
\noindent\begin{minipage}{\linewidth}
\par\medskip\noindent{\large\bfseries Lengths $0\leq\ell(w)\leq10$\par}\nobreak\smallskip
\begin{center}
\fontsize{7.5pt}{9pt}\selectfont
\setlength{\tabcolsep}{1.2pt}
\renewcommand{\arraystretch}{1.05}
\sbox{\AFGDTableBox}{%
\begin{tabular}{r|l|rrrrrrr|rrrrrr|rrrrr|rrrr|rrrr|rrr|rr|rr|r|r|r}
\toprule
$\ell(w)$ & $\lambda(w)$ & \multicolumn{7}{c|}{$|\lambda|=10$} & \multicolumn{6}{c|}{$|\lambda|=9$} & \multicolumn{5}{c|}{$|\lambda|=8$} & \multicolumn{4}{c|}{$|\lambda|=7$} & \multicolumn{4}{c|}{$|\lambda|=6$} & \multicolumn{3}{c|}{$|\lambda|=5$} & \multicolumn{2}{c|}{$|\lambda|=4$} & \multicolumn{2}{c|}{$|\lambda|=3$} & \multicolumn{1}{c|}{$|\lambda|=2$} & \multicolumn{1}{c|}{$|\lambda|=1$} & \multicolumn{1}{c}{$|\lambda|=0$} \\
 &  & \rotatebox{90}{$(8,2)$} & \rotatebox{90}{$(7,3)$} & \rotatebox{90}{$(7,2,1)$} & \rotatebox{90}{$(6,4)$} & \rotatebox{90}{$(6,3,1)$} & \rotatebox{90}{$(5,4,1)$} & \rotatebox{90}{$(5,3,2)$} & \rotatebox{90}{$(7,2)$} & \rotatebox{90}{$(6,3)$} & \rotatebox{90}{$(6,2,1)$} & \rotatebox{90}{$(5,4)$} & \rotatebox{90}{$(5,3,1)$} & \rotatebox{90}{$(4,3,2)$} & \rotatebox{90}{$(7,1)$} & \rotatebox{90}{$(6,2)$} & \rotatebox{90}{$(5,3)$} & \rotatebox{90}{$(5,2,1)$} & \rotatebox{90}{$(4,3,1)$} & \rotatebox{90}{$(6,1)$} & \rotatebox{90}{$(5,2)$} & \rotatebox{90}{$(4,3)$} & \rotatebox{90}{$(4,2,1)$} & \rotatebox{90}{$(6)$} & \rotatebox{90}{$(5,1)$} & \rotatebox{90}{$(4,2)$} & \rotatebox{90}{$(3,2,1)$} & \rotatebox{90}{$(5)$} & \rotatebox{90}{$(4,1)$} & \rotatebox{90}{$(3,2)$} & \rotatebox{90}{$(4)$} & \rotatebox{90}{$(3,1)$} & \rotatebox{90}{$(3)$} & \rotatebox{90}{$(2,1)$} & \rotatebox{90}{$(2)$} & \rotatebox{90}{$(1)$} & \rotatebox{90}{$\varnothing$} \\
\midrule
0 & $\varnothing$ &  &  &  &  &  &  &  &  &  &  &  &  &  &  &  &  &  &  &  &  &  &  &  &  &  &  &  &  &  &  &  &  &  &  &  & 1 \\
\midrule
1 & $(1)$ &  &  &  &  &  &  &  &  &  &  &  &  &  &  &  &  &  &  &  &  &  &  &  &  &  &  &  &  &  &  &  &  &  &  & 1 &  \\
\midrule
2 & $(2)$ &  &  &  &  &  &  &  &  &  &  &  &  &  &  &  &  &  &  &  &  &  &  &  &  &  &  &  &  &  &  &  &  &  & 1 &  &  \\
\midrule
3 & $(2,1)$ &  &  &  &  &  &  &  &  &  &  &  &  &  &  &  &  &  &  &  &  &  &  &  &  &  &  &  &  &  &  &  &  & 1 &  &  &  \\
3 & $(3)$ &  &  &  &  &  &  &  &  &  &  &  &  &  &  &  &  &  &  &  &  &  &  &  &  &  &  &  &  &  &  &  & 1 &  &  &  &  \\
\midrule
4 & $(3,1)$ &  &  &  &  &  &  &  &  &  &  &  &  &  &  &  &  &  &  &  &  &  &  &  &  &  &  &  &  &  &  & 1 &  &  &  &  &  \\
4 & $(4)$ &  &  &  &  &  &  &  &  &  &  &  &  &  &  &  &  &  &  &  &  &  &  &  &  &  &  &  &  &  & 2 &  & 1 &  &  &  &  \\
\midrule
5 & $(4,1)$ &  &  &  &  &  &  &  &  &  &  &  &  &  &  &  &  &  &  &  &  &  &  &  &  &  &  &  & 2 &  &  & 1 &  &  &  &  &  \\
5 & $(3,2)$ &  &  &  &  &  &  &  &  &  &  &  &  &  &  &  &  &  &  &  &  &  &  &  &  &  &  &  &  & 1 &  &  &  &  &  &  &  \\
5 & $(5)$ &  &  &  &  &  &  &  &  &  &  &  &  &  &  &  &  &  &  &  &  &  &  &  &  &  &  & 2 &  &  & 1 &  &  &  &  &  &  \\
\midrule
6 & $(5,1)$ &  &  &  &  &  &  &  &  &  &  &  &  &  &  &  &  &  &  &  &  &  &  & 2 & 2 &  &  & 3 & 1 &  & 1 &  &  &  &  &  &  \\
6 & $(4,2)$ &  &  &  &  &  &  &  &  &  &  &  &  &  &  &  &  &  &  &  &  &  &  &  & 2 & 2 &  &  & 3 & 1 &  & 1 &  &  &  &  &  \\
6 & $(3,3)$ &  &  &  &  &  &  &  &  &  &  &  &  &  &  &  &  &  &  &  &  &  &  &  &  & 2 &  &  &  & 1 &  &  &  &  &  &  &  \\
\midrule
7 & $(4,2,1)$ &  &  &  &  &  &  &  &  &  &  &  &  &  &  &  &  &  &  & 2 & 2 &  & 2 &  & 3 & 3 & 1 &  & 1 & 1 &  &  &  &  &  &  &  \\
7 & $(5,2)$ &  &  &  &  &  &  &  &  &  &  &  &  &  &  &  &  &  &  & 2 & 2 &  &  &  & 3 & 1 &  &  & 1 &  &  &  &  &  &  &  &  \\
7 & $(4,3)$ &  &  &  &  &  &  &  &  &  &  &  &  &  &  &  &  &  &  &  & 2 & 2 &  &  &  & 3 &  &  &  &  &  &  &  &  &  &  &  \\
\midrule
8 & $(5,2,1)$ &  &  &  &  &  &  &  &  &  &  &  &  &  & 4 & 4 &  & 2 &  & 8 & 6 &  & 1 &  & 5 & 2 &  &  & 1 &  &  &  &  &  &  &  &  \\
8 & $(4,3,1)$ &  &  &  &  &  &  &  &  &  &  &  &  &  &  & 2 & 2 & 2 & 2 &  & 5 & 3 & 3 &  &  & 4 &  &  &  &  &  &  &  &  &  &  &  \\
8 & $(5,3)$ &  &  &  &  &  &  &  &  &  &  &  &  &  &  & 2 & 2 &  &  &  & 3 & 1 &  &  &  & 1 &  &  &  &  &  &  &  &  &  &  &  \\
8 & $(4,4)$ &  &  &  &  &  &  &  &  &  &  &  &  &  &  &  & 4 &  &  &  & 2 & 2 &  &  &  & 1 &  &  &  &  &  &  &  &  &  &  &  \\
\midrule
9 & $(5,3,1)$ &  &  &  &  &  &  &  & 4 & 4 & 2 &  & 2 &  &  & 10 & 6 & 3 & 1 &  & 8 & 2 & 1 &  &  & 2 &  &  &  &  &  &  &  &  &  &  &  \\
9 & $(4,4,1)$ &  &  &  &  &  &  &  &  & 4 &  & 4 & 4 &  &  & 2 & 12 & 2 & 2 &  & 5 & 3 & 1 &  &  & 2 &  &  &  &  &  &  &  &  &  &  &  \\
9 & $(5,4)$ &  &  &  &  &  &  &  &  & 4 &  & 4 &  &  &  & 2 & 8 &  &  &  & 3 & 1 &  &  &  & 1 &  &  &  &  &  &  &  &  &  &  &  \\
\midrule
10 & $(5,4,1)$ &  & 8 &  & 8 & 4 & 4 &  & 4 & 24 & 2 & 12 & 8 &  &  & 10 & 22 & 3 & 1 &  & 8 & 2 & 1 &  &  & 2 &  &  &  &  &  &  &  &  &  &  &  \\
10 & $(5,3,2)$ & 4 & 4 & 4 &  & 4 &  & 2 & 12 & 8 & 8 &  & 6 & 1 &  & 13 & 5 & 5 & 2 &  & 6 & 1 & 1 &  &  & 1 &  &  &  &  &  &  &  &  &  &  &  \\
10 & $(4,4,2)$ &  &  &  & 4 & 4 & 4 & 4 &  & 6 & 2 & 6 & 12 & 2 &  & 2 & 9 & 3 & 3 &  & 3 & 1 & 1 &  &  & 1 &  &  &  &  &  &  &  &  &  &  &  \\
\bottomrule
\end{tabular}
}%
\usebox{\AFGDTableBox}%
\end{center}
\end{minipage}
\newgeometry{top=15mm,bottom=15mm,left=18mm,right=18mm}
\begin{landscape}

\par\medskip
\noindent\begin{minipage}{\linewidth}
\par\medskip\noindent{\large\bfseries Lengths $11\leq\ell(w)\leq12$\par}\nobreak\smallskip
\begin{center}
\fontsize{7.5pt}{9pt}\selectfont
\setlength{\tabcolsep}{1.2pt}
\renewcommand{\arraystretch}{1.05}
\sbox{\AFGDTableBox}{%
\begin{tabular}{r|l|rrrrrrrrrrr|rrrrrrrrrr|rrrrrrrr|rrrrrr|rrrr|rrr|r}
\toprule
$\ell(w)$ & $w$ & \multicolumn{11}{c|}{$|\lambda|=12$} & \multicolumn{10}{c|}{$|\lambda|=11$} & \multicolumn{8}{c|}{$|\lambda|=10$} & \multicolumn{6}{c|}{$|\lambda|=9$} & \multicolumn{4}{c|}{$|\lambda|=8$} & \multicolumn{3}{c|}{$|\lambda|=7$} & \multicolumn{1}{c}{$|\lambda|=6$} \\
 &  & \rotatebox{90}{$(9,3)$} & \rotatebox{90}{$(8,4)$} & \rotatebox{90}{$(8,3,1)$} & \rotatebox{90}{$(7,5)$} & \rotatebox{90}{$(7,4,1)$} & \rotatebox{90}{$(7,3,2)$} & \rotatebox{90}{$(6,5,1)$} & \rotatebox{90}{$(6,4,2)$} & \rotatebox{90}{$(6,3,2,1)$} & \rotatebox{90}{$(5,4,3)$} & \rotatebox{90}{$(5,4,2,1)$} & \rotatebox{90}{$(9,2)$} & \rotatebox{90}{$(8,3)$} & \rotatebox{90}{$(8,2,1)$} & \rotatebox{90}{$(7,4)$} & \rotatebox{90}{$(7,3,1)$} & \rotatebox{90}{$(6,5)$} & \rotatebox{90}{$(6,4,1)$} & \rotatebox{90}{$(6,3,2)$} & \rotatebox{90}{$(5,4,2)$} & \rotatebox{90}{$(5,3,2,1)$} & \rotatebox{90}{$(8,2)$} & \rotatebox{90}{$(7,3)$} & \rotatebox{90}{$(7,2,1)$} & \rotatebox{90}{$(6,4)$} & \rotatebox{90}{$(6,3,1)$} & \rotatebox{90}{$(5,4,1)$} & \rotatebox{90}{$(5,3,2)$} & \rotatebox{90}{$(4,3,2,1)$} & \rotatebox{90}{$(7,2)$} & \rotatebox{90}{$(6,3)$} & \rotatebox{90}{$(6,2,1)$} & \rotatebox{90}{$(5,4)$} & \rotatebox{90}{$(5,3,1)$} & \rotatebox{90}{$(4,3,2)$} & \rotatebox{90}{$(6,2)$} & \rotatebox{90}{$(5,3)$} & \rotatebox{90}{$(5,2,1)$} & \rotatebox{90}{$(4,3,1)$} & \rotatebox{90}{$(5,2)$} & \rotatebox{90}{$(4,3)$} & \rotatebox{90}{$(4,2,1)$} & \rotatebox{90}{$(4,2)$} \\
\midrule
11 & $(4,4,2,1)$ &  &  &  &  &  &  &  &  &  &  &  &  &  &  & 8 & 8 & 8 & 12 & 8 & 8 & 4 &  & 12 & 4 & 32 & 34 & 22 & 24 & 2 & 4 & 34 & 10 & 18 & 43 & 6 & 10 & 26 & 8 & 6 & 8 & 2 & 2 & 2 \\
11 & $(5,4,2)$ &  &  &  &  &  &  &  &  &  &  &  &  & 8 &  & 16 & 8 & 8 & 12 & 4 & 4 &  & 4 & 40 & 4 & 48 & 30 & 18 & 8 &  & 16 & 68 & 10 & 28 & 31 & 1 & 23 & 43 & 8 & 3 & 14 & 3 & 2 & 3 \\
11 & $(5,3,3)$ &  &  &  &  &  &  &  &  &  &  &  & 4 & 12 & 4 & 8 & 12 &  & 8 & 6 & 4 &  & 20 & 40 & 16 & 16 & 34 & 12 & 11 &  & 37 & 47 & 21 & 10 & 30 & 2 & 32 & 23 & 11 & 4 & 13 & 2 & 2 & 2 \\
\midrule
12 & $(5,4,2,1)$ &  & 8 & 8 & 16 & 16 & 8 & 8 & 12 & 4 &  & 4 &  & 12 & 4 & 48 & 40 & 28 & 48 & 30 & 18 & 8 & 4 & 44 & 12 & 68 & 68 & 28 & 31 & 1 & 12 & 53 & 13 & 17 & 43 & 3 & 13 & 23 & 6 & 3 & 6 & 1 & 1 & 1 \\
12 & $(5,4,3)$ & 8 & 16 & 8 & 8 & 16 & 8 & 8 & 12 &  & 4 &  & 4 & 40 & 4 & 52 & 40 & 16 & 48 & 24 & 18 &  & 16 & 74 & 12 & 60 & 66 & 28 & 20 &  & 25 & 64 & 13 & 17 & 40 & 1 & 19 & 24 & 6 & 2 & 7 & 1 & 1 & 1 \\
\bottomrule
\end{tabular}
}%
\usebox{\AFGDTableBox}%
\end{center}
\end{minipage}

\par\medskip
\noindent\begin{minipage}{\linewidth}
\par\medskip\noindent{\large\bfseries Length $\ell(w)=13$\par}\nobreak\smallskip
\begin{center}
\fontsize{7.5pt}{9pt}\selectfont
\setlength{\tabcolsep}{0.6pt}
\renewcommand{\arraystretch}{1.05}
\sbox{\AFGDTableBox}{%
\begin{tabular}{r|l|rrrrrrrrrrrrrr|rrrrrrrrrrrrr|rrrrrrrrrr|rrrrrrrr|rrrrrr|rrrr|rrr|r}
\toprule
$\ell(w)$ & $w$ & \multicolumn{14}{c|}{$|\lambda|=13$} & \multicolumn{13}{c|}{$|\lambda|=12$} & \multicolumn{10}{c|}{$|\lambda|=11$} & \multicolumn{8}{c|}{$|\lambda|=10$} & \multicolumn{6}{c|}{$|\lambda|=9$} & \multicolumn{4}{c|}{$|\lambda|=8$} & \multicolumn{3}{c|}{$|\lambda|=7$} & \multicolumn{1}{c}{$|\lambda|=6$} \\
 &  & \rotatebox{90}{$(10,3)$} & \rotatebox{90}{$(9,4)$} & \rotatebox{90}{$(9,3,1)$} & \rotatebox{90}{$(8,5)$} & \rotatebox{90}{$(8,4,1)$} & \rotatebox{90}{$(8,3,2)$} & \rotatebox{90}{$(7,6)$} & \rotatebox{90}{$(7,5,1)$} & \rotatebox{90}{$(7,4,2)$} & \rotatebox{90}{$(7,3,2,1)$} & \rotatebox{90}{$(6,5,2)$} & \rotatebox{90}{$(6,4,3)$} & \rotatebox{90}{$(6,4,2,1)$} & \rotatebox{90}{$(5,4,3,1)$} & \rotatebox{90}{$(10,2)$} & \rotatebox{90}{$(9,3)$} & \rotatebox{90}{$(9,2,1)$} & \rotatebox{90}{$(8,4)$} & \rotatebox{90}{$(8,3,1)$} & \rotatebox{90}{$(7,5)$} & \rotatebox{90}{$(7,4,1)$} & \rotatebox{90}{$(7,3,2)$} & \rotatebox{90}{$(6,5,1)$} & \rotatebox{90}{$(6,4,2)$} & \rotatebox{90}{$(6,3,2,1)$} & \rotatebox{90}{$(5,4,3)$} & \rotatebox{90}{$(5,4,2,1)$} & \rotatebox{90}{$(9,2)$} & \rotatebox{90}{$(8,3)$} & \rotatebox{90}{$(8,2,1)$} & \rotatebox{90}{$(7,4)$} & \rotatebox{90}{$(7,3,1)$} & \rotatebox{90}{$(6,5)$} & \rotatebox{90}{$(6,4,1)$} & \rotatebox{90}{$(6,3,2)$} & \rotatebox{90}{$(5,4,2)$} & \rotatebox{90}{$(5,3,2,1)$} & \rotatebox{90}{$(8,2)$} & \rotatebox{90}{$(7,3)$} & \rotatebox{90}{$(7,2,1)$} & \rotatebox{90}{$(6,4)$} & \rotatebox{90}{$(6,3,1)$} & \rotatebox{90}{$(5,4,1)$} & \rotatebox{90}{$(5,3,2)$} & \rotatebox{90}{$(4,3,2,1)$} & \rotatebox{90}{$(7,2)$} & \rotatebox{90}{$(6,3)$} & \rotatebox{90}{$(6,2,1)$} & \rotatebox{90}{$(5,4)$} & \rotatebox{90}{$(5,3,1)$} & \rotatebox{90}{$(4,3,2)$} & \rotatebox{90}{$(6,2)$} & \rotatebox{90}{$(5,3)$} & \rotatebox{90}{$(5,2,1)$} & \rotatebox{90}{$(4,3,1)$} & \rotatebox{90}{$(5,2)$} & \rotatebox{90}{$(4,3)$} & \rotatebox{90}{$(4,2,1)$} & \rotatebox{90}{$(4,2)$} \\
\midrule
13 & $(5,4,3,1)$ &  & 8 & 8 & 16 & 24 & 16 & 8 & 24 & 32 & 8 & 16 & 12 & 12 & 4 &  & 12 & 4 & 56 & 60 & 60 & 116 & 76 & 52 & 108 & 24 & 18 & 18 & 4 & 56 & 16 & 126 & 154 & 50 & 176 & 124 & 74 & 20 & 16 & 99 & 25 & 112 & 175 & 62 & 75 & 1 & 25 & 81 & 19 & 21 & 82 & 3 & 19 & 28 & 7 & 3 & 7 & 1 & 1 & 1 \\
13 & $(5,4,4)$ & 16 & 32 & 16 & 16 & 32 & 16 &  & 16 & 32 &  & 16 & 16 &  &  & 8 & 88 & 8 & 120 & 88 & 40 & 120 & 64 & 40 & 96 &  & 16 &  & 36 & 188 & 28 & 172 & 180 & 36 & 164 & 88 & 52 &  & 66 & 202 & 38 & 118 & 174 & 56 & 44 &  & 63 & 116 & 25 & 23 & 72 & 1 & 33 & 32 & 8 & 2 & 9 & 1 & 1 & 1 \\
\bottomrule
\end{tabular}
}%
\usebox{\AFGDTableBox}%
\end{center}
\end{minipage}

\par\medskip
\noindent\begin{minipage}{\linewidth}
\par\medskip\noindent{\large\bfseries Length $\ell(w)=14$, $\lambda(w)=(5,4,4,1)$\par}\nobreak\smallskip
\begin{center}
\fontsize{7.5pt}{9pt}\selectfont
\setlength{\tabcolsep}{1.2pt}
\renewcommand{\arraystretch}{1.05}
\sbox{\AFGDTableBox}{%
\begin{tabular}{rrrrrrrrrrrrrrrr|rrrrrrrrrrrrrrr|rrrrrrrrrrrrr}
\toprule
\multicolumn{16}{c|}{$|\lambda|=14$} & \multicolumn{15}{c|}{$|\lambda|=13$} & \multicolumn{13}{c}{$|\lambda|=12$} \\
\rotatebox{90}{$(10,4)$} & \rotatebox{90}{$(10,3,1)$} & \rotatebox{90}{$(9,5)$} & \rotatebox{90}{$(9,4,1)$} & \rotatebox{90}{$(9,3,2)$} & \rotatebox{90}{$(8,6)$} & \rotatebox{90}{$(8,5,1)$} & \rotatebox{90}{$(8,4,2)$} & \rotatebox{90}{$(8,3,2,1)$} & \rotatebox{90}{$(7,6,1)$} & \rotatebox{90}{$(7,5,2)$} & \rotatebox{90}{$(7,4,3)$} & \rotatebox{90}{$(7,4,2,1)$} & \rotatebox{90}{$(6,5,3)$} & \rotatebox{90}{$(6,5,2,1)$} & \rotatebox{90}{$(6,4,3,1)$} & \rotatebox{90}{$(10,3)$} & \rotatebox{90}{$(10,2,1)$} & \rotatebox{90}{$(9,4)$} & \rotatebox{90}{$(9,3,1)$} & \rotatebox{90}{$(8,5)$} & \rotatebox{90}{$(8,4,1)$} & \rotatebox{90}{$(8,3,2)$} & \rotatebox{90}{$(7,6)$} & \rotatebox{90}{$(7,5,1)$} & \rotatebox{90}{$(7,4,2)$} & \rotatebox{90}{$(7,3,2,1)$} & \rotatebox{90}{$(6,5,2)$} & \rotatebox{90}{$(6,4,3)$} & \rotatebox{90}{$(6,4,2,1)$} & \rotatebox{90}{$(5,4,3,1)$} & \rotatebox{90}{$(10,2)$} & \rotatebox{90}{$(9,3)$} & \rotatebox{90}{$(9,2,1)$} & \rotatebox{90}{$(8,4)$} & \rotatebox{90}{$(8,3,1)$} & \rotatebox{90}{$(7,5)$} & \rotatebox{90}{$(7,4,1)$} & \rotatebox{90}{$(7,3,2)$} & \rotatebox{90}{$(6,5,1)$} & \rotatebox{90}{$(6,4,2)$} & \rotatebox{90}{$(6,3,2,1)$} & \rotatebox{90}{$(5,4,3)$} & \rotatebox{90}{$(5,4,2,1)$} \\
\midrule
16 & 16 & 32 & 48 & 32 & 16 & 48 & 80 & 16 & 16 & 64 & 48 & 32 & 32 & 16 & 16 & 24 & 8 & 120 & 128 & 136 & 272 & 184 & 40 & 192 & 376 & 64 & 168 & 144 & 96 & 16 & 8 & 124 & 36 & 308 & 384 & 212 & 584 & 396 & 204 & 592 & 88 & 68 & 52 \\
\bottomrule
\end{tabular}
}%
\usebox{\AFGDTableBox}%
\end{center}
\end{minipage}

\par\medskip
\noindent\begin{minipage}{\linewidth}
\begin{center}
\fontsize{7.5pt}{9pt}\selectfont
\setlength{\tabcolsep}{1.2pt}
\renewcommand{\arraystretch}{1.05}
\sbox{\AFGDTableBox}{%
\begin{tabular}{rrrrrrrrrr|rrrrrrrr|rrrrrr|rrrr|rrr|r}
\toprule
\multicolumn{10}{c|}{$|\lambda|=11$} & \multicolumn{8}{c|}{$|\lambda|=10$} & \multicolumn{6}{c|}{$|\lambda|=9$} & \multicolumn{4}{c|}{$|\lambda|=8$} & \multicolumn{3}{c|}{$|\lambda|=7$} & \multicolumn{1}{c}{$|\lambda|=6$} \\
\rotatebox{90}{$(9,2)$} & \rotatebox{90}{$(8,3)$} & \rotatebox{90}{$(8,2,1)$} & \rotatebox{90}{$(7,4)$} & \rotatebox{90}{$(7,3,1)$} & \rotatebox{90}{$(6,5)$} & \rotatebox{90}{$(6,4,1)$} & \rotatebox{90}{$(6,3,2)$} & \rotatebox{90}{$(5,4,2)$} & \rotatebox{90}{$(5,3,2,1)$} & \rotatebox{90}{$(8,2)$} & \rotatebox{90}{$(7,3)$} & \rotatebox{90}{$(7,2,1)$} & \rotatebox{90}{$(6,4)$} & \rotatebox{90}{$(6,3,1)$} & \rotatebox{90}{$(5,4,1)$} & \rotatebox{90}{$(5,3,2)$} & \rotatebox{90}{$(4,3,2,1)$} & \rotatebox{90}{$(7,2)$} & \rotatebox{90}{$(6,3)$} & \rotatebox{90}{$(6,2,1)$} & \rotatebox{90}{$(5,4)$} & \rotatebox{90}{$(5,3,1)$} & \rotatebox{90}{$(4,3,2)$} & \rotatebox{90}{$(6,2)$} & \rotatebox{90}{$(5,3)$} & \rotatebox{90}{$(5,2,1)$} & \rotatebox{90}{$(4,3,1)$} & \rotatebox{90}{$(5,2)$} & \rotatebox{90}{$(4,3)$} & \rotatebox{90}{$(4,2,1)$} & \rotatebox{90}{$(4,2)$} \\
\midrule
36 & 254 & 66 & 374 & 572 & 104 & 558 & 394 & 204 & 44 & 66 & 265 & 63 & 218 & 445 & 128 & 161 & 1 & 63 & 147 & 33 & 29 & 154 & 3 & 33 & 38 & 9 & 3 & 9 & 1 & 1 & 1 \\
\bottomrule
\end{tabular}
}%
\usebox{\AFGDTableBox}%
\end{center}
\end{minipage}

\par\medskip
\noindent\begin{minipage}{\linewidth}
\par\medskip\noindent{\large\bfseries Length $\ell(w)=15$, $\lambda(w)=(5,4,4,2)$\par}\nobreak\smallskip
\begin{center}
\fontsize{7.5pt}{9pt}\selectfont
\setlength{\tabcolsep}{0.7pt}
\renewcommand{\arraystretch}{1.05}
\sbox{\AFGDTableBox}{%
\begin{tabular}{rrrrrrrrrrrrrrrrrrrr|rrrrrrrrrrrrrrrrrrr|rrrrrrrrrrrrrrrr}
\toprule
\multicolumn{20}{c|}{$|\lambda|=15$} & \multicolumn{19}{c|}{$|\lambda|=14$} & \multicolumn{16}{c}{$|\lambda|=13$} \\
\rotatebox{90}{$(11,4)$} & \rotatebox{90}{$(11,3,1)$} & \rotatebox{90}{$(10,5)$} & \rotatebox{90}{$(10,4,1)$} & \rotatebox{90}{$(10,3,2)$} & \rotatebox{90}{$(9,6)$} & \rotatebox{90}{$(9,5,1)$} & \rotatebox{90}{$(9,4,2)$} & \rotatebox{90}{$(9,3,2,1)$} & \rotatebox{90}{$(8,6,1)$} & \rotatebox{90}{$(8,5,2)$} & \rotatebox{90}{$(8,4,3)$} & \rotatebox{90}{$(8,4,2,1)$} & \rotatebox{90}{$(7,6,2)$} & \rotatebox{90}{$(7,5,3)$} & \rotatebox{90}{$(7,5,2,1)$} & \rotatebox{90}{$(7,4,3,1)$} & \rotatebox{90}{$(6,5,4)$} & \rotatebox{90}{$(6,5,3,1)$} & \rotatebox{90}{$(6,4,3,2)$} & \rotatebox{90}{$(11,3)$} & \rotatebox{90}{$(11,2,1)$} & \rotatebox{90}{$(10,4)$} & \rotatebox{90}{$(10,3,1)$} & \rotatebox{90}{$(9,5)$} & \rotatebox{90}{$(9,4,1)$} & \rotatebox{90}{$(9,3,2)$} & \rotatebox{90}{$(8,6)$} & \rotatebox{90}{$(8,5,1)$} & \rotatebox{90}{$(8,4,2)$} & \rotatebox{90}{$(8,3,2,1)$} & \rotatebox{90}{$(7,6,1)$} & \rotatebox{90}{$(7,5,2)$} & \rotatebox{90}{$(7,4,3)$} & \rotatebox{90}{$(7,4,2,1)$} & \rotatebox{90}{$(6,5,3)$} & \rotatebox{90}{$(6,5,2,1)$} & \rotatebox{90}{$(6,4,3,1)$} & \rotatebox{90}{$(5,4,3,2)$} & \rotatebox{90}{$(11,2)$} & \rotatebox{90}{$(10,3)$} & \rotatebox{90}{$(10,2,1)$} & \rotatebox{90}{$(9,4)$} & \rotatebox{90}{$(9,3,1)$} & \rotatebox{90}{$(8,5)$} & \rotatebox{90}{$(8,4,1)$} & \rotatebox{90}{$(8,3,2)$} & \rotatebox{90}{$(7,6)$} & \rotatebox{90}{$(7,5,1)$} & \rotatebox{90}{$(7,4,2)$} & \rotatebox{90}{$(7,3,2,1)$} & \rotatebox{90}{$(6,5,2)$} & \rotatebox{90}{$(6,4,3)$} & \rotatebox{90}{$(6,4,2,1)$} & \rotatebox{90}{$(5,4,3,1)$} \\
\midrule
16 & 16 & 32 & 64 & 48 & 16 & 80 & 144 & 32 & 32 & 144 & 112 & 80 & 48 & 112 & 64 & 64 & 24 & 40 & 16 & 24 & 8 & 136 & 168 & 168 & 456 & 336 & 56 & 408 & 864 & 160 & 112 & 648 & 488 & 352 & 276 & 156 & 188 & 16 & 8 & 148 & 44 & 428 & 644 & 348 & 1268 & 920 & 76 & 804 & 1988 & 304 & 772 & 698 & 518 & 90 \\
\bottomrule
\end{tabular}
}%
\usebox{\AFGDTableBox}%
\end{center}
\end{minipage}

\par\medskip
\noindent\begin{minipage}{\linewidth}
\begin{center}
\fontsize{7.5pt}{9pt}\selectfont
\setlength{\tabcolsep}{1.2pt}
\renewcommand{\arraystretch}{1.05}
\sbox{\AFGDTableBox}{%
\begin{tabular}{rrrrrrrrrrrrr|rrrrrrrrrr|rrrrrrrr|rrrrrr|rrrr|rrr|r}
\toprule
\multicolumn{13}{c|}{$|\lambda|=12$} & \multicolumn{10}{c|}{$|\lambda|=11$} & \multicolumn{8}{c|}{$|\lambda|=10$} & \multicolumn{6}{c|}{$|\lambda|=9$} & \multicolumn{4}{c|}{$|\lambda|=8$} & \multicolumn{3}{c|}{$|\lambda|=7$} & \multicolumn{1}{c}{$|\lambda|=6$} \\
\rotatebox{90}{$(10,2)$} & \rotatebox{90}{$(9,3)$} & \rotatebox{90}{$(9,2,1)$} & \rotatebox{90}{$(8,4)$} & \rotatebox{90}{$(8,3,1)$} & \rotatebox{90}{$(7,5)$} & \rotatebox{90}{$(7,4,1)$} & \rotatebox{90}{$(7,3,2)$} & \rotatebox{90}{$(6,5,1)$} & \rotatebox{90}{$(6,4,2)$} & \rotatebox{90}{$(6,3,2,1)$} & \rotatebox{90}{$(5,4,3)$} & \rotatebox{90}{$(5,4,2,1)$} & \rotatebox{90}{$(9,2)$} & \rotatebox{90}{$(8,3)$} & \rotatebox{90}{$(8,2,1)$} & \rotatebox{90}{$(7,4)$} & \rotatebox{90}{$(7,3,1)$} & \rotatebox{90}{$(6,5)$} & \rotatebox{90}{$(6,4,1)$} & \rotatebox{90}{$(6,3,2)$} & \rotatebox{90}{$(5,4,2)$} & \rotatebox{90}{$(5,3,2,1)$} & \rotatebox{90}{$(8,2)$} & \rotatebox{90}{$(7,3)$} & \rotatebox{90}{$(7,2,1)$} & \rotatebox{90}{$(6,4)$} & \rotatebox{90}{$(6,3,1)$} & \rotatebox{90}{$(5,4,1)$} & \rotatebox{90}{$(5,3,2)$} & \rotatebox{90}{$(4,3,2,1)$} & \rotatebox{90}{$(7,2)$} & \rotatebox{90}{$(6,3)$} & \rotatebox{90}{$(6,2,1)$} & \rotatebox{90}{$(5,4)$} & \rotatebox{90}{$(5,3,1)$} & \rotatebox{90}{$(4,3,2)$} & \rotatebox{90}{$(6,2)$} & \rotatebox{90}{$(5,3)$} & \rotatebox{90}{$(5,2,1)$} & \rotatebox{90}{$(4,3,1)$} & \rotatebox{90}{$(5,2)$} & \rotatebox{90}{$(4,3)$} & \rotatebox{90}{$(4,2,1)$} & \rotatebox{90}{$(4,2)$} \\
\midrule
44 & 378 & 102 & 682 & 1254 & 366 & 1790 & 1276 & 542 & 2026 & 270 & 191 & 163 & 102 & 519 & 129 & 608 & 1378 & 141 & 1267 & 931 & 461 & 99 & 129 & 414 & 96 & 289 & 850 & 214 & 300 & 1 & 96 & 190 & 42 & 33 & 244 & 3 & 42 & 43 & 10 & 3 & 10 & 1 & 1 & 1 \\
\bottomrule
\end{tabular}
}%
\usebox{\AFGDTableBox}%
\end{center}
\end{minipage}

\par\medskip
\noindent\begin{minipage}{\linewidth}
\par\medskip\noindent{\large\bfseries Length $\ell(w)=16$, $\lambda(w)=(5,4,4,2,1)$\par}\nobreak\smallskip
\begin{center}
\fontsize{7.5pt}{9pt}\selectfont
\setlength{\tabcolsep}{1.2pt}
\renewcommand{\arraystretch}{1.05}
\sbox{\AFGDTableBox}{%
\begin{tabular}{rrrrrrrrrrrrrrrrrrrrrrrrrrr}
\toprule
\multicolumn{27}{c}{$|\lambda|=16$} \\
\rotatebox{90}{$(12,4)$} & \rotatebox{90}{$(12,3,1)$} & \rotatebox{90}{$(11,5)$} & \rotatebox{90}{$(11,4,1)$} & \rotatebox{90}{$(11,3,2)$} & \rotatebox{90}{$(10,6)$} & \rotatebox{90}{$(10,5,1)$} & \rotatebox{90}{$(10,4,2)$} & \rotatebox{90}{$(10,3,2,1)$} & \rotatebox{90}{$(9,7)$} & \rotatebox{90}{$(9,6,1)$} & \rotatebox{90}{$(9,5,2)$} & \rotatebox{90}{$(9,4,3)$} & \rotatebox{90}{$(9,4,2,1)$} & \rotatebox{90}{$(8,7,1)$} & \rotatebox{90}{$(8,6,2)$} & \rotatebox{90}{$(8,5,3)$} & \rotatebox{90}{$(8,5,2,1)$} & \rotatebox{90}{$(8,4,3,1)$} & \rotatebox{90}{$(7,6,3)$} & \rotatebox{90}{$(7,6,2,1)$} & \rotatebox{90}{$(7,5,4)$} & \rotatebox{90}{$(7,5,3,1)$} & \rotatebox{90}{$(7,4,3,2)$} & \rotatebox{90}{$(6,5,4,1)$} & \rotatebox{90}{$(6,5,3,2)$} & \rotatebox{90}{$(6,4,3,2,1)$} \\
\midrule
16 & 16 & 64 & 80 & 48 & 80 & 176 & 208 & 48 & 32 & 176 & 368 & 208 & 160 & 64 & 304 & 368 & 208 & 176 & 192 & 96 & 144 & 208 & 80 & 48 & 56 & 16 \\
\bottomrule
\end{tabular}
}%
\usebox{\AFGDTableBox}%

\par\smallskip
\sbox{\AFGDTableBox}{%
\begin{tabular}{rrrrrrrrrrrrrrrrrrrrrrrr}
\toprule
\multicolumn{24}{c}{$|\lambda|=15$} \\
\rotatebox{90}{$(12,3)$} & \rotatebox{90}{$(12,2,1)$} & \rotatebox{90}{$(11,4)$} & \rotatebox{90}{$(11,3,1)$} & \rotatebox{90}{$(10,5)$} & \rotatebox{90}{$(10,4,1)$} & \rotatebox{90}{$(10,3,2)$} & \rotatebox{90}{$(9,6)$} & \rotatebox{90}{$(9,5,1)$} & \rotatebox{90}{$(9,4,2)$} & \rotatebox{90}{$(9,3,2,1)$} & \rotatebox{90}{$(8,7)$} & \rotatebox{90}{$(8,6,1)$} & \rotatebox{90}{$(8,5,2)$} & \rotatebox{90}{$(8,4,3)$} & \rotatebox{90}{$(8,4,2,1)$} & \rotatebox{90}{$(7,6,2)$} & \rotatebox{90}{$(7,5,3)$} & \rotatebox{90}{$(7,5,2,1)$} & \rotatebox{90}{$(7,4,3,1)$} & \rotatebox{90}{$(6,5,4)$} & \rotatebox{90}{$(6,5,3,1)$} & \rotatebox{90}{$(6,4,3,2)$} & \rotatebox{90}{$(5,4,3,2,1)$} \\
\midrule
24 & 8 & 200 & 192 & 472 & 752 & 448 & 408 & 1248 & 1728 & 304 & 112 & 904 & 2392 & 1336 & 976 & 1040 & 1720 & 920 & 824 & 288 & 520 & 244 & 16 \\
\bottomrule
\end{tabular}
}%
\usebox{\AFGDTableBox}%

\par\smallskip
\sbox{\AFGDTableBox}{%
\begin{tabular}{rrrrrrrrrrrrrrrrrrrr}
\toprule
\multicolumn{20}{c}{$|\lambda|=14$} \\
\rotatebox{90}{$(12,2)$} & \rotatebox{90}{$(11,3)$} & \rotatebox{90}{$(11,2,1)$} & \rotatebox{90}{$(10,4)$} & \rotatebox{90}{$(10,3,1)$} & \rotatebox{90}{$(9,5)$} & \rotatebox{90}{$(9,4,1)$} & \rotatebox{90}{$(9,3,2)$} & \rotatebox{90}{$(8,6)$} & \rotatebox{90}{$(8,5,1)$} & \rotatebox{90}{$(8,4,2)$} & \rotatebox{90}{$(8,3,2,1)$} & \rotatebox{90}{$(7,6,1)$} & \rotatebox{90}{$(7,5,2)$} & \rotatebox{90}{$(7,4,3)$} & \rotatebox{90}{$(7,4,2,1)$} & \rotatebox{90}{$(6,5,3)$} & \rotatebox{90}{$(6,5,2,1)$} & \rotatebox{90}{$(6,4,3,1)$} & \rotatebox{90}{$(5,4,3,2)$} \\
\midrule
8 & 188 & 52 & 860 & 920 & 1372 & 2832 & 1656 & 828 & 3520 & 5684 & 776 & 1292 & 5404 & 3112 & 2244 & 1812 & 972 & 1224 & 114 \\
\bottomrule
\end{tabular}
}%
\usebox{\AFGDTableBox}%
\end{center}
\end{minipage}

\par\medskip
\noindent\begin{minipage}{\linewidth}
\begin{center}
\fontsize{7.5pt}{9pt}\selectfont
\setlength{\tabcolsep}{1.2pt}
\renewcommand{\arraystretch}{1.05}
\sbox{\AFGDTableBox}{%
\begin{tabular}{rrrrrrrrrrrrrrrr|rrrrrrrrrrrrr}
\toprule
\multicolumn{16}{c|}{$|\lambda|=13$} & \multicolumn{13}{c}{$|\lambda|=12$} \\
\rotatebox{90}{$(11,2)$} & \rotatebox{90}{$(10,3)$} & \rotatebox{90}{$(10,2,1)$} & \rotatebox{90}{$(9,4)$} & \rotatebox{90}{$(9,3,1)$} & \rotatebox{90}{$(8,5)$} & \rotatebox{90}{$(8,4,1)$} & \rotatebox{90}{$(8,3,2)$} & \rotatebox{90}{$(7,6)$} & \rotatebox{90}{$(7,5,1)$} & \rotatebox{90}{$(7,4,2)$} & \rotatebox{90}{$(7,3,2,1)$} & \rotatebox{90}{$(6,5,2)$} & \rotatebox{90}{$(6,4,3)$} & \rotatebox{90}{$(6,4,2,1)$} & \rotatebox{90}{$(5,4,3,1)$} & \rotatebox{90}{$(10,2)$} & \rotatebox{90}{$(9,3)$} & \rotatebox{90}{$(9,2,1)$} & \rotatebox{90}{$(8,4)$} & \rotatebox{90}{$(8,3,1)$} & \rotatebox{90}{$(7,5)$} & \rotatebox{90}{$(7,4,1)$} & \rotatebox{90}{$(7,3,2)$} & \rotatebox{90}{$(6,5,1)$} & \rotatebox{90}{$(6,4,2)$} & \rotatebox{90}{$(6,3,2,1)$} & \rotatebox{90}{$(5,4,3)$} & \rotatebox{90}{$(5,4,2,1)$} \\
\midrule
52 & 614 & 146 & 1866 & 2352 & 2078 & 5624 & 3204 & 606 & 4644 & 9208 & 1012 & 3574 & 2950 & 2188 & 324 & 146 & 1101 & 231 & 2328 & 3561 & 1657 & 6223 & 3503 & 1881 & 6948 & 687 & 469 & 415 \\
\bottomrule
\end{tabular}
}%
\usebox{\AFGDTableBox}%

\par\smallskip
\sbox{\AFGDTableBox}{%
\begin{tabular}{rrrrrrrrrr|rrrrrrrr|rrrrrr|rrrr|rrr|r}
\toprule
\multicolumn{10}{c|}{$|\lambda|=11$} & \multicolumn{8}{c|}{$|\lambda|=10$} & \multicolumn{6}{c|}{$|\lambda|=9$} & \multicolumn{4}{c|}{$|\lambda|=8$} & \multicolumn{3}{c|}{$|\lambda|=7$} & \multicolumn{1}{c}{$|\lambda|=6$} \\
\rotatebox{90}{$(9,2)$} & \rotatebox{90}{$(8,3)$} & \rotatebox{90}{$(8,2,1)$} & \rotatebox{90}{$(7,4)$} & \rotatebox{90}{$(7,3,1)$} & \rotatebox{90}{$(6,5)$} & \rotatebox{90}{$(6,4,1)$} & \rotatebox{90}{$(6,3,2)$} & \rotatebox{90}{$(5,4,2)$} & \rotatebox{90}{$(5,3,2,1)$} & \rotatebox{90}{$(8,2)$} & \rotatebox{90}{$(7,3)$} & \rotatebox{90}{$(7,2,1)$} & \rotatebox{90}{$(6,4)$} & \rotatebox{90}{$(6,3,1)$} & \rotatebox{90}{$(5,4,1)$} & \rotatebox{90}{$(5,3,2)$} & \rotatebox{90}{$(4,3,2,1)$} & \rotatebox{90}{$(7,2)$} & \rotatebox{90}{$(6,3)$} & \rotatebox{90}{$(6,2,1)$} & \rotatebox{90}{$(5,4)$} & \rotatebox{90}{$(5,3,1)$} & \rotatebox{90}{$(4,3,2)$} & \rotatebox{90}{$(6,2)$} & \rotatebox{90}{$(5,3)$} & \rotatebox{90}{$(5,2,1)$} & \rotatebox{90}{$(4,3,1)$} & \rotatebox{90}{$(5,2)$} & \rotatebox{90}{$(4,3)$} & \rotatebox{90}{$(4,2,1)$} & \rotatebox{90}{$(4,2)$} \\
\midrule
231 & 1191 & 225 & 1707 & 3268 & 395 & 3451 & 2101 & 995 & 200 & 225 & 794 & 138 & 655 & 1724 & 386 & 559 & 1 & 138 & 312 & 52 & 49 & 419 & 3 & 52 & 60 & 11 & 3 & 11 & 1 & 1 & 1 \\
\bottomrule
\end{tabular}
}%
\usebox{\AFGDTableBox}%
\end{center}
\end{minipage}
\end{landscape}

\newgeometry{left=20mm,right=20mm,top=25mm,bottom=25mm}
\subsection{Sign-normalized coefficient tables for type $\mathrm{B}_3$}
The following table records the positive integers $b_{w,\lambda}$; a blank entry means
$b_{w,\lambda}=0$. For $\ell(w)\leq 12$, the rows are labeled by $\lambda(w)\in \Y^3_{\B{}}$; at length $13$, the Weyl
group element is displayed in the heading and the table contains only its
coefficient row. The columns are ordered by decreasing $|\lambda|$ and,
for fixed size, by decreasing lexicographic order. Horizontal rules separate
distinct values of $\ell(w)$, and vertical rules separate distinct values of
$|\lambda|$. Columns that vanish throughout a displayed table are omitted.
\par\medskip
\noindent\begin{minipage}{\linewidth}
\par\medskip\noindent{\large\bfseries Lengths $0\leq\ell(w)\leq10$\par}\nobreak\smallskip
\begin{center}
\fontsize{7.5pt}{9pt}\selectfont
\setlength{\tabcolsep}{0.6pt}
\renewcommand{\arraystretch}{1.05}
\sbox{\AFGQTableBox}{%
\begin{tabular}{r|l|rrrrrrrr|rrrrrrr|rrrrr|rrrr|rrr|rrr|rr|rr|r|r|r}
\toprule
$\ell(w)$ & $\lambda(w)$ & \multicolumn{8}{c|}{$|\lambda|=10$} & \multicolumn{7}{c|}{$|\lambda|=9$} & \multicolumn{5}{c|}{$|\lambda|=8$} & \multicolumn{4}{c|}{$|\lambda|=7$} & \multicolumn{3}{c|}{$|\lambda|=6$} & \multicolumn{3}{c|}{$|\lambda|=5$} & \multicolumn{2}{c|}{$|\lambda|=4$} & \multicolumn{2}{c|}{$|\lambda|=3$} & \multicolumn{1}{c|}{$|\lambda|=2$} & \multicolumn{1}{c|}{$|\lambda|=1$} & \multicolumn{1}{c}{$|\lambda|=0$} \\
 &  & \rotatebox{90}{$(8,2)$} & \rotatebox{90}{$(7,3)$} & \rotatebox{90}{$(7,2,1)$} & \rotatebox{90}{$(6,4)$} & \rotatebox{90}{$(6,3,1)$} & \rotatebox{90}{$(5,4,1)$} & \rotatebox{90}{$(5,3,2)$} & \rotatebox{90}{$(4,3,2,1)$} & \rotatebox{90}{$(8,1)$} & \rotatebox{90}{$(7,2)$} & \rotatebox{90}{$(6,3)$} & \rotatebox{90}{$(6,2,1)$} & \rotatebox{90}{$(5,4)$} & \rotatebox{90}{$(5,3,1)$} & \rotatebox{90}{$(4,3,2)$} & \rotatebox{90}{$(7,1)$} & \rotatebox{90}{$(6,2)$} & \rotatebox{90}{$(5,3)$} & \rotatebox{90}{$(5,2,1)$} & \rotatebox{90}{$(4,3,1)$} & \rotatebox{90}{$(6,1)$} & \rotatebox{90}{$(5,2)$} & \rotatebox{90}{$(4,3)$} & \rotatebox{90}{$(4,2,1)$} & \rotatebox{90}{$(5,1)$} & \rotatebox{90}{$(4,2)$} & \rotatebox{90}{$(3,2,1)$} & \rotatebox{90}{$(5)$} & \rotatebox{90}{$(4,1)$} & \rotatebox{90}{$(3,2)$} & \rotatebox{90}{$(4)$} & \rotatebox{90}{$(3,1)$} & \rotatebox{90}{$(3)$} & \rotatebox{90}{$(2,1)$} & \rotatebox{90}{$(2)$} & \rotatebox{90}{$(1)$} & \rotatebox{90}{$\varnothing$} \\
\midrule
0 & $e$ &  &  &  &  &  &  &  &  &  &  &  &  &  &  &  &  &  &  &  &  &  &  &  &  &  &  &  &  &  &  &  &  &  &  &  &  & 1 \\
\midrule
1 & $(1)$ &  &  &  &  &  &  &  &  &  &  &  &  &  &  &  &  &  &  &  &  &  &  &  &  &  &  &  &  &  &  &  &  &  &  &  & 1 &  \\
\midrule
2 & $(2)$ &  &  &  &  &  &  &  &  &  &  &  &  &  &  &  &  &  &  &  &  &  &  &  &  &  &  &  &  &  &  &  &  &  &  & 1 &  &  \\
\midrule
3 & $(2,1)$ &  &  &  &  &  &  &  &  &  &  &  &  &  &  &  &  &  &  &  &  &  &  &  &  &  &  &  &  &  &  &  &  &  & 1 &  &  &  \\
3 & $(3)$ &  &  &  &  &  &  &  &  &  &  &  &  &  &  &  &  &  &  &  &  &  &  &  &  &  &  &  &  &  &  &  &  & 2 &  & 1 &  &  \\
\midrule
4 & $(3,1)$ &  &  &  &  &  &  &  &  &  &  &  &  &  &  &  &  &  &  &  &  &  &  &  &  &  &  &  &  &  &  &  & 2 &  & 1 &  &  &  \\
4 & $(4)$ &  &  &  &  &  &  &  &  &  &  &  &  &  &  &  &  &  &  &  &  &  &  &  &  &  &  &  &  &  &  & 2 &  & 1 &  &  &  &  \\
\midrule
5 & $(4,1)$ &  &  &  &  &  &  &  &  &  &  &  &  &  &  &  &  &  &  &  &  &  &  &  &  &  &  &  & 2 & 2 &  & 3 & 1 & 1 &  &  &  &  \\
5 & $(3,2)$ &  &  &  &  &  &  &  &  &  &  &  &  &  &  &  &  &  &  &  &  &  &  &  &  &  &  &  &  & 2 & 2 &  & 3 &  &  &  &  &  \\
\midrule
6 & $(3,2,1)$ &  &  &  &  &  &  &  &  &  &  &  &  &  &  &  &  &  &  &  &  &  &  &  &  & 2 & 2 & 2 &  & 3 & 3 &  & 1 &  &  &  &  &  \\
6 & $(4,2)$ &  &  &  &  &  &  &  &  &  &  &  &  &  &  &  &  &  &  &  &  &  &  &  &  & 2 & 2 &  &  & 3 & 1 &  & 1 &  &  &  &  &  \\
6 & $(3,3)$ &  &  &  &  &  &  &  &  &  &  &  &  &  &  &  &  &  &  &  &  &  &  &  &  &  & 4 &  &  & 2 & 2 &  & 1 &  &  &  &  &  \\
\midrule
7 & $(4,2,1)$ &  &  &  &  &  &  &  &  &  &  &  &  &  &  &  &  &  &  &  &  & 4 & 4 &  & 2 & 8 & 6 & 1 &  & 5 & 2 &  & 1 &  &  &  &  &  \\
7 & $(3,3,1)$ &  &  &  &  &  &  &  &  &  &  &  &  &  &  &  &  &  &  &  &  &  & 4 & 4 & 4 & 2 & 12 & 2 &  & 3 & 3 &  & 1 &  &  &  &  &  \\
7 & $(4,3)$ &  &  &  &  &  &  &  &  &  &  &  &  &  &  &  &  &  &  &  &  &  & 4 & 4 &  & 2 & 8 &  &  & 3 & 1 &  & 1 &  &  &  &  &  \\
\midrule
8 & $(4,3,1)$ &  &  &  &  &  &  &  &  &  &  &  &  &  &  &  &  & 8 & 8 & 4 & 4 & 4 & 24 & 12 & 8 & 8 & 22 & 1 &  & 5 & 2 &  & 1 &  &  &  &  &  \\
8 & $(3,3,2)$ &  &  &  &  &  &  &  &  &  &  &  &  &  &  &  &  &  & 4 & 4 & 4 &  & 6 & 6 & 6 &  & 7 &  &  &  &  &  &  &  &  &  &  &  \\
\midrule
9 & $(3,3,2,1)$ &  &  &  &  &  &  &  &  &  &  & 8 & 8 & 8 & 12 & 8 &  & 12 & 32 & 22 & 22 &  & 30 & 18 & 17 &  & 20 &  &  &  &  &  &  &  &  &  &  &  \\
9 & $(4,3,2)$ &  &  &  &  &  &  &  &  &  & 8 & 16 & 8 & 8 & 12 & 4 & 4 & 40 & 48 & 24 & 18 & 12 & 66 & 28 & 20 & 13 & 40 & 1 &  & 6 & 2 &  & 1 &  &  &  &  &  \\
\midrule
10 & $(4,3,2,1)$ &  & 8 & 8 & 16 & 16 & 8 & 12 & 4 &  & 12 & 48 & 28 & 28 & 48 & 18 &  & 40 & 68 & 34 & 28 &  & 45 & 17 & 15 &  & 18 &  &  &  &  &  &  &  &  &  &  &  \\
10 & $(4,3,3)$ & 16 & 32 & 16 & 16 & 32 & 16 & 16 &  & 8 & 88 & 120 & 64 & 40 & 96 & 16 & 28 & 180 & 164 & 88 & 52 & 38 & 174 & 56 & 44 & 25 & 72 & 1 &  & 8 & 2 &  & 1 &  &  &  &  &  \\
\bottomrule
\end{tabular}
}%
\usebox{\AFGQTableBox}%
\end{center}
\end{minipage}
\newgeometry{top=15mm,bottom=15mm,left=18mm,right=18mm}
\begin{landscape}

\par\medskip
\noindent\begin{minipage}{\linewidth}
\par\medskip\noindent{\large\bfseries Lengths $11\leq\ell(w)\leq12$\par}\nobreak\smallskip
\begin{center}
\fontsize{7.5pt}{9pt}\selectfont
\setlength{\tabcolsep}{1.0pt}
\renewcommand{\arraystretch}{1.05}
\sbox{\AFGQTableBox}{%
\begin{tabular}{r|l|rrrrrrrrrrrr|rrrrrrrrrr|rrrrrrrr|rrrrrr|rrrr|rrr|r}
\toprule
$\ell(w)$ & $w$ & \multicolumn{12}{c|}{$|\lambda|=12$} & \multicolumn{10}{c|}{$|\lambda|=11$} & \multicolumn{8}{c|}{$|\lambda|=10$} & \multicolumn{6}{c|}{$|\lambda|=9$} & \multicolumn{4}{c|}{$|\lambda|=8$} & \multicolumn{3}{c|}{$|\lambda|=7$} & \multicolumn{1}{c}{$|\lambda|=6$} \\
 &  & \rotatebox{90}{$(9,3)$} & \rotatebox{90}{$(9,2,1)$} & \rotatebox{90}{$(8,4)$} & \rotatebox{90}{$(8,3,1)$} & \rotatebox{90}{$(7,5)$} & \rotatebox{90}{$(7,4,1)$} & \rotatebox{90}{$(7,3,2)$} & \rotatebox{90}{$(6,5,1)$} & \rotatebox{90}{$(6,4,2)$} & \rotatebox{90}{$(6,3,2,1)$} & \rotatebox{90}{$(5,4,3)$} & \rotatebox{90}{$(5,4,2,1)$} & \rotatebox{90}{$(9,2)$} & \rotatebox{90}{$(8,3)$} & \rotatebox{90}{$(8,2,1)$} & \rotatebox{90}{$(7,4)$} & \rotatebox{90}{$(7,3,1)$} & \rotatebox{90}{$(6,5)$} & \rotatebox{90}{$(6,4,1)$} & \rotatebox{90}{$(6,3,2)$} & \rotatebox{90}{$(5,4,2)$} & \rotatebox{90}{$(5,3,2,1)$} & \rotatebox{90}{$(8,2)$} & \rotatebox{90}{$(7,3)$} & \rotatebox{90}{$(7,2,1)$} & \rotatebox{90}{$(6,4)$} & \rotatebox{90}{$(6,3,1)$} & \rotatebox{90}{$(5,4,1)$} & \rotatebox{90}{$(5,3,2)$} & \rotatebox{90}{$(4,3,2,1)$} & \rotatebox{90}{$(7,2)$} & \rotatebox{90}{$(6,3)$} & \rotatebox{90}{$(6,2,1)$} & \rotatebox{90}{$(5,4)$} & \rotatebox{90}{$(5,3,1)$} & \rotatebox{90}{$(4,3,2)$} & \rotatebox{90}{$(6,2)$} & \rotatebox{90}{$(5,3)$} & \rotatebox{90}{$(5,2,1)$} & \rotatebox{90}{$(4,3,1)$} & \rotatebox{90}{$(5,2)$} & \rotatebox{90}{$(4,3)$} & \rotatebox{90}{$(4,2,1)$} & \rotatebox{90}{$(4,2)$} \\
\midrule
11 & $(4,3,3,1)$ &  &  &  &  &  &  &  &  &  &  &  &  &  & 16 & 16 & 32 & 48 & 16 & 48 & 48 & 32 & 16 & 24 & 120 & 80 & 136 & 224 & 112 & 144 & 16 & 108 & 292 & 152 & 128 & 332 & 68 & 190 & 276 & 132 & 96 & 155 & 49 & 45 & 50 \\
\midrule
12 & $(4,3,3,2)$ & 16 & 16 & 32 & 64 & 16 & 80 & 80 & 32 & 96 & 32 & 24 & 24 & 24 & 136 & 96 & 168 & 360 & 56 & 336 & 336 & 220 & 76 & 124 & 396 & 224 & 324 & 768 & 340 & 438 & 22 & 262 & 550 & 264 & 188 & 690 & 83 & 289 & 356 & 161 & 102 & 169 & 42 & 41 & 42 \\
\bottomrule
\end{tabular}
}%
\usebox{\AFGQTableBox}%
\end{center}
\end{minipage}

\par\medskip
\noindent\begin{minipage}{\linewidth}
\par\medskip\noindent{\large\bfseries Length $\ell(w)=13$, $\lambda(w)=(4,3,3,2,1)$\par}\nobreak\smallskip
\begin{center}
\fontsize{7.5pt}{9pt}\selectfont
\setlength{\tabcolsep}{2.0pt}
\renewcommand{\arraystretch}{1.05}
\sbox{\AFGQTableBox}{%
\begin{tabular}{rrrrrrrrrrrrrrr|rrrrrrrrrrrrr}
\toprule
\multicolumn{15}{c|}{$|\lambda|=13$} & \multicolumn{13}{c}{$|\lambda|=12$} \\
\rotatebox{90}{$(10,3)$} & \rotatebox{90}{$(10,2,1)$} & \rotatebox{90}{$(9,4)$} & \rotatebox{90}{$(9,3,1)$} & \rotatebox{90}{$(8,5)$} & \rotatebox{90}{$(8,4,1)$} & \rotatebox{90}{$(8,3,2)$} & \rotatebox{90}{$(7,6)$} & \rotatebox{90}{$(7,5,1)$} & \rotatebox{90}{$(7,4,2)$} & \rotatebox{90}{$(7,3,2,1)$} & \rotatebox{90}{$(6,5,2)$} & \rotatebox{90}{$(6,4,3)$} & \rotatebox{90}{$(6,4,2,1)$} & \rotatebox{90}{$(5,4,3,1)$} & \rotatebox{90}{$(10,2)$} & \rotatebox{90}{$(9,3)$} & \rotatebox{90}{$(9,2,1)$} & \rotatebox{90}{$(8,4)$} & \rotatebox{90}{$(8,3,1)$} & \rotatebox{90}{$(7,5)$} & \rotatebox{90}{$(7,4,1)$} & \rotatebox{90}{$(7,3,2)$} & \rotatebox{90}{$(6,5,1)$} & \rotatebox{90}{$(6,4,2)$} & \rotatebox{90}{$(6,3,2,1)$} & \rotatebox{90}{$(5,4,3)$} & \rotatebox{90}{$(5,4,2,1)$} \\
\midrule
16 & 16 & 64 & 80 & 80 & 176 & 144 & 32 & 176 & 288 & 80 & 160 & 144 & 112 & 48 & 24 & 200 & 120 & 472 & 656 & 408 & 1128 & 880 & 544 & 1312 & 304 & 288 & 244 \\
\bottomrule
\end{tabular}
}%
\usebox{\AFGQTableBox}%
\end{center}
\end{minipage}

\par\medskip
\noindent\begin{minipage}{\linewidth}
\begin{center}
\fontsize{7.5pt}{9pt}\selectfont
\setlength{\tabcolsep}{2.0pt}
\renewcommand{\arraystretch}{1.05}
\sbox{\AFGQTableBox}{%
\begin{tabular}{rrrrrrrrrr|rrrrrrrr|rrrrrr|rrrr|rrr|r}
\toprule
\multicolumn{10}{c|}{$|\lambda|=11$} & \multicolumn{8}{c|}{$|\lambda|=10$} & \multicolumn{6}{c|}{$|\lambda|=9$} & \multicolumn{4}{c|}{$|\lambda|=8$} & \multicolumn{3}{c|}{$|\lambda|=7$} & \multicolumn{1}{c}{$|\lambda|=6$} \\
\rotatebox{90}{$(9,2)$} & \rotatebox{90}{$(8,3)$} & \rotatebox{90}{$(8,2,1)$} & \rotatebox{90}{$(7,4)$} & \rotatebox{90}{$(7,3,1)$} & \rotatebox{90}{$(6,5)$} & \rotatebox{90}{$(6,4,1)$} & \rotatebox{90}{$(6,3,2)$} & \rotatebox{90}{$(5,4,2)$} & \rotatebox{90}{$(5,3,2,1)$} & \rotatebox{90}{$(8,2)$} & \rotatebox{90}{$(7,3)$} & \rotatebox{90}{$(7,2,1)$} & \rotatebox{90}{$(6,4)$} & \rotatebox{90}{$(6,3,1)$} & \rotatebox{90}{$(5,4,1)$} & \rotatebox{90}{$(5,3,2)$} & \rotatebox{90}{$(4,3,2,1)$} & \rotatebox{90}{$(7,2)$} & \rotatebox{90}{$(6,3)$} & \rotatebox{90}{$(6,2,1)$} & \rotatebox{90}{$(5,4)$} & \rotatebox{90}{$(5,3,1)$} & \rotatebox{90}{$(4,3,2)$} & \rotatebox{90}{$(6,2)$} & \rotatebox{90}{$(5,3)$} & \rotatebox{90}{$(5,2,1)$} & \rotatebox{90}{$(4,3,1)$} & \rotatebox{90}{$(5,2)$} & \rotatebox{90}{$(4,3)$} & \rotatebox{90}{$(4,2,1)$} & \rotatebox{90}{$(4,2)$} \\
\midrule
164 & 828 & 380 & 1332 & 2084 & 540 & 2372 & 1908 & 1296 & 354 & 474 & 1678 & 646 & 1662 & 3148 & 1322 & 1616 & 45 & 747 & 1792 & 620 & 563 & 2090 & 156 & 678 & 898 & 319 & 180 & 336 & 71 & 69 & 71 \\
\bottomrule
\end{tabular}
}%
\usebox{\AFGQTableBox}%
\end{center}
\end{minipage}
\end{landscape}
\restoregeometry

\end{document}